%% file: geometry-higher-rank.tex
\documentclass[11pt, a4paper]{amsart}
\usepackage{amscd,amsmath,amsthm,amssymb}
\usepackage{mathtools}
\usepackage{mathrsfs}
\usepackage{comment}
\usepackage[all,cmtip]{xy}
\usepackage{xcolor}
\usepackage{bbm}
\usepackage{enumitem} 
\usepackage{array}
\usepackage{caption}
\usepackage{accents}
\usepackage{scalerel}
\usepackage{a4wide}

\usepackage{circuitikz}
\ctikzset{bipoles/resistor/height=0.15}
\ctikzset{bipoles/resistor/width=0.4}

\definecolor{cadmiumgreen}{rgb}{0.0, 0.42, 0.24}
\usepackage[
colorlinks, citecolor=cadmiumgreen,
pdfauthor={Omid Amini, Hernan Iriarte}, 
pdftitle={Geometry of higher rank valuations},
pdfstartview ={FitV},
]{hyperref}

\usepackage[
alphabetic,
msc-links,
nobysame,
lite,
]{amsrefs} 

\usepackage{tikz, float} 
\usepackage{tikz-cd}
\usetikzlibrary {positioning}

\usetikzlibrary{calc,decorations.markings}
\usetikzlibrary{shapes,snakes}

\renewcommand{\div}{\mathrm{div}}

\DeclareMathOperator{\Supp}{Supp}

\DeclareMathOperator{\Hom}{Hom}

\DeclareMathOperator{\id}{Id}

\newcommand{\Der}{\mathscr D}

\newcommand{\circsm}{{\scaleto{\circ}{2pt}}}
\newcommand{\bethsm}{{\scaleto{\beth}{4pt}}}
\newcommand{\bethd}{\mathpalette\raisebeth\relax}
\newcommand{\raisebeth}[2]{{\scaleto{\beth}{4pt}}\!\!\!\raisebox{0.7pt}{$#1\circsm$}\:}

\newcommand{\dalethd}{\text{\reflectbox{$\bethd$}}}

\newcommand{\LSP}{\mathbf{LSP}}
\newcommand{\CLSP}{\mathbf{CLSP}}
\newcommand{\LSC}{\mathbf{LSC}}

\usepackage[toc]{appendix}

\newcommand{\rest}[1]{\raisebox{-1pt}{$\vert$}_{#1}}

\begin{document} 
	
\theoremstyle{plain}
\newtheorem{thm}{Theorem}[section]
\newtheorem{lem}[thm]{Lemma}
\newtheorem{cor}[thm]{Corollary}
\newtheorem{prop}[thm]{Proposition}
\newtheorem{prop-defn}[thm]{Proposition-Definition}
\newtheorem{question}[thm]{Question}
\newtheorem{claim}[thm]{Claim}
\newtheorem{nota}[thm]{Notation}
\newenvironment{notation}
  {\pushQED{\qed}\renewcommand{\qedsymbol}{$\diamond$}\nota}
  {\popQED\endnota}
  
\theoremstyle{definition}
\newtheorem{defii}[thm]{Definition}
\newenvironment{defn}
  {\pushQED{\qed}\renewcommand{\qedsymbol}{$\diamond$}\defii}
  {\popQED\enddefii}
\newtheorem{remm}[thm]{Remark}
\newenvironment{remark}
  {\pushQED{\qed}\renewcommand{\qedsymbol}{$\diamond$}\remm}
  {\popQED\endremm}
\newtheorem{exx}[thm]{Example}
\newenvironment{example}
  {\pushQED{\qed}\renewcommand{\qedsymbol}{$\diamond$}\exx}
  {\popQED\endexx}

\newtheorem{ex}[thm]{Example}
\newtheorem{conj}[thm]{Conjecture}
\newtheorem{constr}[thm]{Construction}
\newenvironment{construction}
  {\pushQED{\qed}\renewcommand{\qedsymbol}{$\diamond$}\constr}
  {\popQED\endconstr}

\numberwithin{equation}{section}
\newcommand{\eq}[2]{\begin{equation}\label{#1}#2 \end{equation}}
\newcommand{\ml}[2]{\begin{multline}\label{#1}#2 \end{multline}}
\newcommand{\ga}[2]{\begin{gather}\label{#1}#2 \end{gather}}
\newcommand{\mc}{\mathcal}
\newcommand{\mb}{\mathbb}
\newcommand{\surj}{\twoheadrightarrow}
\newcommand{\inj}{\hookrightarrow}
\newcommand{\red}{{\rm red}}
\newcommand{\codim}{{\rm codim}}
\newcommand{\rank}{{\rm rank}}
\newcommand{\Pic}{{\rm Pic}}
\newcommand{\Div}{{\rm Div}}
\newcommand{\im}{{\rm im}}
\newcommand{\Spec}{{\rm Spec \,}}
\newcommand{\Sing}{{\rm Sing}}
\newcommand{\Char}{{\rm char}}
\newcommand{\Tr}{{\rm Tr}}
\newcommand{\Gal}{{\rm Gal}}
\newcommand{\Min}{{\rm Min \ }}
\newcommand{\Max}{{\rm Max \ }}
\newcommand{\ti}{\times }
\newcommand{\sA}{{\mathcal A}}
\newcommand{\sB}{{\mathcal B}}
\newcommand{\sC}{{\mathcal C}}
\newcommand{\sD}{{\mathcal D}}
\newcommand{\sE}{{\mathcal E}}
\newcommand{\sF}{{\mathcal F}}
\newcommand{\sG}{{\mathcal G}}
\newcommand{\sH}{{\mathcal H}}
\newcommand{\sI}{{\mathcal I}}
\newcommand{\sJ}{{\mathcal J}}
\newcommand{\sK}{{\mathcal K}}
\newcommand{\sL}{{\mathcal L}}
\newcommand{\sM}{{\mathcal M}}
\newcommand{\sN}{{\mathcal N}}
\newcommand{\sO}{{\mathcal O}}
\newcommand{\sP}{{\mathcal P}}
\newcommand{\sQ}{{\mathcal Q}}
\newcommand{\sR}{{\mathcal R}}
\newcommand{\sS}{{\mathcal S}}
\newcommand{\sT}{{\mathcal T}}
\newcommand{\sU}{{\mathcal U}}
\newcommand{\sV}{{\mathcal V}}
\newcommand{\sW}{{\mathcal W}}
\newcommand{\sX}{{\mathcal X}}
\newcommand{\sY}{{\mathcal Y}}
\newcommand{\sZ}{{\mathcal Z}}
\newcommand{\TP}{{\mathbb{TP}}}
\newcommand{\A}{{\mathbb A}}
\newcommand{\B}{{\mathbb B}}
\newcommand{\C}{{\mathbb C}}
\newcommand{\D}{{\mathbb D}}
\newcommand{\E}{{\mathbb E}}
\newcommand{\F}{{\mathbb F}}
\newcommand{\G}{{\mathbb G}}
\renewcommand{\H}{{\mathbb H}}
\newcommand{\J}{{\mathbb J}}
\newcommand{\M}{{\mathscr M}}
\newcommand{\N}{{\mathbb N}}
\renewcommand{\P}{{\mathbb P}}
\newcommand{\Q}{{\mathbb Q}}
\newcommand{\R}{{\mathbb R}}
\newcommand{\T}{{\mathbb T}}
\newcommand{\U}{{\mathbb U}}
\newcommand{\V}{{\mathbb V}}
\newcommand{\W}{{\mathbb W}}
\newcommand{\Z}{{\mathbb Z}}
\newcommand{\pic}{{\text{Pic}(C,\sD)[E,\nabla]}}
\newcommand{\ocd}{{\Omega^1_C\{\sD\}}}
\newcommand{\oc}{{\Omega^1_C}}
\newcommand{\al}{{\alpha}}
\newcommand{\be}{{\beta}}
\newcommand{\ta}{{\theta}}
\newcommand{\ve}{{\varepsilon}}
\newcommand{\lr}[2]{\langle #1,#2 \rangle}
\newcommand{\nnn}{\newline\newline\noindent}
\newcommand{\nn}{\newline\noindent}

\newcommand{\onote}[1]{{\color{blue} Omid: #1}}
\newcommand{\hnote}[1]{{\color{orange} Hernan: #1}}
\newcommand{\vs}{\mathcal V} 
\newcommand{\ws}{\mathcal W} 
\newcommand{\spanningf}{\mathcal{SF}}
\newcommand{\gr}{\mathrm{gr}}

\newcommand{\bY}{\mathbf{Y}} 
\newcommand{\ok}{R} 
\newcommand{\ovar}{Y} 
\newcommand{\Diag}{\mathrm{Diag}}
\newcommand{\green}{\mathfrak g}
\newcommand{\p}{\mathbf{p}}
\newcommand{\ps}{\mathbf{p}}
\newcommand{\oA}{\mathfrak A}
\newcommand{\oB}{\mathfrak B}
\newcommand{\oD}{\mathfrak D}

\newcommand{\ad}{\operatorname{ad}}
\newcommand{\K}{\mathbb K}
\newcommand{\an}{\operatorname{an}}
\newcommand{\val}{\mathrm{val}}
\newcommand{\Ar}{\mathrm{Ar}}
\newcommand{\onto}{\twoheadrightarrow}  
\newcommand{\oL}{\mathrm{L}}
\newcommand{\rl}{\mathit{l}}

\newcommand{\TC}{T\mathcal{C}}
\newcommand{\cent}{\mathrm{center}}
\newcommand{\I}{\mathcal I}
\newcommand{\ord}{\mathrm{ord}}
\newcommand{\trop}{\mathrm{Trop}}
\newcommand{\strop}{\mathrm{trop}}
\newcommand{\bir}{\mathrm{bir}}
\newcommand{\e}{\mathfrak e}
\renewcommand{\rank}{\operatorname{rank\,}}

\newcommand{\sgr}{\mathrm{SGrp}}
\newcommand{\TS}{\R}
\newcommand{\contr}[1]{/\!#1}

\newcommand{\geqcw}{\geq_{_{\mathrm{cw}}}}
\newcommand{\leqcw}{\leq_{_{\mathrm{cw}}}}
\newcommand{\minS}[1]{A_{#1}}
\newcommand{\clsp}[1]{\overline{\bf{#1}}}

\newcommand{\lex}{\mathrm{lex}}
\newcommand{\leqlex}{\preceq_{\mathrm{lex}}}
\newcommand{\geqlex}{\succeq_{\mathrm{lex}}}

\newcommand{\Pder}{T}

\newcommand{\coef}{\widetilde \kappa}

\newcommand{\init}{\mathrm{in}} 
\newcommand{\logsmooth}{\mathscr{S}}
\newcommand{\transpose}{{\scaleto{\textrm{t}}{4.4pt}}}
\newcommand{\ret}[1]{{r_{_{\hspace{-.05cm}#1}}}}
\newcommand{\ev}[1]{\mathrm{ev}_{{\hspace{-.05cm}#1}}}
\newcommand{\centre}[1]{{\mathrm{c}_{#1}}}
\newcommand{\centrebis}[1]{{\mathrm{c'}_{#1}}}
\newcommand{\suppaux}[2]{\scalebox{1}[1.4]{$#1\lvert$}#2\scalebox{1}[1.4]{$#1\rvert$}}
\newcommand{\supp}[1]{\mathpalette\suppaux{#1}}
  
  \newcommand{\basefield}{\kappa}

  \newcommand{\st}{\bigm|} 
\newcommand{\Bigst}{\Bigm|} 

  
  \newcommand{\x}{\scaleto{\mathrm{X}}{5.8pt}}
  \newcommand{\y}{\scaleto{\mathrm{Y}}{5.8pt}}

\title[Geometry of higher rank valuations]{Geometry of higher rank valuations}
\author[Omid Amini]{Omid Amini}
\address{CNRS - CMLS, \'Ecole polytechnique, France}
\email{omid.amini@polytechnique.edu}

\author[Hernan Iriarte]{Hernan Iriarte}
\address{Department of Mathematics, University of Texas at Austin, USA.}
\email{iriarte@utexas.edu}

\begin{abstract} The aim of this paper is to introduce a certain number of tools and results suitable for the study of valuations of higher rank on function fields of algebraic varieties. This will be based on a study of higher rank quasi-monomial valuations taking values in the lexicographically ordered group $\R^k$. 
 
We prove a \emph{duality theorem} that gives a geometric realization of higher rank quasi-monomial valuations as \emph{tangent cones of dual cone complexes}. Using this duality,  we provide an analytic description  of quasi-monomial valuations as \emph{multi-directional derivative operators} on tropical functions. 

We consider moreover a refined notion of tropicalization in which we remember the initial terms of power series on each cone of a dual complex, and  prove a tropical analogue of the  \emph{weak approximation theorem} in number theory by showing that any compatible collection of initial terms on cones of a dual cone complex is the refined tropicalization of a rational function in the function field of the variety.

 Endowing the value group $\mathbb{R}^k$ with its Euclidean topology, we study then a natural topology on spaces of higher rank valuations that we call the \emph{tropical topology}. By using the approximation theorem we provide an explicit description of the tropical topology on tangent cones of dual cone complexes. 
 
Finally, we show that tangent cones of dual complexes provide a notion of \emph{skeleton} in higher rank non-archimedean geometry. That is, generalizing the picture in rank one to higher rank, we construct retraction maps to tangent cones of dual cone complexes, and use them to obtain limit formulae in which we reconstruct higher rank non-archimedian spaces with their tropical topology as the projective limit of their higher rank skeleta. 

We conjecture that these higher rank skeleta provide appropriate bases for the study of variations of Newton-Okounkov bodies.

 \end{abstract}

\date{\today}

\maketitle

\setcounter{tocdepth}{1}

\tableofcontents

\section{Introduction} \label{sec:introduction}

The aim of this paper is to introduce a certain number of tools and results suitable for the study of valuations of higher rank on function fields of algebraic varieties. This will be based on finite type approximations of the  valuation spaces  under consideration via a theory of higher rank skeleta developed in this paper, which provides a geometric interpretation of higher rank valuations in terms of tangent cones of cone complexes.

  The motivations behind the study undertaken here are multifold and come 
  \begin{itemize}
  \item on one side, from the theory of Newton-Okounkov bodies and their variations~\cite{Ok96, KK12,LM, Bouck12, BC11, Ami-w, CGGKRUUX, CFKLRS, KM19, RW19, EH, CMM21, Bos21, HKW20}, where one desires to understand continuity and wall-crossing behavior of convex bodies associated to big line bundles when the corresponding defining valuations vary. 
    \item on the other side, from the recent works~\cite{AN20, AN22} involving the hybrid geometry of curves and their moduli spaces, in the constructions of higher rank hybrid and tropical compactifications, and the development of a function theory in higher rank tropical non-archimedean geometry. 
  \end{itemize}

We propose an approach to the first by providing a base space for the study of families of Newton-Okounkov bodies. The second highlights the importance of higher rank non-archimedean tropical geometry in the study of the asymptotic geometry of multiparameter dependent families of complex varieties. In this regard, a framework for higher rank polyhedral geometry intimately related to the content of this paper is developed in the forthcoming companion~\cite{Iri22}.

In the rest of this introduction, we provide an overview of our results and comment on the links to the related works.

All through this paper we fix a field $\basefield$ that we can assume to be algebraically closed. For a positive integer $k\in \N$, we set $[k] \coloneqq \{1, \dots, k\}$.

\subsection{Valuations} We start by explaining the kind of valuations we consider in this paper. 

Let $(\Gamma,\preceq)$ be a totally ordered abelian group and let $K/\basefield$ be a field extension. A \emph{valuation} $\nu$ on $K/\basefield$ with values in $\Gamma$ is a map $\nu\colon  K \to \Gamma \cup\{\infty\}$ which verifies the following properties for any pair of elements $a, b\in K$.
\begin{enumerate}
	\item $\nu(a) = \infty \Longleftrightarrow a=0$.
	
	\item   $\nu(a+b) \geq \min\{\nu(a), \nu(b)\} \qquad \textrm{and} \qquad \nu(ab) = \nu(a)+\nu(b).$
	
	\item $\nu(a)=0$ provided that $a\in \basefield$.
\end{enumerate}

 In this paper we consider the additive group $\R^k$, for a fixed $k \in \N$, endowed with the lexicographic order $\leqlex$ that we simplify to $\preceq$. This is the order defined by saying $x\preceq y$, $x=(x_1, \dots, x_k)$ and $y =(y_1, \dots, y_k)$, if either $x=y$ or 
there is $i\in[k]$ such that $x_j=y_j$ for $j<i$ and $x_i<y_i$. Moreover, we will suppose that $K$ has finite transcendence degree over $\basefield$, that is, we assume the existence of a smooth connected variety $X$ over $\basefield$ such that $K$ is the function field $K(X)$ of $X$. The integer number $k$ will be regarded as an upper bound for the rank for the valuations considered in this paper. The idea to consider valuations of different ranks simultaneously comes from practical situations in the study of degenerations of families of algebraic varieties over higher dimensional bases.

Basic examples of valuations in this setting are the followings:

\, -- (Monomial valuations). Let $X=\A^2 = \Spec(\basefield[\x,\y])$ and  $K=\basefield(\x,\y)$. For $x,y\in \R_+$, there is a unique valuation \[\nu_{x,y}\colon K \rightarrow \R\cup \{\infty\}\]
called \emph{monomial valuation} with respect to $(x,y)$ and 
given by \[\nu_{x,y}(f)\coloneqq\min\bigl\{ix+jy\st c_{ij}\neq 0\bigr\}, \qquad f=\sum_{(i,j)\in \Z^2} c_{ij}\x^i\y^j\in \basefield[\x,\y].\]
Here and all through the paper, $\R_+$ is the set of non-negative real numbers.

\, -- (Divisorial and flag valuations). Suppose $X/\basefield$ is a normal irreducible  variety and let $F\subsetneq X$ be a closed irreducible subvariety of codimension one. The \emph{order of vanishing along $F$} denoted by $\ord_F$ is a rank one valuation on $K=K(X)$, and any positive scalar multiple of $\ord_F$ is called a \emph{divisorial valuation}. More generally, we can consider a flag of normal irreducible subvarieties 
\[\sF\colon \qquad F_0\supsetneq F_1 \supsetneq \dots \supsetneq F_k\]
where $F_0=X$ and $\codim_X(F_\ell)=\ell$, $\ell \in [k]$. Each $F_{\ell}$ thus defines a discrete valuation $\ord_{F_\ell}$ over $K(F_{\ell-1})$. This gives rise to a \emph{flag valuation} $\nu_{\sF}$ of rank $k$ defined as 
	\begin{align}
	\begin{split}\nu_{\sF}\colon K(X)^*&\rightarrow \mathbb{R}^k\\
	f&\mapsto (\mathrm{ord}_{F_1}(f_1), \mathrm{ord}_{F_2}(f_2),\dots,\mathrm{ord}_{F_k}(f_k))\label{val-flag}
	\end{split}
	\end{align}
where $f_1=f$ and $f_{\ell+1} \in K(F_{\ell})$ is the restriction of $f_\ell\cdot t_{\ell+1}^{-\mathrm{ord}_{F_{\ell+1}}(f_\ell)}$ to $F_{\ell+1}$ for $t_{\ell+1}$ a uniformizer for the valuation $\ord_{F_{\ell+1}}$, see for example~\cite{LM, KK12} for more details and for the link to the theory of Newton-Okounkov bodies.

\, -- (Quasi-monomial valuations) We can generalize the first example above by replacing $\A^2$ by any normal irreducible variety $X$ and taking a simple normal crossing (SNC) divisor $D=D_1\cup\dots \cup D_r$ on $X$. This leads to the concept of \emph{quasi-monomial valuations}, which generalizes monomial, divisorial, and flag valuations, as we will see later in Theorem \ref{thm:flagvaluations}.

Consider the \emph{dual cone complex} of the divisor $D$. This is a simplicial cone complex $\Sigma(X,D)$, that sometime we abbreviate to $\Sigma(D)$, in which there is a ray $\rho_i$ corresponding to each component $D_i$ of $D$, and for each subset $I\subseteq [r]$, each  connected component (if any) of the intersection $D_I\coloneqq\bigcap_{i\in I}D_i$ gives rise to a face $\sigma$ with generating rays $\{\rho_i\}_{i\in I}$. More details can be found in Construction \ref{constr:dualcone}.  Each face $\sigma$ of $\Sigma(D)$ thus corresponds to  a connected component of $D_I$, for $I\subset [r]$, that we denote by $D_\sigma$. In this case, we set $I_\sigma=I$ identified as the set of elements $i\in [r]$ such that $D_i$ contains $D_\sigma$. The divisor $D$ being SNC, $D_\sigma$ is normal irreducible and has a generic point $\eta_\sigma$. Moreover, we can 
choose local equations $\{z_i\}_{i\in I}$ for the components $\{D_i\}_{i\in I}$ around $\eta_\sigma$.  

Just as we did for the case of monomial valuations, for the totally ordered abelian group $(\Gamma,\preceq)$, we can pick a vector $\underline{\alpha}=(\alpha_i)_{i\in I}$ with $\alpha_i\in \Gamma_{\succeq 0}$, and define a unique valuation $\nu_{\underline{\alpha}}$ on $K=K(X)$ by requiring 
\[\nu_{\underline{\alpha}}(\prod_{i\in I} z_i^{\gamma_i})\coloneqq\sum_{i\in I}\alpha_i\gamma_i\]
 for any $\gamma = (\gamma_i) \in\Z^I_+$. We can then naturally extend this, first, to the local ring $\sO_{X,\eta_\sigma}$ by taking the minimum over terms of a power series expansion (after passing to the local completion), and then to the full function field. This is the quasi-monomial valuation associated to $D$ and the \emph{weights} $\underline{\alpha}$. Further details can be found in Section~\ref{sec:intro_refined} and Section~\ref{sec:monomial-valuations}.

We denote by $\M^k(D) = \M^k(X, D)$ the set of all \emph{quasi-monomial  valuations of rank bounded by $k$} with  $(\Gamma, \preceq)=(\R^k, \preceq)$. For $k=1$, we further simplify the notation  to $\M(D)$. From the above description, it follows that elements of $\M(D)$ are in bijection with the pairs $(\sigma, \underline\alpha)$ with $\underline \alpha=(\alpha_i)_{i\in I_\sigma} \in \R^{I_\sigma}_+$. This means $\M(D)$ can be naturally identified with $\Sigma(D)$. The above sets come with a natural tower of projection maps 
\[\M(D) \leftarrow \M^2(D)\leftarrow \dots \leftarrow \M^{k-1}(D) \leftarrow \M^k(D) \leftarrow \dots\]
induced by the projection maps to the first $j-1$ coordinates $\R^{j} \rightarrow \R^{j-1}$, $j=2, \dots, k$.

\subsection{Tropicalization}  Let $X$ be a normal irreducible variety and let $D$ be an SNC divisor on $X$. The elements of the dual cone complex $\Sigma(D)$ correspond to quasi-monomial valuations of rank bounded by one on the function field $K(X)$ of $X$. For each rational function $f\in K=K(X)$, we thus get by evaluation a function 
\[\strop(f)\colon \Sigma(D)\rightarrow \R, \qquad \underline\alpha \in \sigma \to \nu_{\underline\alpha}(f), \]
called the tropicalization of $f$. This is a piecewise integral linear function on each cone $\sigma$ of $\Sigma(D)$. 

In this paper we provide an extension of this picture to the case of higher rank quasi-monomial valuations. 
This will be based on a duality theorem we state in the next section which will allow to give a geometric meaning to the space $\M^k(D)$ and the tropicalization map 
\[\strop(f)\colon \M^k(D)\rightarrow \R^k, \qquad \underline\alpha \in \sigma \to \nu_{\underline\alpha}(f), \]
for $\underline \alpha \in ((\R^k)_{_{\succeq 0}})^{I_\sigma}$, where $\succeq$ refers to the lexicographical order.

\subsection{Tangent cone bundles and duality theorem}
The first contribution of this paper is the duality theorem below which provides a geometric realization of  the set $\M^k(D)$ of quasi-monomial valuations  of rank bounded by $k$ as a \emph{tangent cone bundle} on $\Sigma(D)$.  

Consider the projection map $\M^k(D)\rightarrow \M(D)$ which allows to view $\M^k(D)$ as a bundle over $\M(D)=\Sigma(D)$. We have the following geometric characterization of this bundle.

\begin{thm}[Duality theorem]\label{thm:introduality} There is an isomorphism of bundles over $ \M(D) \simeq {\Sigma(D)}$
	\begin{equation}
	\begin{tikzcd}
	\M^k(D)\arrow{r}{\simeq} \arrow[swap]{d}{ } & \TC^{k-1}\Sigma(D) \arrow{d} \\%
	\M(D)\arrow{r}{\simeq} &{\Sigma(D)}
	\end{tikzcd}
	\end{equation}
where $\TC^{k-1}\Sigma(D)$ is defined as the set of all elements of the form $(\underline x;\underline w_1,\dots,\underline w_{k-1})$ where 
\begin{itemize}
\item[-] the \emph{base point} $\underline x$ is a point of $\Sigma(D)$, and
\item[-]  $\underline w_1,\dots,\underline w_{k-1}$ is an ordered set of tangent vectors to $\Sigma(D)$ at $\underline x$ such that we have
\[x+ \varepsilon \underline w_1+\varepsilon^2 \underline w_2+\dots+\varepsilon^r \underline w_r\in \Sigma(D),\]
for any $r\in [k-1]$ and any small enough $\varepsilon>0$.
\end{itemize}
\end{thm}
For a more precise meaning to the above taken sum, we refer to  Section~\ref{sec:tangent_cones}. We call $\TC^{k-1}\Sigma(D)$ the \emph{tangent cone bundle of $\Sigma(D)$ of order $k-1$}.

Using the above correspondence, we give an explicit realization of higher rank quasi-monomial valuations as \emph{directional derivative operators} defined in terms of the corresponding tangent vectors. In order to do this, we equip the cone complex $\Sigma(D)$ with its \emph{structure sheaf} $\sO_{\Sigma(D)}$ which is the sheaf of \emph{tropical functions}. These are continuous functions whose restrictions on each cone $\sigma$ of $\Sigma$ coincide with a piecewise integral linear function defined on that cone.

A rational function $f\in K=K(X)$ induces a global section $\strop(f)$ of the structure sheaf.

\begin{thm}[Duality Theorem, analytic form] Let $(x;\underline w)$ be an element of the tangent cone $\TC^{k-1}\Sigma(D)$ with $\underline w = (w_1,\dots,w_{k-1})$. The valuation $\nu_{x;\underline{w}}$ given by the duality theorem  above is described as
	\begin{align*}
	v_{x;\underline{w}}\colon K(X)&\longrightarrow \mathbb{R}^k \\
	f &\longmapsto \bigl(\strop(f)(x),D_{w_1}\strop(f)(x),\dots, D_{w_1,\dots,w_k}\strop(f)(x)\bigr)
	\end{align*}
where 
\begin{itemize}
\item $D_{w_1}\strop(f)(x)$ is the directional derivative of the function $\strop(f)$ at $x$ in the direction $w_1$, and 
\item recursively, $D_{w_1,\dots, w_{r+1}}\strop(f)(x)$ is the directional derivative of $D_{w_1,\dots,w_{r}}\strop(f)(x)$ seen as a function on the variable $w_r$ in the direction $w_{r+1}$.
\end{itemize}
\end{thm}

\subsection{Topologies on the tangent cone of dual complexes} Notations as before, let $k$ be an integer and consider the tangent cone bundle $\TC^{k-1}\Sigma(D)$. There are four natural topologies one can define on the tangent cone of a dual complex. They all coincide in the case $k=1$, but differ fundamentally for larger values of $k$. We now discuss these topologies.

First note that by the definition of the monomial valuations, we have an injection $\M^k(D) \hookrightarrow (\R^k)_{\succeq 0}^r$. Via the duality theorem, the two first topologies on the tangent cone $\TC^{k-1}\Sigma(D)$ are induced by this injection. Namely,

 \noindent $\bullet$ (\emph{Ordered topology}) This is the  topology on $\TC^{k-1}\Sigma(D) \simeq \M^k(D)$ induced by the ordered topology of  $(\R^k)_{\succeq 0}$.

\noindent $\bullet$ (\emph{Euclidean topology}) This is the topology on $\TC^{k-1}\Sigma(D)$ induced by the Euclidean topology of $(\R^k)_{\succeq 0} \subseteq \R^k$. Equivalently, this is the topology induced by the Euclidean topology of $\Sigma(D)$.

\noindent $\bullet$ (\emph{Hahn-Berkovich topology}) This is the natural topology which appears usually in the context of non-archimedean geometry, that is the coarsest topology which makes continuous all the tropicalization maps 
\[\strop(f) \colon \TC^{k-1}\Sigma(D) \to \R^k_\lex, \qquad f\in K(X)\] 
where $\R^k_\lex$ refers to $\R^k$ equipped with its lexicographically ordered topology. Note that this makes sense for any ordered abelian group $\Gamma$ as the value group for the space of valuations.

\noindent $\bullet$ (\emph{Tropical topology}) This is arguably the most interesting topology one can define on the tangent cone, as it happens to mix the properties of the Euclidean topology on $\R^k$ with those coming from the lexicographic order used in  defining the valuations (see also Section~\ref{sec:okounkov_bodies}). By definition, this is the coarsest topology which makes continuous all the tropicalization maps 
\[\strop(f) \colon \TC^{k-1}\Sigma(D) \to \R^k, \qquad f\in K(X)\] 
in which $\R^k$ is equipped with its Euclidean topology. This topology might be called as well the Hahn-Euclidean topology.

In this paper we provide an explicit description of the tropical topology. This is obtained as a consequence of the tropical weak approximation theorem proved below.

\subsection{Refined tropicalization and tropical weak approximation}\label{sec:intro_refined} Let $D$ be an SNC divisor on $X$. For each cone $\sigma\in \Sigma(D)$ and for each $i\in I_\sigma$, consider a local equation $z_i$ for $D_i$ around $\eta_\sigma$. The family $\{z_i\}_{i\in I_\sigma}$ provides a system of local parameters for the local ring $\widehat{\sO}_{X,\eta_\sigma}$ obtained as the completion of $\sO_{X,\eta_\sigma}$ at its maximal idea. Each element of the local ring 
$\widehat{\sO}_{X,\eta_\sigma}$ admits an \emph{admissible expansion} in the terminology of~\cite{JM12}, that is, an expansion of the form  
\begin{equation}\label{eq:admissible-intro}
	f=\sum_{\beta \in \Z_+^{I_\sigma}}c_\beta z^\beta, \; c_\beta \in \widehat{\sO}_{X,\eta_\sigma},
	\end{equation}
	in which the right hand side is a convergent series with each  coefficient $c_\beta$ either zero or a unit element in $\widehat{\sO}_{X,\eta_\sigma}$. (Here and in what follows, for $\beta\in \Z^r$ with coordinates $\beta_1, \dots, \beta_r$, the notation $z^\beta$ stands for the product $z_1^{\beta_1}\dots z_r^{\beta_r}$.)

	The \emph{support} of the admissible expansion is the set of all $\beta\in\Z_+^{I_\sigma}$ such that $c_\beta$ is not zero.
 
Although an element $f$ has generally infinitely many admissible expansions, we will show later that the set of initial terms of $f$ is invariant under the choice of the expansion and the local parameters. Here an initial term is an element of the support which is minimal for the partial order $\leqcw$ in which a vector $x=(x_i)_{i\in I_\sigma}$ is less than or equal to $y=(y_i)_{i\in I_\sigma}$ if coordinate-wise we have $x_i \leq y_i, i\in I_\sigma$. We denote the initial terms of $f$ by $\minS{f}^\sigma$. Note that the terms in $\minS{f}^\sigma$ form an \emph{antichain} for the partial order $\leqcw$, that is any pair of distinct elements in $\minS{f}^\sigma$ is incomparable relative to $\leqcw$. The antichain $\minS{f}^\sigma$ determines the restriction of $\strop(f)\rest{\sigma}$, that is, we have
\[\strop(f)(x) =\min_{\beta\in \minS{f}^\sigma}\langle x, \beta\rangle\]
where the notation $\langle x, y \rangle$, $x,y\in \R^{I_\sigma}$, stands for the inner product $\sum_{i\in I_\sigma} x_i y_i$, and the minimum is taken over the finite set $A^\sigma_f$.

For a rational function $f \in K(X)$ which belongs to all the local rings $\widehat{\sO}_{X,\eta_\sigma}$, $\sigma \in \Sigma$, we thus get the \emph{refined tropicalization} of $f$ given by the collection $\sA_f\coloneqq\bigl\{\minS{f}^\sigma \,\bigl|\, \sigma \in \Sigma(D) \textrm{ with }f\in \sO_{X, \eta_\sigma}\bigr\}$, the \emph{family of antichains attached to $f$}. Such a family verifies the following:

\noindent $\bullet$ (\emph{Coherence property}) For any inclusion of faces $\tau \subseteq \sigma$, we have the relation 
	\[A_f^\tau=\min_{\leqcw}\, \mathrm{pr}_{_{\hspace{-.05cm} \sigma \succ \tau}}(A_f^\sigma).\]
Here $\mathrm{pr}_{_{\hspace{-.05cm} \sigma \succ \tau}}$ is the projection map $\R^{I_\sigma} \to \R^{I_\tau}$.

\begin{thm}[Tropical weak approximation theorem]\label{thm:approx-intro}  Let $X$ be a smooth quasi-projective variety over a field $k$ and let $D$ be an SNC divisor on $X$. 
	Let $\mathcal{A}=\{A^\sigma\,|\, \sigma \in \Sigma(D)\}$ be a family consisting of finite sets $A^\sigma \subset \Z_+^{I_\sigma}$ such that
	\begin{itemize}
	\item each $A^\sigma$ is an antichain for the partial order $\leqcw$, for $\sigma \in \Sigma(D)$
	\item the family $\mathcal A$ verifies the coherence property, that is for inclusion of faces $\tau \subseteq \sigma$, we have $A^\tau=\min_{\leqcw}\, \mathrm{pr}_{_{\hspace{-.05cm} \sigma \succ \tau}}(A^\sigma)$.
\end{itemize}	 
	Then, there exists a rational function $f\in K(X)$ such that for each cone $\sigma$ of $\Sigma(D)$, we have $f\in \sO_{X,\eta_\sigma}$ and $A^\sigma=A^\sigma_f$.
\end{thm} 

This result might be regarded as a tropical analogue of the weak approximation theorem in number theory.

From the above theorem we deduce the following result.

\begin{cor}\label{cor:approx} Let $X$ be a smooth quasi-projective variety over a field $\basefield$ and let $D$ be a simple normal crossing divisor on $X$. Any tropical function $F$ on the support of the dual cone complex $\Sigma(D)$  is the tropicalization of a rational function $f\in K(X)$.
\end{cor}

As a consequence of the above result and our analytic description of higher rank valuations as multidirectional derivatives of tropical functions, we infer that both the Hahn-Berkovich and tropical topology are intrinsic to the cone complex $\Sigma(D)$, that is, they can be defined more generally for any rational cone complex $\Sigma$. (The intrinsic nature of the two other topologies, the ordered and the Euclidean, is obvious from the definition.)

The following theorem provides a description of the tropical topology. Let $\Sigma$ be a rational cone complex and suppose $\widetilde{\Sigma}$ is a rational subdivision of it. Let $k$ be a positive integer. A set $U\subset \TC^{k-1}\Sigma$ is called a \emph{$\widetilde \Sigma$-open} if $U\cap \TC^{k-1}\sigma$ is open in $\TC^{k-1}\sigma$ with respect to the Euclidean topology  for every cone $\sigma$ of $\widetilde\Sigma$.

\begin{thm}[Characterization of the tropical topology] Notations as above, we have	\begin{enumerate}
		\item For each rational subdivision $\widetilde\Sigma$ of $\Sigma$, the $\widetilde\Sigma$-open sets of $\TC^{k-1}\Sigma$ are open with respect to the tropical topology. 
		\item The union of all $\widetilde\Sigma$-open sets, $\widetilde \Sigma$ a rational subdivision of $\Sigma$, form a basis of opens sets for the tropical topology.
	\end{enumerate}
\end{thm}

\subsection{Spaces of higher rank valuations} \label{sec:spaces-valuations-intro}

 Given a variety $X$ over $\basefield$, the \emph{birational analytification  of $X$ of bounded rank $k$} is the set
		\[X^{\bir,k}\coloneqq\bigl\{\nu\colon K(X)^*\rightarrow \R^k\mid  \nu \text{ is a valuation}\bigr\}\]
	that we endow  with the coarsest topology which makes continuous all the evaluation maps, 	for any $f\in K(X)^*$,
		\begin{align*}
		\ev{f}\colon X^{\bir,k}&\longrightarrow \R^k \\
		\nu &\longmapsto \nu(f).
		\end{align*}
Here, we equip $\R^k$ with its Euclidean topology. Moreover, we define the following subspaces of $X^{\bir, k}$
		\begin{align*}X^{\bethd,k}\coloneqq\bigl \{\nu\in X^{\bir,k}\mid \nu \text{ has a center in } X \bigr\} \\ 
		X^{\dalethd,k} \coloneqq \bigl\{\nu \in X^{\bir,k}\mid \nu \text{ does not have any center in } X\bigr\}
		\end{align*}
		that we endow with the topology induced by that of $X^{\bir, k}$. Recall that for a variety $X$ and a valuation $\nu\colon K(X)\rightarrow \Gamma$, the center of $\nu$, if it exists, is the unique point $x\in X$ such that $\nu$ is non-negative over $\sO_{X,x}$ and strictly positive over its maximal ideal.

		 Notice that $X^{\bir,k}=X^{\bethd,k}\sqcup X^{\dalethd,k}$ and $X^{\bir,k}=X^{\bethd,k}$ if $X$ is proper. In the terminology of Foster and Ranganathan~\cite{FR16}, the space $X^{\bir,k}$ coincides with the subspace of all valuations defined over the generic point in the Hahn analytification of $X$ endowed with the \emph{extended Euclidean topology}. Note that $X^{\bir,1}$ coincides with the birational part of the Berkovich analytification  $X^{\mathrm{an}}$ of $X$. Moreover, the notation $X^{\bethd,k}$ is used in analogy to the analytic space $X^\bethsm$ of Berkovich~\cite{Ber-formal} and Thuillier~\cite{Th}, where the circle is a reminder that we are considering only the birational parts.

We can actually go further and introduce a flag of subspaces of $X^{\bir,k}$ called the \emph{centroidal flag} which interpolates between $X^{\bethd,k}$ and $X^{\bir,k}$. This is done as follows. For $0\leq r\leq k$, we consider the set
		\[\mathscr{F}^rX^{\bir,k}\coloneqq \bigl \{\nu\in X^{\bir,k}\mid \mathrm{proj}_r(\nu) \text{ has a center in }X\bigr\}\]
		where $\mathrm{proj}_r(v)$ is the composition of $v$ with the projection $\R^k\rightarrow \R^r$ to the first $k$ coordinates. In other words, 
		\[\mathscr{F}^rX^{\bir,k}\coloneqq\mathrm{proj}_r^{-1}X^{\bethd,k}.\]
		This gives a decreasing filtration
		\[X^{\bir,k}=\mathscr{F}^0X^{\bir,k}\supseteq \mathscr{F}^1X^{\bir,k} \supseteq \dots \supseteq \mathscr{F}^kX^{\bir,k}=X^{\bethd,k}.\]

	For an SNC divisor $D$ on the variety $X$, the space $\TC^{k-1}\Sigma(D)$ endowed with its tropical topology naturally fits inside $X^{\bethd,k}$. As we will next explain, tangent cone bundles provide a higher rank notion of skeleton for the above spaces of valuations.

\subsection{Tangent cone bundles as higher rank skeleta}  We start by recalling some basic definitions in birational geometry. Let $X$ be a smooth variety over $\basefield$.    A \emph{log-smooth compactification of $X$} is a proper variety $Y$ containing $X$ as an open subvariety such that $Y\setminus X$ is a simple normal crossing divisor on $Y$. A morphism between log-smooth compactifications $Y'$ and $Y$ is a morphism $f\colon Y'\rightarrow Y$ between the underlying varieties such that $f^{-1}(X)=X$ and $f\rest{X}$ is an isomorphism. The category of log-smooth compactifications of $X$ will be denoted by $\LSC_X$. 

A \emph{compactified log-smooth pair} is the data of a pair $\overline{\mathbf{Y}}=(Y,D)$ consisting of a proper variety $Y$ and a simple normal crossing divisor $D\subset Y$ together with a birational map $\varphi\colon Y\dashrightarrow X$  such that
the divisor $D$ can be decomposed as $D = D^\circ + D^\infty$ where $D^\circ$ and $D^\infty$ do not have any component in common, and such that
\begin{itemize}\item[$(i)$]  the domain of definition of $\varphi$ is $Y\setminus D_\infty$, that is, 
	\[\varphi\colon  Y\setminus D^\infty\longrightarrow X\]
	is well-defined and $Y\setminus D^\infty$ is the maximum open set with this property.
	\item[$(ii)$] the pair $(Y\setminus D^\infty, D^\circ\rest{Y\setminus D^\infty})$ is a log-smooth pair for $X$, i.e., $\varphi\rest{Y\setminus D^\infty}$ is a proper morphism from $Y\setminus D^\infty$ to $X$ and the restriction
	\[Y\setminus (D^\circ\cup D^\infty)\longrightarrow X \setminus \varphi(D^\circ)\]
	is an isomorphism.
\end{itemize}

Morphisms between compactified log-smooth pairs can be defined in a natural way. The category of compactified log smooth pairs will be denoted by $\CLSP_X$.

For a compactified log-smooth pair $\overline{\mathbf{Y}}=(Y,D)$, we denote by $\Sigma(\overline{\mathbf{Y}})=\Sigma(Y,D)$ the dual cone complex associated to the divisor $D$ on $Y$. We denote by $\TC^{k-1}\Sigma(\overline{\mathbf{Y}})$ the corresponding tangent cone bundle that we endow with the tropical topology.

Given a compactified log-smooth pair $\clsp{Y}=(Y,D)$ over $X$, as above, the decomposition $D=D^\circ\cup D^\infty$ gives the subcomplex $\Sigma(D^\circ)$ inside $\Sigma(\clsp{Y})$ which we denote by $\Sigma(\clsp{Y}^\circ)$.

	The \emph{centroidal filtration} of $\TC^{k-1}\Sigma{(\clsp{Y})}$ is by definition the filtration
		\[\mathscr{F}^0\TC^{k-1}\Sigma(\clsp{Y})\supseteq \mathscr{F}^1\TC^{k-1}\Sigma(\clsp{Y})\supseteq \dots \supseteq \mathscr{F}^k\TC^{k-1}\Sigma(\clsp{Y})\]
		given for $0\leq r \leq k$ by
		\begin{equation*} \mathscr{F}^r\TC^{k-1}\Sigma(\clsp{Y})\coloneqq \Bigl\{\bigl(x;(w_1, \dots, w_{k-1})\bigr)\in \TC^{k-1}\Sigma(\clsp{Y})\,\bigl|\, \bigl(x;(w_1, \dots, w_{r-1})\bigr)\in \TC^{r-1}\Sigma(\clsp{Y}^\circ)\Bigr\}.
		\end{equation*}

We prove the following theorem. 

\begin{thm} Notations as above, for each compactified log-smooth pair $\overline{\mathbf{Y}}=(Y,D)$ over $X$,  there is a continuous retraction
	\[\ret{{\overline{\mathbf{Y}}}}\colon X^{\bir,k}\longrightarrow \TC^{k-1}\Sigma(\overline{\mathbf{Y}}).\]
	Moreover, the deduced continuous map 
	\[\ret{}\colon  X^{\bir,k}\longrightarrow \lim_{\substack{\longleftarrow\\ \overline{\bf Y}\in \CLSP_X}} \TC^{k-1}\Sigma(\overline{{\bf Y}})\]
	is a homeomorphism. In addition, the limit are compatible with the centroidal filtration on the analytic spaces and on tangent cone bundles. That is, for each $0\leq r \leq k$, we get a homeomorphism
		\begin{equation*} 
		\mathscr{F}^rX^{\mathrm{bir},k}\longrightarrow \lim_{\substack{\longleftarrow\\ \clsp{Y}\in \CLSP_X}} \mathscr{F}^r\TC^{k-1}\Sigma(\clsp{Y}).
		\end{equation*} 
\end{thm}

We note that the above theorem shows that tangent cones with their tropical topology should be regarded as the \emph{higher rank analogue of skeletons} in non-archimedean geometry. Also, remark that the theorem suggests that the space $X^{\bir,k}$ can be regarded as the tangent cone $\TC^{k-1}X^{\bir, 1}$ of $X^{\bir, 1}$ in the Berkovich analytification $X^{\mathrm{an}}$.

The statements of the above theorem hold as well in the case where the spaces in consideration are equipped with the Hahn-Berkovich topology. Due to \emph{mixed nature} of the tropical topology, the arguments in the proof in the case of the tropical topology become more subtle. In particular, we discover a somehow surprising Topology-Mixing Lemma~\ref{lem:mixing-lemma}.

\subsection{Variations of Newton-Okounkov bodies} \label{sec:okounkov_bodies} The above spaces of valuations and the tropical topology seem to be the right topological spaces for the problem of understanding the variations of Newton-Okounkov bodies, as we explain now.

Let $X$ be a smooth projective variety of dimension $d$ and let $L = \sO(E)$ be a big line bundle over $X$. Consider the graded algebra
\[H_\bullet=\bigoplus_{n\geq 0} H_n\]
where $H_n: = H^0(X, \sO(nE))$  is a finite dimensional $\basefield$-vector subspace of $K(X)$.

Each valuation $\nu\in X^{\bir, d}$ gives rise to the corresponding Newton-Okounkov body in $\R^d$ denoted by $\Delta_{\nu}$ and defined by 
\[\Delta_\nu\coloneqq\overline{\bigcup_{n\geq 0} \left \{\frac{\nu(f)}{n}\,\bigl|\, f\in H_n \right \}}.\]

Let $D$ be a simple normal crossing divisor on $X$ and consider the tangent cone $\TC^{d-1}\Sigma(D)$. Consider the space $\mathrm{BC}(\mathbb{R}^d)$ of compact subsets of $\mathbb{R}^d$ endowed with the Hausdorff distance. We get a map 
\begin{align}\begin{split}\label{map:okounkov}
\Delta\colon \TC^{d-1}\Sigma(D)&\longrightarrow \mathrm{BC}(\mathbb{R}^d) \\
(x;\underline w)&\longmapsto \Delta_{\nu_{x,\underline w}}.
\end{split}\end{align}

The study undertaken in this paper has as objective to ultimately prove that the above map $\Delta$ is continuous when $\TC^{d-1}\Sigma(D)$ is endowed with the tropical topology. This topology is actually the only natural one on $\TC^{d-1}\Sigma(D)$ for which one can expect this statement to be both true and non-trivial, as can be verified through basic examples. 

\begin{conj}\emph{Let $L$ be a big line bundle on a projective variety $X$ of dimension $d$. Let $D$ be a simple normal crossing divisor on $X$ with dual cone complex $\Sigma(D)$ of pure dimension $d$. The variation of Newton-Okounkov bodies on $\TC^{d-1}(\Sigma(D))$ is continuous.} 
\end{conj}

We provide in Section~\ref{sec:heuristic} a heuristic argument for the validity of this conjecture. In particular, on those cones of the dual cone complex whose augmented semigroups are finitely generated, the conjecture holds.

\subsection{Related work} In this final section, we make a comparison of our results with the existing ones in the literature.

 The contributions of this paper can be regarded as part of the recent attempts to generalize the framework of tropical and non-archimedean geometry to higher rank valuations.

Analytification of varieties based on valuations has been developed in the pioneering works of Berkovich~\cite{Ber12} and Huber~\cite{huber}. Both spaces are intimately linked with tropical geometry, in the former by means of usual tropicalization and in the latter by means of adic tropicalization~\cite{F16}. More recently, Kedlaya~\cite{Kedlaya} and Foster-Ranganathan~\cite{FR16, FR16b} introduced an alternative analytification directly linked to the one of Berkovich based on higher rank valuations. This last point of view is similar to the one we have adopted as the setting for formulating our results in this paper.

Higher rank tropicalization has been studied by Aroca~\cite{Aroca10},  Banerjee~\cite{Ban15}, Foster-Ranganathan \cite{FR16, FR16b}, Kaveh-Manon \cite{KM19, KM19b}, Escobar-Harada \cite{EH}, and Joswig-Smith~\cite{JS18}. Our work can be regarded as the geometric version of higher rank tropicalization. A framework for higher rank polyhedral and tropical geometry related to the set-up introduced in this paper will appear in the forthcoming paper~\cite{Iri22}. Tangent cone bundles we introduce in this paper and their refinements play a central role in that work.

Geometric tropicalization in rank one has been studied by Hacking-Keel-Tevelev \cite{HKT}, Thuillier \cite{Th}, Abramovich-Caporaso-Payne~\cite{ACP}, and more recently by Ulirsch \cite{Uli17} and Gross~\cite{Gross18}, among others.

A more general framework for tropicalization has been developed in the work of Lorscheid on blueprints~\cite{Lor15}, and in the works of Giansiracusa-Giansiracusa~\cite{Gians, Gians2} and Maclagan-Rinc{\'o}n~\cite{MR}. Because of the level of generality in those works, higher rank tropicalization can be treated using any of the former two frameworks. Tropicalization with values in hyperfields is studied by Viro~\cite{Viro}, Jun~\cite{jun} and Jell-Scheiderer-Yu~\cite{JSY}.

The link between skeletons and tropicalizations in rank one has been thoroughly studied in the works of Gubler-Rabinoff-Werner~\cite{GRW, GRW2}, Macpherson~\cite{macpherson}, and Baker-Payne-Rabinoff~\cite{BPR}. Since skeletons play a central role in connecting complex and non-archimedean geometry, in the study of one-parameter families of complex manifolds, we expect that higher rank analogues of skeleta introduced in this paper, and their polyhedral counterparts further developed in~\cite{Iri22}, will play a central role in the study of multiparameter families of complex manifolds. A systematic study of multiparameter families of Riemann surfaces is undertaken in the series of works~\cite{AN20, AN22}.

The link between tangent cones of dual cone complexes and higher rank valuations established in the paper allows to illuminate the recent work of Kaveh and Manon~\cite{KM19} on Khovanskii bases. In that work, the authors show how to associate to prime cones appearing in the tropicalization of subvarieties of affine spaces higher rank valuations on the coordinate ring of the variety with a finitely generated semigroup. By the work of Gubler-Rabinoff-Werner~\cite{GRW, GRW2}, a prime cone appearing in the tropicalization of a variety can be viewed naturally in the Berkovich analytification of that variety. Moreover, a prime cone can be embedded in a dual cone complex associated to a simple normal crossing divisor. In this regard, the valuations defined by Kaveh and Manon are examples of quasi-monomial valuations studied in this paper. Combined with the duality theorem, this allows to view Kaveh-Manon valuations as points living in the tangent cone of appropriate dual cone complexes, giving them an analytic description.  We refer to~\cite{RW19, Bos21, BLMM17, BFFHL18, IW20, EH} for further results on the connection between tropical geometry, toric degenerations and Khovanskii bases.

The origin of limit theorems  goes back to the work of Zariski~\cite{Zar1, Zar2} on resolution of singularities in dimension two and three using Riemann-Zariski spaces. For tropicalizations, this has been shown in rank one by Payne~\cite{payne} and Foster-Gross-Payne \cite{FGP}. For geometric tropicalizations, this appears in the work by Kontsevich and Tschinkel \cite{KT} (unpublished), Jonsson-Musta\c{t}\u{a} \cite{JM12}, Boucksom-Favre-Jonsson \cite{BFJ16}, and Boucksom-Jonsson \cite{BJ18}. We have been particularly inspired by the work of~\cite{JM12} in establishing our limit theorems. A higher rank version of~\cite{FGP} has been obtained by Foster-Ranganathan~\cite{FR16} in the situation with the Hahn-Berkovich topology. Our limit theorem suggests the statement of the limit theorem in~\cite{FR16} should remain valid for tropicalizations also with respect to the tropical topology of the space of higher rank valuations.  Relative Riemann-Zariski spaces are studied by Temkin~\cite{Temkin, Temkin2}. We refer to the book of Fujiwara-Kato~\cite{FK} for a  detailed discussion of Riemann-Zariski spaces and their applications in rigid geometry.

A version of the duality theorem for the valuative tree was proved by Favre and Jonsson in~\cite{FJ04}. For curves over non-trivially valued fields, this theorem should be compared with the description of tangent directions at points of type 2 in the Berkovich analytification as valuations of rank two on the function field of the curve, a result which can be traced back to Bosch-L{\"u}tkebohmert~\cite{BL85} and Berkovich~\cite{Ber12}. This is also the main ingredient in Thuillier's non-archimedean version of Poincar\'e-Lelong formula for curves~\cite{Thuillier-thesis} and its reformulation as a slope-formula by Baker-Payne-Rabinoff~\cite{BPR2}.

Finally, let us mention that a version of the approximation theorem for curves for non-trivially valued base fields is proved by Baker-Rabinoff~\cite{Br15}. We expect that our theorem should be true in the non-trivially valued case  in any dimension and plan to come back to this setting in a future work.  

\subsection{Organization of the paper} Here is the plan of the paper. In Section~\ref{sec:polyhedral} we introduce dual cone complexes endowed with the sheaf of tropical functions and their tangent cones. We attach to a simple normal crossing divisor on a variety $X$ its corresponding dual cone complex and its tangent cone. In Section~\ref{sec:tropicalization} we recall the definition of tropicalization of rational functions, and explain how to attach a system of antichains to a rational function, leading to a refinement of the definition of tropicalization.

Section~\ref{sec:monomial-valuations} introduces quasi-monomial valuations of a given rank and studies their basic properties. This section contains the proof of the duality theorem and an analytic description of the monomial valuations in terms of directional derivatives along elements in the tangent cones. Section~\ref{sec:approximation} contains the proof of our approximation theorem.  In Section~\ref{sec:tropical-topology}, we study the tropical topology on the tangent cone, and provide an explicit basis of this topology. Section~\ref{sec:spaces-valuations} introduces several spaces of higher rank valuations on the function field of a smooth variety over $\basefield$. The results are used in Sections~\ref{sec:retraction}, \ref{sec:log-smooth} and~\ref{sec:limit} to prove the continuity of the retraction map and the limit formulae.

\subsection{Acknowledgments} We would like to thank S\'ebastien Boucksom, Charles Favre, Alex K\"uronya, Mirko Mauri and Enrica Mazzon for their  constructive comments and  remarks on the first version of this paper. It was pointed to us independently by S\'ebastien Boucksom and by Mirko Mauri and Enrica Mazzon that Corollary~\ref{cor:approx} can be alternatively obtained by  more direct  methods.

 We are grateful to J\'er\^ome Poineau for his very careful reading of the manuscript and all the suggestions which helped us improve the presentation. We especially thank Marco Maculan for his involvement and all the discussions we had during the early stage of this work, and for his helpful comments on the presentation. Finally, we thank Noema Nicolussi for  ongoing collaboration on higher rank non-archimedean and hybrid geometry, as well as for helpful discussions and comments  on the content of this paper. 
\subsection*{Basic notations}

Along the text we work with varieties over an algebraically closed field $\basefield$, that is, integral schemes of finite type over $\basefield$. Points on varieties are not necessarily closed.

We use the notations $\R_+=\{a\in \R\mid a\geq 0\}$ and $\Z_+=\{a\in \Z\mid a\geq 0\}$. For positive integer $d$, we denote $[d]:=\{1, \dots, d\}$.

In the following, we will denote by $\leqcw$ the coordinate-wise partial order on $\Z^I$, that is, given elements $\beta, \beta'\in \Z^I$, we have $\beta\leqcw\beta'$ if and only if $\beta_i\leqcw \beta'_i$ for each $i\in I$. Sometimes, we only use $\leq$ if the partial order is understood from the context.

We write the symbol $a  \gg b$ to indicate that $a$ is large enough compared to $b$.

For a ring $R$, we denote by $R^\times$ the set of invertible elements of $R$.

\section{Cone complexes and tangent cones} \label{sec:polyhedral}
This section introduces the polyhedral geometry concepts used throughout the document. This includes the notion of cone complexes, their sheaf of tropical functions and tangent cones, as well as dual cone complexes associated with simple normal crossing divisors.

\subsection{Cone complexes}
All through this section, the letter $N$ is used for a free $\Z$-module of finite rank, and $M$ denotes the dual of $N$, that is $M = N^\vee \coloneqq\Hom(N, \Z)$. We denote by $N_\R$ and $M_\R$ the corresponding real vector spaces that are dual to each other. The duality pairing between $M$ and $N$ is denoted by $\langle\,, \rangle$. Recall that a saturated sublattice of $N$ is a subgroup $N'$ with the property that $N'_\R\cap N = N'$.
\begin{defn}[Cones and cone complexes]\label{defn:cones}
	\begin{enumerate}
		\item[]
		\item  A \emph{rational polyhedral cone} in $N_\R$ is a set of the form
		\[\sigma=\Bigl\{x\in N_\R\,\bigl|\, \langle x, u_1\rangle\geq 0,\dots,\langle x, u_k \rangle \geq 0 \Bigr\}\]
		for some $u_1,\dots,u_k\in M=N^\vee$. We say that $\sigma$ is \emph{strictly convex} if it does not contain any line in $N_\R$. A \emph{face} of $\sigma$ is a non-empty subset of the form \[\tau=\sigma\cap \Bigl\{x\in N_\R\,\bigl|\,\langle x, u \rangle=0\Bigr\}\] for some $u\in N^\vee$ non-negative over $\sigma$.
		\item A \emph{rational polyhedral cone complex with \emph{(}weak\emph{)} integral structure} is a pair ($\Sigma, \supp{\Sigma}$) where $\supp{\Sigma}$ is a topological space and $\Sigma$ is a family of closed subsets of $\supp{\Sigma}$ such that: 
		\begin{enumerate}
			\item Each $\sigma \in \Sigma$ is enriched with a lattice $N_\sigma$ and an identification of $\sigma$ with a full dimensional rational strictly convex polyhedral cone in $N_{\sigma,\R}$.
			\item These identifications are compatible in the sense that for each element $\sigma\in\Sigma$, faces of $\sigma$ seen as a cone in $N_{\sigma, \R}$ correspond to elements $\tau$ of $\Sigma$. Under this identification, the lattice $N_\tau$ is identified with a saturated sublattice of $N_\sigma$.
			\item As a set we have $\supp{\Sigma}=\bigsqcup_{\sigma\in \Sigma} \accentset{\circ}{\sigma}$ where $\accentset{\circ}{\sigma}$ is the relative interior of $\sigma$.
			\item The intersection of two elements in $\Sigma$ can be written as a union of elements in $\Sigma$.
		\end{enumerate}
			\end{enumerate}
	We call $\supp{\Sigma}$ the \emph{support} of the cone complex, and the elements of $\Sigma$ are called the cones or faces of the cone complex. By an abuse of notation, we will only use $\Sigma$ to refer to the pair $(\Sigma,\supp{\Sigma})$. For each cone $\sigma$, the lattice $N_\sigma$ is called its  \emph{underlying integral structure} and we identify $\sigma$ with its image in $N_{\sigma,\R}$.
	\begin{enumerate}
		\item[(3)] A cone of dimension one in $\Sigma$ is called a \emph{ray} and a cone of maximal dimension is called a \emph{facet}. The set of all rays of $\Sigma$ (resp. of a cone $\sigma$ in $\Sigma$) is denoted by $\Sigma_1$ (resp. $\sigma_1$). More generally, for any integer $k$, we denote by $\Sigma_k$ (resp. $\sigma_k$) the set of all faces of $\Sigma$ (resp. of $\sigma$) of dimension $k$.
	\end{enumerate}
\end{defn}

\begin{figure}[!t]
\centering
   \scalebox{.45}{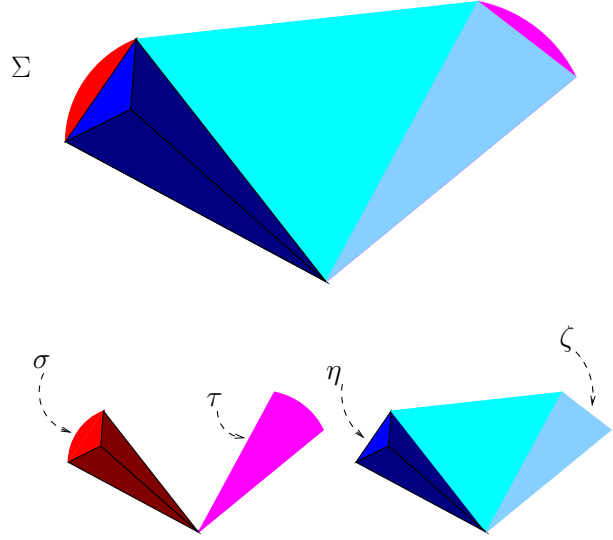}
\caption{Example of a cone complex. 3-dimensional cones $\sigma$ and $\eta$ share two 2-dimensional faces. Similarly, 2-dimensional cones $\tau$ and $\zeta$ share two rays.}
\label{fig:conecomplex}
\end{figure}

\subsection*{Convention} In what follows, a rational strictly convex  polyhedral cone will be simply called  a cone as these are the only kind of cones we will deal with in this paper. Similarly, a rational polyhedral cone complex with a (weak) integral structure will be called simply a cone complex.

\begin{remark} Definition \ref{defn:cones} resembles the notion of a fan used in toric geometry but it differs from it in several ways. First, the lattices $N_\sigma$ may not come simultaneously from a global ambient lattice $N$. Also, condition (d) allows an intersection of two cones to be a union of multiple faces instead of a single face as in the case of fans. In such a situation, the cone complex will have \emph{parallel faces}, that is, two different faces $\tau$ and $\sigma$ in $\Sigma$ with the same set of rays $\tau_1=\sigma_1$. An example of a cone complex is given in Figure~\ref{fig:conecomplex}.
\end{remark}

\begin{defn}[Subdivision] A \emph{rational subdivision} of a cone complex $\Sigma$ is a rational cone complex $\widetilde{\Sigma}$ such that $\supp{\Sigma}=\supp{\widetilde{\Sigma}}$ and for each cone $\tilde{\sigma}\in\widetilde{\Sigma}$, there is a cone $\sigma\in \Sigma$ such that $\tilde{\sigma}\subseteq \sigma$ and $N_{\tilde{\sigma}}$ is a saturated sublattice of $N_\sigma$. 
\end{defn}
It follows from the definition that $\tilde{\sigma}$ is a rational cone  in $N_{\sigma, \R}$.
Again, rational subdivisions are the only ones appearing in this paper, so we drop sometimes the word rational and simply talk about subdivisions.

\subsection{Dual complexes}

\begin{defn}[SNC divisor and stratum]\label{defn:snc} Let $X$ be a smooth variety and $D$ a divisor on it.	
	\begin{enumerate}
		\item The divisor $D$ on $X$ is called \emph{simple normal crossing}, \emph{SNC} in short, if 
		\begin{itemize}
			\item $D$ is reduced, and 
			\item for each point $x\in X$, there is a Zariski neighborhood $U_x$ of $x$ and a regular system of parameters $z_1,\dots,z_r\in \sO_{X,x}$ with $r=\codim\overline{\{x\}}$ such that the zero set of the product $z_1\dots z_{j}$ over $U_x$ coincides with $D\cap U_x$ for some non-negative integer $j=j_x\leq r$.
		\end{itemize}
		
		\item Given an SNC divisor $D$ on $X$, we can write $D$ as a sum $\sum_{i\in \I}D_i$ where $D_i$ are the irreducible components of $D$. A connected component of an intersection of the form \[D_I \coloneqq \bigcap_{i\in I}D_i\]
		for some $I\subseteq \I$ is called a \emph{stratum} of $D$.
	\end{enumerate}
\end{defn}

\begin{remark}
	Notice that each SNC divisor is a Cartier divisor. Moreover, the SNC condition implies that each $D_I$ appearing above is smooth, and in particular, has disjoint irreducible components, coinciding with its connected components.
\end{remark}

\begin{construction}[Dual complex]\label{constr:dualcone} Given a divisor $D=\sum_{i\in \I}D_i$ on a variety $X$ we construct its \emph{dual cone complex} $\Sigma(D)$ as follows. To each stratum $S$ of $D$ which is an irreducible component of $D_I$ for a subset $I\subseteq \I$, one associates a cone $\sigma_S$ which is a copy of $\R_+^I\subseteq \R^\I$ with its natural integral structure given by the lattice $\Z^I \subseteq \Z^{\I}$. If a stratum $S$ is included in another stratum $T$, then the subset $I\subseteq \I$ which corresponds to $S$ should contain the subset $J\subseteq \I$ which corresponds to $T$. In particular, one can naturally identify the cone $\sigma_T$ as a face of $\sigma_S$ via the identification $\R_+^I\subseteq \R_+^J$, as the set of all points with zero coordinates corresponding to elements of $J \setminus I$. The topological space $\supp{\Sigma(D)}$ is defined as the gluing of all $\sigma_S$ along these identifications and the set $\Sigma(D)$ is given as the image of the family $\{\sigma_S\}$ in the space $\supp{\Sigma(D)}$. Sometimes we use the notation $\Sigma(X,D)$
 to emphasize that $D$ is a divisor in $X$. \end{construction}

\begin{prop} The pair $(\Sigma(D),\supp{\Sigma(D)})$ constructed above is a cone complex in the terminology of Definition \ref{defn:cones}.
\end{prop}
\begin{proof} This is straightforward.
\end{proof}

\begin{notation}\rm Notations as above, given a cone $\sigma\in \Sigma(D)$, we denote by $S_\sigma$ the associated stratum. The generic point of $S_\sigma$ is denoted by $\eta_\sigma$. If the divisor is given by $D=\sum_{i\in \I}D_i$, we denote by $I_\sigma$ the subset $I\subseteq \I$ such that $S_\sigma$ is a connected component of $D_I=\bigcap_{i\in I}D_i$.
\end{notation}

\subsection{Tropical functions}
We endow a cone complex $\Sigma$ with its \emph{structure sheaf} $\sO_\Sigma$ which is the sheaf of \emph{tropical functions}.

\begin{defn}[Tropical functions and the structure sheaf] Let $\Sigma$ be a cone complex with an integral structure and let $U$ be an open subset of $\supp{\Sigma}$.
	A function 
	\[F\colon U\rightarrow \R\]
	is called \emph{tropical} if there is a rational subdivision $\widetilde{\Sigma}$ of $\Sigma$ such that for each $\sigma\in \widetilde{\Sigma}$, the restriction $F\rest{\sigma\cap U}$ is integral linear, i.e., viewing $\sigma$ in $N_{\sigma, \R}$,  $F\rest{\sigma\cap U}$ coincides with the restriction to $\sigma \cap U$ of an element in $M_\sigma \subseteq N_{\sigma,\R}^\vee$.  The \emph{structure sheaf} $\sO_\Sigma$ is defined as the one whose sections on an open set $U$ are given by the set of tropical functions on $U$. A tropical function on $\Sigma$ is a global section of $\sO_\Sigma$.
\end{defn}

\begin{remark}\label{remark:tropicalization1} Let $X$  be a smooth variety and $D$ an SNC divisor on $X$. As we will see later, the tropicalization of a rational function on $X$ is a tropical function on $\Sigma(D)$, and any tropical function is of this form.
\end{remark}

\subsection{Tangent cones} \label{sec:tangent_cones}
We now explain how to deal with tangent vectors in cone complexes. We are specially interested in those that \emph{point inward} the cone complex. We start by introducing them in the case in which the cone complex is a single cone and we glue this construction to obtain the general case.

\begin{defn}[The tangent cone] 
	\begin{enumerate}
		\item[]
		\item Let $\sigma\subseteq N_\R$ be a cone and $x\in \sigma$. The \emph{tangent cone at $x$} denoted by $\TC_x\sigma$ is the set of all $w\in N_\R$ for which $x+ \varepsilon w\in \sigma$ provided that $\varepsilon>0$ is small enough.
		\item In the same setting, given an integer $k\geq 1$, we introduce the \emph{$k$-tangent cone at $x$}, denoted by $\TC^k_x\sigma$, as the set of all tuples $\underline{w}=(w_1,\dots,w_k)$ of vectors in $(N_\R)^k$ for which we have the following property:

		For any $r\in [k]$ and for $\varepsilon_i>0$, $i\in[r]$, we have \[x+\varepsilon_1 w_1+\dots+\varepsilon_r w_r\in \sigma\] provided that $\varepsilon_1$ is sufficiently small and $\varepsilon_j$ is sufficiently small with respect to $\varepsilon_{j-1}$ for $1<j\leq r$. Equivalently, if for any small enough $\varepsilon>0$, we have 
		\[x+ \varepsilon w_1+\varepsilon^2 w_2+\dots+\varepsilon^r w_r\in \sigma.\]
		\item The \emph{$k$-tangent cone bundle} is the set $\TC^k\sigma\coloneqq \bigsqcup_{x\in \sigma} \TC^k_x\sigma$. It comes with a natural projection map $\TC^k\sigma\rightarrow \sigma$ and its elements are denoted by $(x;w_1,\dots,w_k)$ or $(x; \underline{w})$, to make reference to the base point explicit. 
		\item We can generalize these constructions to a cone complex $\Sigma$. A face map $\tau\hookrightarrow \sigma$ gives  an
		inclusion $\TC^k\tau\hookrightarrow \TC^k\sigma$.  Therefore, using the fact that \[|\Sigma|=\mathop{\mathrm{colim}}_{\sigma \in \Sigma}\sigma\] where the colimit is taken in the category of sets and goes over all the face maps of $\Sigma$, we introduce the \emph{tangent cone bundle $\TC^k\Sigma$ of $\Sigma$} as
		\[\TC^k\Sigma\coloneqq\mathop{\mathrm{colim}}_{\sigma \in \Sigma} \TC^k\sigma.\]
		The projection maps glue in a natural way to a projection
		\[\TC^k\Sigma\rightarrow |\Sigma|.\]		
	\end{enumerate}
\end{defn}

The topology on $\TC^k\Sigma$ will be discussed in more detail in Section \ref{sec:tropical-topology}.
The next proposition tell us that the tangent cone of a polyhedral complex remains invariant under subdivisions.

\begin{figure}[!t]
		\centering
		\includegraphics[width=1\linewidth]{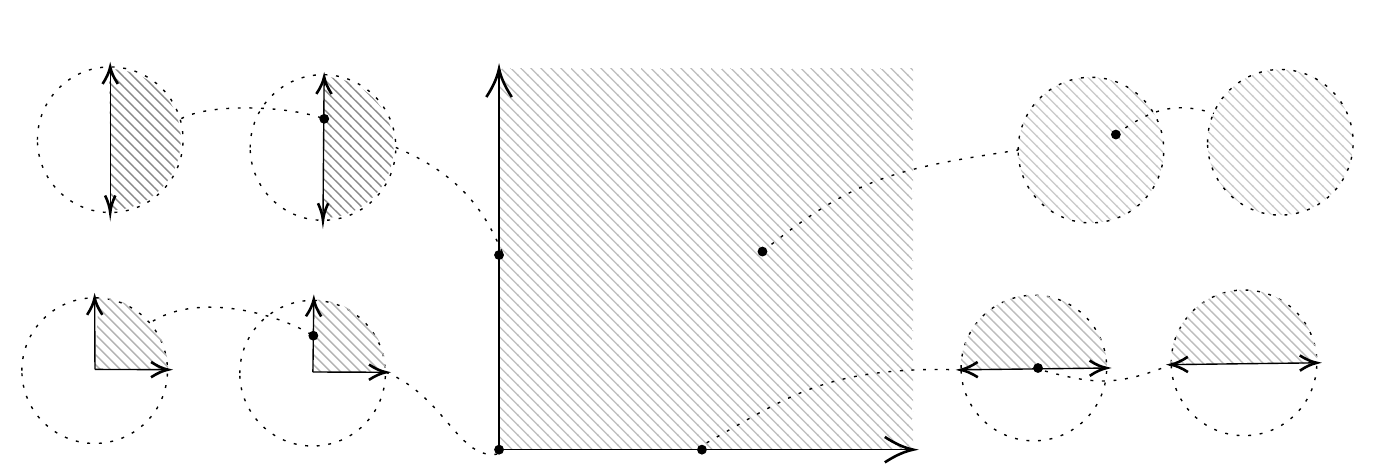}
		\caption{The $2$-tangent cone bundle of a 2-dimensional cone represented via the fibers of the chain of maps $\TC^2\sigma\rightarrow \TC\sigma \rightarrow \sigma$.}
		\label{fig:TCdraw}
	\end{figure}

\begin{example} The $2$-tangent cone bundle of the cone $\sigma=\{(x,y)\in \R^2\mid x,y\geq 0\}$ come equipped with projections \[\TC^2\sigma\rightarrow \TC\sigma \rightarrow \sigma\]
whose fibers are represented in Figure~\ref{fig:TCdraw}.
\end{example}

\begin{prop}\label{prop:subd} If $\widetilde{\Sigma}$ is a subdivision of $\Sigma$, then $\TC^k\Sigma=\TC^k\widetilde{\Sigma}.$
\end{prop}
\begin{proof} Each cone $\sigma\in\Sigma$ becomes the support of the polyhedral complex $\widetilde{\sigma}=\{\tau\in \widetilde{\Sigma}\mid \tau\subseteq \sigma\}.$
	It is enough to prove that $\TC^k\sigma=\TC^k\tilde{\sigma}$, or in other words, 
	\[\TC^k\sigma=\bigcup_{\tau\in \tilde{\sigma}}\TC^k\tau.\]

	The inclusion $\supseteq$ is clear as $\sigma\supseteq \tau$ implies $\TC^k\sigma\supseteq \TC^k\tau$. The other inclusion can be obtained by induction on $k$. For the base case, take $k=0$ and consider $\TC^k\sigma$ as $|\sigma|$. Now, assume the inclusion for $k-1$ and let us prove it for $k$.

	If $(x;w_1,\dots,w_k)\in \TC^k\sigma$, then $x+\frac{w_1}{n}\in \sigma$ for $n$ big enough, so $(x+\frac{w_1}{n};w_2,\dots,w_k)\in \TC^{k-1}\sigma$. By the induction hypothesis, there is some $\tau_n\in \tilde{\sigma}$, depending on $n$, such that \[(x+\frac{w_1}{n};w_2,\dots,w_k)\in \TC^{k-1}\tau_n.\]
	As there are finitely many cones in $\tilde{\sigma}$, there is one $\tau$ on it such that $\tau=\tau_n$ for infinitely many $n$. Now $x+\frac{w_1}{n}\rightarrow x\in \tau$, and then by convexity, if we take $\delta_0, \delta_1>0$ with $\delta_0+\delta_1=1$, for each $\ve_1\gg \dots\gg \ve_{k-1}>0$ small enough we have
	\[\delta_0x+\delta_1(x+\frac{w_1}{n}+\ve_1w_2+\dots+\ve_{k-1} w_k)=x+\delta'_1w_1+\dots+\delta'_kw_k\in \tau.\]
	By choosing $\delta_0,\delta_1$ appropriately we can obtain any $\delta'_1\gg \dots\gg\delta'_k>0$ small enough, hence $(x;\underline{w})\in \TC\tau$ which shows the inclusion we wanted. 
\end{proof}

\section{Tropicalization of rational functions} \label{sec:tropicalization}
We  now recall how to tropicalize rational functions on a variety into tropical functions on cone complexes. 
This is based on the idea that, given a point $x$ in a variety $X$ and a fixed set of local parameters in $\sO_{X,x}$ at the point, the completion $\widehat{\sO}_{X,x}$ of the local ring at $x$ becomes isomorphic to a power series ring in the local parameters. This isomorphism allows  to see each rational function regular at $x$ as a power series. We can then use the usual tropicalization procedure with respect to the trivial valuation on the base field. Following this procedure, given an SNC divisor $D$, we can use the local equations of its components as local parameters to obtain for each rational function a tropical function over $\Sigma(D)$. 

\subsection{Admissible expansions}

The following notion is useful to understand power series expansions directly in the ring $\widehat{\sO}_{X,x}$. It is borrowed from~\cite{JM12}.

\begin{defn}[Admissible expansion]\label{defn:admissible} Let $R$ be a complete regular local $\kappa$-algebra and $z_1,\dots,z_r$ with $r=\dim(R)$ a system of parameters for it. Given $f\in R$, an \emph{admissible expansion for $f$} is an expression of the form
	\begin{equation}\label{eq:admissible}
	f=\sum_{\beta \in \Z_+^r}c_\beta z^\beta, \; c_\beta \in R,
	\end{equation}
	in which the right hand side is a convergent series in which each  coefficient $c_\beta$ is either zero or a unit on $R$. The \emph{support} of the admissible expansion is the set of all $\beta\in\Z^r_+$ with $c_\beta\neq0$.
\end{defn}

Here and in what follows, the notation $z^\beta$ stands for the product $z_1^{\beta_1}\dots z_r^{\beta_r}$ where $\beta_1, \dots, \beta_r$ denote the coordinates of $\beta\in \Z^r$.
\begin{remark} We will be essentially interested in the case in which $R$ is equal to the completion $\widehat \sO_{X,x}$ of the local ring of a point $x$ in a smooth variety $X$. For technical reasons however we have defined it in this generality (see the proof of Proposition \ref{prop:compatibility}).
\end{remark}
\begin{remark}	An element $f\in R$ has several admissible expansions and the support of these admissible expansions may vary. As an example, the identity $1=(1-z^\beta)\cdot 1 +1\cdot z^\beta$ shows two different admissible expansions with different supports for the constant function $1$. Although admissible expansions are not unique, they always exist and as we will see next, the minimal terms of their supports form a uniquely determined set.
\end{remark}

\begin{prop}[Existence of admissible expansions and uniqueness of the minimal elements of the support]\label{prop:mins} Notations as in Definition \ref{defn:admissible}, consider an element $f\in R$.
	\begin{enumerate}
		\item There is an admissible expansion for $f$.
		\item In the notation of \eqref{eq:admissible}, the set
		\[\minS{f}\coloneqq\min_{\leqcw}\{\beta\in \Z^r\mid c_\beta \neq 0\}\]
		depends only on f and not on the choice of the admissible expansion.
		\item The set $\minS{f}$ does not change if we change the local parameters $z_1,\dots,z_r$ for some local parameters $z'_1,\dots,z'_r$ such that we have $z'_i=z_iu_i$ for some unit $u_i\in R^\times$ for each $1\leq i\leq r$. 
	\end{enumerate}
\end{prop}
\begin{remark} A slightly weaker version of this proposition is stated in~\cite{JM12}, where it is shown that the piecewise linear function defined by the admissible expansion is well-defined. Note that it might happen that two power series with different sets of minimal elements give the same piecewise linear function. The above proposition claims the uniqueness of minimal elements in different admissible expansions of a given rational function. 
\end{remark}
\begin{proof}[Proof of Proposition~\ref{prop:mins}]
	(1) Denote by $\basefield(R)$ the residue field of $R$. By Cohen structure theorem~\cite[Theorem 9]{Cohen46}, the ring $R$ contains a coefficient field, that is, a field $\coef \subseteq R$ such that the projection map $R\rightarrow \basefield(R)$ restricts to an isomorphism $\coef\overset{\sim}{\rightarrow} \basefield(R)$. Moreover, this coefficient field induces a continuous isomorphism 
	\begin{equation}\label{eq:cohen-iso}\varphi\colon \basefield(R)[[\x_1,\dots,\x_r]]\overset{\sim}{\longrightarrow}R
	\end{equation}
	which extends the isomorphism between $\basefield(R)$ and $\coef$ by sending $\x_i$ to $z_i$. Writing $\varphi^{-1}(f)=\sum_{\beta\in \Z^r_+}c_\beta \x^\beta$, we get an admissible expansion for $f$ of the form \[f=\sum_{\beta\in \Z^r_+}\varphi(c_\beta) z^\beta.\]

	(2) Notations as in (1), let $f=\sum_{\beta\in \Z^r_+}a_\beta z^\beta$ be a second admissible expansion for $f$. 
		Using the isomorphism \eqref{eq:cohen-iso} above, we can see each $a_\beta$ for $\beta \in \Z_+^{r}$ as a power series with coefficients in $\basefield(R)$, that is as 
	\[\varphi^{-1}(a_\beta)= \sum_{\substack{\gamma \in \Z_{+}^{r}}}a_{\beta,\gamma} \x^\gamma\in \basefield(R)[[\x_1,\dots,\x_r]].\]
	We infer that 
	\begin{align*}
	\sum_{\beta\in \Z^r_+}c_\beta \x^\beta&=\varphi^{-1}(f)
	=\sum_{\beta\in \Z^r_+}\varphi^{-1}(a_\beta) \x^\beta	\\
	&=\sum_{\beta\in \Z^r_+}\left(\sum_{\substack{\gamma \in \Z_{+}^{r}}}a_{\beta,\gamma} \x^\gamma\right) \x^\beta
	=\sum_{\beta\in \Z^r_+}\left(\sum_{0\leqcw \gamma\leqcw \beta} a_{\gamma, \beta-\gamma}\right)\x^\beta
	\end{align*}
	which implies that $c_\beta=\sum_{0 \leqcw \gamma\leqcw \beta} a_{\gamma, \beta-\gamma}$. Now if $\beta$ is a minimal element with $c_\beta\neq 0$, then $a_{\gamma, \beta -\gamma}$ is nonzero for some $\gamma \leqcw \beta$, and therefore $a_\gamma$ is nonzero. Conversely, if $\beta$ is minimal among those $\beta'$ such that $a_{\beta'} \neq 0$, then we have on one side $c_\beta=a_{\beta,0}$, and on the other side, we have $a_{\beta, 0}=\varphi^{-1}(a_\beta) \neq 0$ because $a_\beta$ is a unit. Combined together, we have  shown that any minimal element in the support of one admissible expansion dominates a minimal element in the support of the second.   This proves the statement in the proposition.

	(3) The last point is straightforward. \end{proof}

\begin{remark}\label{remark:finite} 
	\begin{enumerate}
		\item[]
		\item Recall that a subset $A$ of a partially ordered set is called an \emph{antichain} if any pair of distinct elements in $A$ are not comparable in the partial order. It is not hard to prove that an antichain in $(\Z^r_+,\leqcw)$ is necessarily finite. Since the sets $\minS{f}$ considered above are all antichains, we conclude that they must be finite.
		\item For $f,g\in R$, by manipulating admissible expansions, we can see that
		\begin{align*}\min_{\leqcw}\Bigl(\minS{f+g}\cup \minS{f}\cup\minS{g}\Bigr) &= \min_{\leqcw}\bigl(\minS{f}\cup\minS{g}\bigr) \\
		\min_{\leqcw}\Bigl(\minS{f\cdot g}\cup\min_{\leqcw}(\minS{f}+\minS{g})\Bigr)&=\min_{\leqcw}\bigl(\minS{f}+\minS{g}\bigr).\end{align*}
	\end{enumerate}
\end{remark}
\begin{cor}\label{cor:finite}Any function $f\in R$ admits an admissible expansion with finite support.
\end{cor}
\begin{proof} Let $f\in R$ be an admissible expansion. By Proposition \ref{prop:mins},  $f$ admits an admissible expansion $f=\sum_{\beta\in \Z^r_+}c_\beta z^\beta.$
	Rearranging terms, we can rewrite this in the form $f=\sum_{\beta\in \minS{f}}\widetilde{c}_\beta z^\beta$
	where each coefficient $\widetilde{c}_\beta$ can be written in form $\widetilde{c}_\beta= c_\beta+\sum_{\gamma >_{\mathrm{cw}}\beta} c'_\gamma z^{\gamma-\beta}$, for $c'_{\gamma}$ either $0$ or equal to $c_\gamma$, and is still invertible. By Remark \ref{remark:finite}, the set $\minS{f}$ is finite.
\end{proof}

\subsection{Conewise antichains associated to rational functions} Let $D$ be an SNC divisor on $X$. For each cone $\sigma\in \Sigma(D)$ and for each $i\in I_\sigma$, consider a local equation $z_i$ for $D_i$ around $\eta_\sigma$. Then, the family $\{z_i\}_{i\in I_\sigma}$ provides a system of local parameters for the local ring $\widehat{\sO}_{X,\eta_\sigma}$. For a function $f\in K(X)$ with $f\in \sO_{X,\eta_{\sigma}}$, $\sigma\in \Sigma(D)$, we define the set 
\[\minS{f}^\sigma\coloneqq\min_{\leqcw}\{\beta\in \Z^{I_\sigma}\mid c_\beta \neq 0\}\]
for a given (and so for any) admissible expansion $f=\sum_{\beta\in \Z^r_+}c_\beta z^\beta$.

\begin{defn}[Antichains attached to a rational function] Notations as above, for a rational function $f$ on $X$, we call the family $\sA_f\coloneqq\bigl\{\minS{f}^\sigma \,\bigl|\, \sigma \in \Sigma(D) \textrm{ with }f\in \sO_{X, \eta_\sigma}\bigr\}$ the \emph{family of antichains attached to $f$}.
\end{defn}
\begin{remark}
	In practice, we reduce to rational functions $f$ which belong to \emph{all} local rings $\sO_{X, \eta_\sigma}$ for $\sigma \in \Sigma(D)$. In this case, the family of antichains has an element $\minS{f}^\sigma$ for any $\sigma \in \Sigma(D)$. Any more general rational function $h$ on $X$ can be written as the ratio $h=f_1/f_2$ of two such rational functions, i.e., with $f_1, f_2$ belonging both to all the local ring $\sO_{X, \eta_\sigma}$ for $\sigma\in \Sigma(D)$.
\end{remark}
\begin{prop}[Compatibility of the antichains]\label{prop:compatibility-antichain}
	Let $D$ be an SNC divisor on $X$. Fix a cone $\sigma\in \Sigma(D)$ and a face $\tau$ of $\sigma$. Consider the projection
	\begin{align*} \mathrm{pr}_{_{\hspace{-.05cm}\sigma \succ \tau}}\colon \R^{I_\sigma}&\longrightarrow \R^{I_\tau}\\ (x_i)_{i\in I_\sigma}&\longmapsto (x_i)_{i\in I_\tau}.
	\end{align*}
	For each $f\in \sO_{X,\eta_\sigma}$, we then have $f\in \sO_{X,\eta_\tau}$ and an equality of the form
	\[\minS{f}^\tau=\min_{\leqcw}\bigl(\mathrm{pr}_{_{\hspace{-.05cm}\sigma \succ \tau}}(\minS{f}^\sigma)\bigr).\]
\end{prop}
\begin{proof} Consider the diagram

	\[ \begin{tikzcd}  & \widehat{\sO}_{X,\eta_\sigma}\arrow[hook]{dr}{\iota_3} &\\[-12pt]
	\sO_{X,\eta_\sigma} \arrow[hook]{ur}{\iota_1} \arrow[hook]{rd}{\iota_2} & & \widehat{\left(\widehat{\sO}_{X,\eta_\sigma}\right)}_{\mathfrak{p}_\tau}  \\[-12pt]
	& \widehat{\sO}_{X,\eta_\tau} \arrow[hook]{ur}{{\iota_4}} &
	\end{tikzcd}
	\]
	Here $\mathfrak{p}_\tau$ is the prime ideal in $\widehat{\sO}_{X,\eta_\sigma}$ generated by $\{z_i\mid i\notin I_\tau\}$ and each completion is taken with respect to the maximal ideal. Moreover $\iota_1$ is the inclusion in the completion, $\iota_2$ and $\iota_3$ are the compositions of a localization with an inclusion into the corresponding completion, and $\iota_4$ is obtained by functoriality by localizing $\iota_1$ at $\mathfrak{p}_\tau$ and completing with respect to the maximal ideal. This is a commutative diagram of $\kappa$-algebras.
	
	Given an element $f\in \sO_{X,\eta_\sigma}$, by Corollary \ref{cor:finite} we can find finite admissible expansions
	\begin{align*} \iota_1(f)=\sum_{\beta\in \Z^{I_\sigma}}a_\beta z^\beta & & \iota_2(f)=\sum_{\gamma\in \Z^{I_\tau}}b_\gamma z^\gamma
	\end{align*}
	in $\sO_{X,\eta_\sigma}$ and $\sO_{X, \eta_\tau}$, respectively. We then get 
	\begin{align}\label{eq:3.4} \iota_3(\iota_1(f))&=\sum_{\beta\in \Z^{I_\sigma}}\iota_3(a_\beta z^{\beta})=\sum_{\gamma\in \Z^{I_\tau}}\left(\sum_{\mathrm{pr}_{_{\hspace{-.05cm}\sigma \succ \tau}}(\beta)=\gamma}\iota_3(a_\beta z^{\beta- \gamma})\right)z^\gamma \\ \label{eq:3.5} \iota_4(\iota_2(f))&=\sum_{\gamma\in \Z^{I_\tau}}\iota_4(b_\gamma) z^\gamma.
	\end{align}
	Now, we observe that $h_\gamma=\sum_{\mathrm{pr}_{_{\hspace{-.05cm}\sigma \succ \tau}}(\beta)=\gamma}a_\beta z^{\beta- \gamma}\notin \mathfrak{p}_\tau$. Indeed, otherwise, we could construct an admissible expansion for $h_\gamma$ with nonzero terms in some monomials supported in $I_\tau$, and then, these two different admissible expansions would have different set of minimal terms, contradicting Proposition \ref{prop:mins}. Hence, the image of $h_\gamma$ by $\iota_3$ is invertible.

	In this way, Equations \eqref{eq:3.4} and \eqref{eq:3.5} give us two admissible expansions for $\iota_3(\iota_1(f))=\iota_4(\iota_2(f))$ inside $\widehat{\left(\widehat{\sO}_{X,\eta_\sigma}\right)}_{\mathfrak{p}_\tau}$. By Proposition \ref{prop:mins} again, we get
	\[\min_{\leqcw}\Bigl\{\gamma\in\Z^{I_\tau} \,\bigl|\, \iota_4(b_\gamma)\neq 0\Bigr\}=\min_{\leqcw}\textstyle\Bigl\{\gamma\in \Z^{I_\tau}\, \bigl|\, \iota_3(\sum_{\mathrm{pr}_{_{\hspace{-.05cm}\sigma \succ \tau}}(\beta)=\gamma}a_\beta z^{\beta- \gamma})\neq 0\Bigr\}.\]
	Since $\iota_3$ and $\iota_4$ are both injective, we infer $\minS{f}^\tau=\min_{\leqcw}\bigl(\mathrm{pr}_{_{\hspace{-.05cm}\sigma \succ \tau}}(\minS{f}^\sigma)\bigr)$, as required.
\end{proof}

\subsection{Tropicalization} We now define the tropicalization of rational functions.

\begin{construction}[Tropicalization] Let $X$ be a variety and let $D\subseteq X$ be an SNC divisor. Let $\sigma\in \Sigma(D)$ and let $x\in \sigma$.
	
	$\bullet$ For $f\in \sO_{X,\eta_\sigma}$, we define 
	\[\strop(f)(x)\coloneqq \min\bigl\{ \langle x,\beta \rangle \mid \beta\in \minS{f}^\sigma\bigr\}.\] 
	
	$\bullet$ For two elements $f_1, f_2\in \sO_{X,\eta_\sigma}$, we have \[\strop(f_1f_2)(x)=\strop(f_1)(x)+\strop(f_2)(x).\] This allows to extend the above definition to an arbitrary $g\in K(X)$. In this case, we write $g=f_1/f_2$ for $f_1,f_2\in \sO_{X,\eta_\sigma}$ and define for each $x\in \sigma$
	\[\strop(g)(x)\coloneqq \strop(f_1)(x)-\strop(f_2)(x).\]
	
	$\bullet$	Finally, as $\strop(f)$ depends entirely on the family of antichains $\sA_f$, by the compatiblity shown in Proposition \ref{prop:compatibility-antichain} above, $\strop(f)(x)$ is independent of the choice of the face of $\Sigma(D)$ which contains $x$. Hence, we obtain a well defined map
	\[\strop(f)\colon |\Sigma(D)|\rightarrow \R\]
	which we call the \emph{tropicalization of $f$} with respect to $D$.
\end{construction}
\begin{remark}\label{remark:valuation-property} In order to prove the second property, namely, that 
	\[\strop(f_1f_2)(x)=\strop(f_1)(x)+\strop(f_2)(x)\] for $f_1, f_2\in \sO_{X,\eta_\sigma}$, let $\mathrm{in}_x(\minS{f_i}^\sigma)$ be the subset of $\minS{f_i}^\sigma$ consisting of all $\beta$ with $\strop(f_i)(x) =\langle x, \beta\rangle.$ Then, we get $\init_{x}(\minS{f_1f_2}^\sigma) \cap \Bigl(\init_x(\minS{f_1}^\sigma)+ \init_x(\minS{f_2}^\sigma)\Bigr) \neq \emptyset$. Combined with the second part of Remark~\ref{remark:finite}, this gives the result. 
\end{remark}

\begin{prop} The tropicalization of a rational function is a tropical function.
\end{prop}

\begin{proof} For $\sigma\in \Sigma(D)$ and $f\in \sO_{X,\eta_\sigma}$ the tropicalization $\strop(f)\rest\sigma$ is the minimum of finitely many linear functions with integral coefficients. Therefore, this is an integral  piecewise linear function  on $\sigma$. More generally, for any element $f\in K(X)$, the tropicalization $\strop(f)$ can be written as the difference of two integral piecewise linear functions over each cone $\sigma$, and so it is itself integral piecewise linear on each cone. It follows that tropicalization of $f$ is a tropical function. 
\end{proof}

\section{Quasi-monomial valuations of higher rank}~\label{sec:monomial-valuations}

In this section, we define quasi-monomial valuations as certain Krull valuations attached to a given SNC divisor. We study their basic properties and then relate their combinatorial structure with the one of the dual complex in the case the values are taken in $\R^k$ with its lexicographic order.

\subsection{Definition}
We start by giving the definition in the more general setting of totally ordered abelian group. The one important for us in this paper will be the additive group $\R^k$ endowed with the \emph{lexicographic order} $\preceq_\lex$ that we sometimes simply denote by $\preceq$. This is the order defined by $x\preceq_\lex y$ iff $x=y$ or there is an $1\leq i\leq k$ such that $x_j=y_j$ for $j<i$ and $x_i<y_i$. This ordered group has specific properties, depicted in the presence of its two different natural topologies,  which are exploited in this work.

Let $(\Gamma,\preceq)$ be a totally ordered abelian group and consider $\Gamma_+=\{\alpha\in \Gamma\mid \alpha\succeq 0\}$. 	Let $D$ be an SNC divisor in a smooth variety $X$. To a given cone $\sigma\in \Sigma(D)$ and a tuple $\underline{\alpha}\in \Gamma^{I_\sigma}_+$, we associate the valuation $\nu_{\sigma,\underline{\alpha}}$ by defining its value first at an element $f\in \sO_{X,\eta_\sigma}$ by
\begin{equation}\label{eq:monomialdef}
\nu_{\sigma,\underline{\alpha}}(f)\coloneqq\min_{\preceq}\left\{\textstyle\sum_{i\in I_\sigma}\beta_i\alpha_i\in \Gamma\mid \beta\in \minS{f}^{\sigma}\right\}.
\end{equation}
By Remark \ref{remark:finite} and the argument used in~\ref{remark:valuation-property}, it is straightforward to see that 
\begin{align*}
\nu_{\sigma,\underline{\alpha}}(fg)&=\nu_{\sigma,\underline{\alpha}}(f)+\nu_{\sigma,\underline{\alpha}}(g), \quad \textrm{and} \\ \nu_{\sigma,\underline{\alpha}}(f+g)&\succeq \min\{\nu_{\sigma,\underline{\alpha}}(f),\nu_{\sigma,\underline{\alpha}}(g)\}.
\end{align*}
This shows that $\nu_{\sigma,\underline{\alpha}}$  verifies the properties of a valuation on $\sO_{X,\eta_\sigma}$ and so 
uniquely extends to a valuation on $K(X)$, the fraction field of $\sO_{X,\eta_\sigma}$. 

\begin{defn}[Quasi-monomial valuations] Notations as above, the valuation $\nu_{\sigma,\underline{\alpha}}$  is  called the \emph{$\Gamma$-quasi-monomial valuation with respect to $\sigma$ and $\underline{\alpha}$}. The set of all $\Gamma$-quasi-monomial valuations for a given cone $\sigma\in \Sigma(D)$ is denoted by $\M^\Gamma_\sigma(D)$. The set of all $\Gamma$-quasi-monomial valuations coming from any cone of $\Sigma(D)$ is denoted by $\M^\Gamma(D)$.
	
	When the ordered group is the additive group $\R^k$ endowed with the lexicographic order, for a natural number $k$, we call the valuation $\nu_{\sigma,\underline{\alpha}}$ a \emph{quasi-monomial valuation of rank bounded by $k$}. We denote simply by $\M^k_\sigma(D)$ and $\M^k(D)$ the corresponding sets of quasi-monomial valuations $\M^{\R^k}_\sigma(D)$ and $\M^{\R^k}(D)$, respectively. For $k=1$, we further simplify $\M^1_\sigma(D)$ and $\M^1(D)$ to $\M_\sigma(D)$ and $\M(D)$, respectively. 
\end{defn}
In the rest of this paper, we will only consider quasi-monomial valuations of rank bounded by $k$ for some positive integer $k$.

\begin{remark} 
The integer $k$ used in the definition of the quasi-monomial valuation makes reference to the rank of the codomain of the valuation. This should not be confused with the Krull dimension of the valuation ring of $\nu_{\sigma,\underline{\alpha}}$, neither with the rank of the value group of the valuation, as we allow the value group $\nu_{\sigma,\underline{\alpha}}(K(X))$ to be of rank strictly smaller than $k$.  The idea of studying valuations of different ranks all together, simultaneously, is motivated from practical situations appearing in the study of multi-parameter degenerations of complex varieties, see for example~\cite{AN20, AN22}.
\end{remark}

\subsection{The duality theorem} In this section, we provide a dual description of the set of quasi-monomial valuations of rank bounded by $k$. 

Recall that for a variety $X$ and a valuation $\nu\colon K(X)\rightarrow \Gamma$, the \emph{center of $\nu$}, if it exists, is the unique point of $X$ denoted by $\centre{\nu}$ such that $\nu$ is non-negative over $\sO_{X,x}$ and strictly positive over its maximal ideal. The center of a quasi-monomial valuation always exists.

\begin{prop}\label{thm:monomial} Let $D$ be an SNC divisor on a variety $X$ and let $\Gamma$ be a totally ordered abelian group. For $\sigma\in \Sigma(D)$ and $\underline{\alpha}\in \Gamma_{+}^{I_\sigma}$, consider the unique face $\tau$ of $\sigma$ given by the rays $I_\tau=\{i\in I_\sigma\mid \alpha_i\succ 0\}$. Let $\underline \alpha_\tau =\mathrm{pr}_{_{\hspace{-.05cm}\sigma \succ \tau}}(\underline{\alpha})$ be the element $\Gamma_{+}^{I_\tau}$ whose coordinates are given by those of $\alpha$. 
	
	Then, we have $\nu_{\sigma,\underline{\alpha}}=\nu_{\tau,\underline \alpha_\tau}$. Moreover, the center of $\nu_{\sigma,\underline{\alpha}}$ exists and is equal to $\eta_\tau$.
\end{prop}
\begin{proof} The first assertion follows directly from Proposition \ref{prop:compatibility-antichain} and the definition of valuations $\nu_{\sigma,\underline{\alpha}}$ and $\nu_{\tau,\underline \alpha_\tau}$ given in~\eqref{eq:monomialdef}.
		To prove the second, notice that $\nu_{\tau,\underline \alpha_\tau}(f)\succeq0$ for each $f\in \sO_{X,\eta_\tau}$.
	Moreover, $\nu_{\tau,\underline \alpha_\tau}(f)=0$ if and only if $0\in \minS{f}^\tau$, i.e., in the case $f$ is invertible. This shows that the center of  $\nu_{\tau,\underline \alpha_\tau}$ is $\eta_\tau$.
\end{proof}

Consider now the case $\Gamma =\R$. In this case, the elements of $\R_+^{I_\sigma}$ can be naturally identified with the points of $\sigma$. From the compatibility in the above proposition,  we get a natural bijection
\begin{equation}\label{eq:map-iso}
\supp{\Sigma(D)} \longrightarrow \M(D)
\end{equation}
obtained by sending a point $x\in |\Sigma(D)|$ to the  quasi-monomial valuation $\nu_{\sigma,\underline{\alpha}}\in \M(D)$ for $\sigma$ any cone of $\Sigma(D)$ containing $x$, and $\underline{\alpha}$ the coordinates of $x$ in $\sigma$.

We now generalize this bijection to higher rank quasi-monomial valuations.
First observe that there is a natural projection map $\pi\colon \M^k(D)\longrightarrow \M(D)$ defined as follows. Take a point $\underline{\alpha}=(\alpha_i)_{i\in I_\sigma}\in (\R^k)_{+}^{I_\sigma}$.  Each $\alpha_i$ is an element of $(\R^k)_{+} = (\R^k)_{\succeq 0}$ and we denote its coordinates by $\alpha_i=(\alpha_i^1, \dots, \alpha_i^k)$. Consider the projection to the first coordinate denoted by an abuse of the notation by $\pi$ and given by 
\[\pi\colon (\R^k)_{+}^{I_\sigma} \to \R_+^{I_\sigma}, \qquad \pi(\underline{\alpha})=(\alpha_i^1)_{i\in I_\sigma}.\]
The projection map $\pi$ is then defined by   
\begin{equation}\label{eq:definition-pi}\pi(\nu_{\sigma,\underline{\alpha}}) \coloneqq \nu_{\sigma,\pi(\underline{\alpha})} \end{equation}
over each cone $\sigma$ in $\Sigma(D)$. This allows to view $\M^k(D)$ fibered over $\M(D)$.

\begin{thm}[Duality theorem] \label{thm:duality} Notations as above, there is an isomorphism of bundles over $ \M(D) \simeq \supp{\Sigma(D)}$
	\begin{equation}
	\begin{tikzcd}
	\M^k(D)\arrow{r}{\phi} \arrow[swap]{d}{\pi} & \TC^{k-1}\Sigma(D) \arrow{d} \\%
	\M(D)\arrow{r} &\supp{\Sigma(D)}
	\end{tikzcd}
	\end{equation}
	where:
	\begin{itemize}
		\item the map $\M(D) \to \supp{\Sigma(D)}$ on the base is  the inverse of the isomorphism \eqref{eq:map-iso}, and

		\item the map $\phi$ is defined by a compatible family of maps
		\[\phi_\sigma\colon \M_\sigma(D)\longrightarrow\TC^{k-1}\sigma, \qquad \sigma \in \Sigma(D).\]
	 For $\sigma \in \Sigma(D)$, the map $\phi_\sigma$ is defined as follows. Take a point $\underline{\alpha}=(\alpha_i)_{i\in I_\sigma}$ in $(\R^k)^{I_\sigma}_{+}$, let $x= \pi(\underline \alpha)=(\alpha_i^1)_{i\in I_\sigma} \in \R_+^{I_\sigma}$, and for each $j=2, \dots, k$, define 
		\[w_{j-1} \coloneqq (\alpha_i^j)_{i\in I_\sigma} \in \R^{I_\sigma}.\]
		Then, the point $(x;w_1, \dots, w_{k-1})$ belongs to $\TC^{k-1}\sigma$, and we set
		\[\phi_\sigma(\nu_{\sigma,\underline{\alpha}}) \coloneqq (x;w_1, \dots, w_{k-1}).\]
	\end{itemize}
\end{thm}

\begin{remark} In a nutshell, the proof of the duality theorem reduces to the following statement in coordinates. A real matrix $A\in \mathrm{Mat}_{k,r}(\R)$ has columns in $(\R^k)_{+}$, with respect to the lexicographic order on $\R^k$,  if and only if the family $(A^\transpose_{\bullet,1};A^\transpose_{\bullet,2},\dots,A^\transpose_{\bullet,k})$, given by the columns of the transpose $A^\transpose$ of $A$,  belongs to the tangent cone $\TC^{k-1}\Bigl((\R_{+})^r\Bigr)$. This justifies the name given to the theorem. 
\end{remark}

\begin{proof}[Proof of Theorem~\ref{thm:duality}] We verify that each $\phi_\sigma$ is a bijection. Let $\sigma$ be a cone in $\Sigma(D)$. By definition, an element $(\alpha_i)_{i\in I_\sigma} \in (\R^k)^{I_\sigma}$ gives a valuation $\nu_{\sigma,\underline{\alpha}}$ in the domain of $\phi_\sigma$ provided that for each $i\in I_\sigma$, the vector $\alpha_i$ belongs to $(\R^k)_{+}$, that is, it is non-negative with respect to the lexicographic order. Denoting by $(\alpha_i^1, \dots, \alpha_i^k)$ the coordinates of $\alpha_i$, this means that for each $i \in I_\sigma$, we must have 
	\begin{align} 
	\begin{split} \textrm{Either, } &\alpha_i^1> 0, \qquad 
	\text{Or, }  \alpha_i^1=0 \text{ and } \alpha_i^2>0, \qquad   \text{Or, }  \alpha_i^1=\alpha_i^2=0 \text{ and } \alpha_i^3>0 \\ &\dots \qquad  \text{Or, } \alpha_i^1=\dots=\alpha_i^{k-1}=0 \text{ and } \alpha_{i}^k\geq 0. \label{con1}
	\end{split}
	\end{align}
	
	On the other hand, for a collection of vectors $x$, $w_1, \dots, w_{k-1}$  in $\R^{I_\sigma}$, by definition, the family $(x;w_1, \dots, w_{k-1})$  belongs to $\TC^{k-1}\sigma$ if and only if we have
	\begin{align}\label{con2bis}
	\begin{split}
	& x \in \sigma \text{ and }\\ 
	&x +\varepsilon_1 w_1 \in \sigma \text{ for }\varepsilon_1>0 \text{ small enough, and }\\
	&x+\varepsilon_1 w_1+\varepsilon_2 w_2\in \sigma \text{ for } \varepsilon_1 \gg \varepsilon_2>0 \text{ small  enough, and }\\
	&\hspace{4mm}\vdots\\
	&x+\varepsilon_1 w_1+\dots+\varepsilon_{k-1}w_{k-1} \in \sigma \text{ for }\epsilon_{k-2} \gg \varepsilon_{k-1}>0 \text{ small enough}.
	\end{split}
	\end{align}
	
	Specifying the collection of vectors $x, w_1, \dots, w_{k-1}$ to the ones given in the statement of the theorem, the conditions in \eqref{con2bis} above can be rephrased as follows. For each $i\in I_\sigma$, 
	\begin{align}
	\begin{split}
	& \alpha_i^1\geq 0 \text{ and } \\
	& \alpha^1_i+\varepsilon_1 \alpha^2_i\geq 0 \text{ for }\varepsilon_1>0 \text{ small enough, and } \\
	&\alpha^1_i+\varepsilon_1 \alpha^2_i+\varepsilon_2\alpha^3_i\geq 0 \text{ for  }\varepsilon_1 \gg \varepsilon_2>0 \text{ small  enough, and }\\
	&\hspace{4mm}\vdots\\
	&\alpha^1_i+\varepsilon_1 \alpha^2_i+\dots+\varepsilon_{k-1}\alpha^k_i\geq 0 \text{ for } \varepsilon_{k-2} \gg \varepsilon_{k-1}>0 \text{ small enough}.
	\end{split}\label{con2}
	\end{align}
	
	Clearly, Conditions (\ref{con1}) and (\ref{con2}) are equivalent, and  we infer that $\varphi_\sigma$ is a bijection.

	Now to conclude, note that the family of maps $\{\phi_\sigma\}_\sigma$ is compatible with the descriptions of  $\M^k(D)$ and $\TC^{k-1}(D)$ as the unions  $\M^k(D)=\bigcup_\sigma \M^k_\sigma(D)$ and $\TC^{k-1}\Sigma(D)=\bigcup_{\sigma}\TC^{k-1}\sigma$, respectively,  and so they can be glued together  to define a map $\phi\colon \M^k(D)\rightarrow \TC^{k-1}\Sigma(D)$. Since each $\phi_\sigma$ is a bijection, so is $\phi$.
\end{proof}

\subsection{An analytic description of quasi-monomial valuations}

We now explain how to understand higher quasi-monomial valuations from an analytic point of view, directly from tangent vectors, by taking directional derivatives. This leads to a description of the inverse $\phi^{-1}$ of the map $\phi$ appearing in the Duality Theorem.

We need to introduce the notion of \emph{derivative of a function with respect to a tuple of inward tangent vectors} in the tangent cone.

\begin{defn}[Directional derivatives] Given a polyhedral complex $\Sigma$ and a function $F\colon \supp{\Sigma}\rightarrow \R$, the derivative of $F$ at a point $x\in \supp{\Sigma}$ along an inward vector $w\in \TC_x\Sigma$ is the limit
	\[D_wF(x)\coloneqq\lim_{\varepsilon \rightarrow 0^+} \frac{F(x+\varepsilon w)-F(x)}{\varepsilon},\]
	whenever this limit exists. 	More generally, we inductively define the derivative of $F$ at a point $x \in |\Sigma|$ and with respect to the tuple $\underline{w}=(w_1\dots,w_k)\in \TC^k_x\Sigma$ as the limit 
	
	\begin{equation}\label{eq:derivative} D_{{(w_1, \dots, w_k)}}F(x)\coloneqq \lim_{\varepsilon \rightarrow 0^+} \frac{D_{(w_1,\dots,w_{k-1}+\varepsilon w_k)}F(x)-D_{(w_1,\dots,w_{k-1})}F(x)}{\varepsilon},
	\end{equation}
	whenever the directional derivatives  $D_{(w_1,\dots,w_{k-1}+\varepsilon w_k)}F(x)$, for $\varepsilon\geq 0$ small enough, and the above limit exist.
	
	In the case these limits exist for all points $x\in |\Sigma|$ and $\underline w \in \TC^k\Sigma$, we denote by $\Der^kF$ the corresponding \emph{derivative function} from $\TC^k\Sigma \to \R^{k+1}$. This is the function which to a point $x\in |\Sigma| $ and $\underline w = (w_1, \dots, w_k) \in \TC_x^k\Sigma$ associates the point 
	\[\Der^kF(x; \underline w) \coloneqq \bigl(F(x), D_{w_1}F(x), D_{(w_1, w_2)}F(x), \dots, D_{(w_1, \dots, w_k)}F(x)\bigr) \in \R^{k+1}.\]
\end{defn}

\begin{remark}	We make a few remarks.
	\begin{enumerate}
		\item When $(w_1,\dots,w_k)\in \TC ^{k}\Sigma$,  the points $(w_1, \dots, w_{k-2}, w_{k-2} + \epsilon w_{k-1})$, for $\varepsilon\geq 0$ small enough, all belong to $\TC^{k-1}\Sigma$. So the limit in~\refeq{eq:derivative} is well-posed.
		\item At a smooth point $x\in \supp{\Sigma}$, when $x$ lies in the relative interior of a facet of $\Sigma$, and for a function $F\colon \supp{\Sigma}\rightarrow \R$ which is smooth on a neighborhood of $x$, the definition of $D_{\underline{w}}F(x)$ for $\underline w = (w_1,\dots,w_k)\in \TC ^{k}\Sigma$ coincides with the evaluation at the $k$-tuple of tangent vectors $\underline w$ of the $k$-th derivative of $F$ at $x$. The definition is thus a natural extension to the case where $F$ is not necessarily a smooth function and $x$ is an arbitrary point of $\supp{\Sigma}$.  
	\end{enumerate}
\end{remark}
The following proposition provides an alternative way of computing $D_{\underline{w}}F(x)$ when it exists.

\begin{prop} Consider a point $x\in \Sigma$ and a tuple $\underline{w}\in \TC^k_x\Sigma$. Let $F\colon \Sigma\rightarrow \R$ be a function for which $D_{\underline{w}}F(x)$ exists.  Then we have
	\begin{align*}
	D_{\underline{w}}F(x)=\lim_{\varepsilon_k\rightarrow 0^+}\dots \lim_{\varepsilon_1\rightarrow 0^+} \frac1{\varepsilon_1\cdots\varepsilon_k} \Bigl(F(x+&\varepsilon_1w_1+\dots+\varepsilon_1\cdots\varepsilon_kw_k) \\
	&-F(x+\varepsilon_1w_1+\dots+\varepsilon_1\cdots\varepsilon_{k-1}w_{k-1})\Bigr).
	\end{align*}
\end{prop}
\begin{proof} For $k=1$, this is the definition of $D_{\underline{w}}F(x)$. The general case can be obtained by induction.
\end{proof}

In this paper we are mainly interested in directional derivatives of tropical functions. In this case the derivatives always exist as the following proposition shows.

\begin{prop}\label{prop:derivative} For any piecewise linear function $F\colon |\Sigma|\rightarrow \R$ and any $k\geq 0$ the derivative $\Der^kF$ exists. Moreover, if $F$ is linear on a cone $\tau$ such that the point $(x, \underline w) \in \TC^k\Sigma$, $\underline w=(w_1, \dots, w_k)$, belongs to $\TC^k\tau$, then we have
	\begin{equation}\label{eq:linear-derivative}
	\Der^k F(x;w_1,\dots,w_k)=(F_\sigma(x),F_\sigma(w_1),\dots,F_\sigma(w_k)).
	\end{equation}
\end{prop}
\begin{proof} Let $\widetilde{\Sigma}$ be a subdivision of $\Sigma$ such that $F$ is linear on each cone $\sigma\in \widetilde{\Sigma}$. By Proposition~\ref{prop:subd}, we have that $\TC^k\Sigma=\bigcup_{\sigma\in \widetilde{\Sigma}}\TC^k\sigma$, so given $(x;\underline{w})\in \TC^k\Sigma(D)$, there is a cone $\sigma \in \widetilde\Sigma$ such that $(x;\underline{w})\in \TC^k\sigma$. Denote by $F_\sigma$ the linear function which is equal to the restriction of $F$ to $\sigma$. A direct calculation shows that $\Der^k F(x;\underline{w})$ exists and is given by \eqref{eq:linear-derivative}, the which proves the proposition. 
\end{proof}
We now come back to the tropicalization of rational functions and its link to quasi-monomial valuations. From the very definition, it is clear that we can retrieve rank one quasi-monomial valuations by evaluating tropical functions at their corresponding point, that is, given $f\in K(X)^*$ and $x\in |\Sigma(D)|$, if $\nu_x$ denotes the valuation corresponding to $x$ under the map in \eqref{eq:map-iso}, then 
\[\nu_{x}(f)=\strop(f)(x).\]

The following result extends this relation to higher rank quasi-monomial valuations.

\begin{thm}[Quasi-monomial valuations using derivatives] \label{thm:analytic} Let $k\in \N$ be a natural number. Given $(x;\underline{w})\in \TC^{k-1}\Sigma(D)$, consider the evaluation map 
	\begin{align*}
	\nu_{(x;\underline{w})} \colon K(X)^*&\longrightarrow \mathbb{R}^{k} \\
	f&\longmapsto \Der^k\strop(f)(x; \underline w).
	\end{align*}
	Then $\nu_{x;\underline{w}}$ is a well-defined function and it coincides with the valuation $\phi^{-1}(x;\underline{w})$ given by the Duality Theorem~\ref{thm:duality}.
\end{thm}

\begin{proof} Fix a point $(x;\underline{w})\in \TC^{k-1}\Sigma$ and let $\sigma$ be a face of $\Sigma(D)$ containing $x$ such that $(x;\underline{w})\in \TC^{k-1}\sigma$. Then, by definition, for any $f\in \sO_{X,\eta_\sigma}$, we have
	\begin{equation}\label{eq:tropmin}
	\strop(f)(x)=\min\Bigl\{\langle x,\beta\rangle\,\bigl|\, \beta \in \minS{f}^\sigma\Bigr\}.
	\end{equation}
	As $\strop(f)$ is piecewise linear, there is a subdivision $\Sigma_f$ of $\Sigma(D)$ such that $\strop(f)$ is linear on each face of $\Sigma_f$. By Proposition \ref{prop:subd}, there is a cone $\tau$ in $\Sigma_f$ such that $(x;\underline{w})\in \TC^{k-1}_x\tau$. Let $\beta_\tau\in \minS{f}^\sigma$ be the exponent such that $\strop(f)(y)=\langle y,\beta_\tau \rangle$ for any $y\in\tau$. By Proposition \ref{prop:derivative}, we get
	\begin{align}
	\nu_{x;\underline{w}}(f)=\bigl(\langle x,\beta_\tau\rangle,\langle w_{1},\beta_\tau\rangle, \dots,\langle w_{k-1},\beta_\tau \rangle\bigr).
	\label{eq:flag}
	\end{align}
	We now show that $\nu_{x;\underline{w}}=\phi^{-1}(x;\underline w)$, that is,  $\nu_{x;\underline{w}}=\nu_{\sigma,\underline{\alpha}}$ where $\underline{\alpha}=(\alpha_i)_{i\in I_\sigma}$ and $\alpha_i=(x^i,w_1^ i,\dots,w_{k-1}^i)$.

\noindent Note that here for $i\in I_\sigma$, $x^i$ and $w_j^i$ are the $i$-th coordinate of $x$ and $w_j$, respectively. So with our previous notation, we have $\alpha_i^j = x_i$ for $j=1$ and $\alpha_i^j = w^i_{j-1}$ for $j=2, \dots, k$.

\noindent	To show the above claim, note that for any $f\in \sO_{\eta_\sigma}$, we have
	\begin{align}\label{eq:mon}\begin{split}v_{\sigma,\underline{\alpha}}(f)=&\min_{\leqlex}\left\{\sum_{i\in I_\sigma} \alpha_i \beta_i\mid \beta\in \minS{f}^\sigma\right\} \\ =&\sum_{i\in I_\sigma} \alpha_i {\beta_{\underline\alpha,i}} = \sum_{i\in I_\sigma} \left(x^i,w_1^i,\dots, w_{k-1}^i\right) {\beta_{\underline\alpha,i}}  = \bigl(\langle x,\beta_{\underline \alpha}\rangle, \langle w_1,\beta_{\underline \alpha}\rangle, \dots, \langle w_{k-1},\beta_{\underline \alpha}\rangle\bigr)\end{split} \end{align}
	where $\beta_{\underline \alpha}$ is an exponent in $\minS{f}^\sigma$ which gives the minimum in the first equation above, and $\beta_{\underline \alpha, i}$ is the $i$-th coordinate of $\beta_{\underline \alpha}$ for $i\in I_\sigma$. 	We thus need to prove that the two expressions in $(\ref{eq:flag})$ and $(\ref{eq:mon})$ are equal.

	We will prove this by induction. The first entry in both expressions $(\ref{eq:flag})$ and $(\ref{eq:mon})$ coincide as they are both equal to $\strop(f)(x)$. Assuming the two expressions have the same $j$-entries for all $1\leq j \leq \ell-1$, we will prove that the $\ell$-entries are also equal. The first $\ell-1$ entries being equal,
	\begin{align}\label{eq:int-equality}
	\langle x,\beta_\tau \rangle=\langle x,\beta_{\underline \alpha} \rangle, \langle w_1,\beta_\tau \rangle=\langle w_1,\beta_{\underline \alpha} \rangle, \dots,\langle w_{\ell-1},\beta_\tau \rangle=\langle w_{\ell-1},\beta_{\underline \alpha} \rangle,
	\end{align}
	we infer  that  
	\begin{align*}
	&\langle x+\varepsilon_1w_1+\dots+\varepsilon_1\cdots\varepsilon_{\ell}w_{\ell},\beta_\tau\rangle \\ &\qquad = \strop(f) \bigl(x+\varepsilon_1w_1+\dots+\varepsilon_1\cdots\varepsilon_{\ell}w_{\ell}\bigr)\\
	&\qquad =\min_{\leqlex}\Bigl\{\langle x+\varepsilon_1w_1+\dots+\varepsilon_1\cdots\varepsilon_{\ell}w_{\ell},\beta\rangle\,\bigl|\, \beta \in \minS{f}^\sigma\Bigr\}\\
	&\qquad \overset{\star}{=}\min_{\leqlex}\Bigl\{\langle x+\varepsilon_1w_1+\dots+\varepsilon_1\cdots\varepsilon_{\ell}w_{\ell},\beta\rangle\,\bigl|\, \beta \in \minS{f}^\sigma\\
	&\hspace{3cm}\text{such that } \langle x, \beta\rangle=\langle x, \beta_\tau\rangle,\langle w_j, \beta\rangle=\langle w_j, \beta_\tau\rangle \text{ for } 1\leq j \leq \ell-1 \Bigr\}\\
	&=\langle x+\varepsilon_1w_1+\dots+\varepsilon_1\cdots\varepsilon_{\ell}w_{\ell},\beta_{\underline \alpha}\rangle.
	\end{align*}
	Here, in $\overset{\star}{=}$ we used the fact that to minimize $\langle x+\varepsilon_1w_1+\dots+\varepsilon_1\cdots\varepsilon_{\ell}w_{\ell},\beta \rangle$ for $\beta \in \minS{f}^\sigma$ and for $\varepsilon_1 \gg \varepsilon_2 \gg \dots \varepsilon_\ell >0$ small enough, we need to first minimize $\langle x,\beta\rangle$, then minimize $\langle w_1,\beta\rangle$ and so on. By the hypothesis of our induction, $\beta_\tau$ does exactly this as it behaves like $\beta_{\underline \alpha}$ in those entries. From this equality, using Equation~\eqref{eq:int-equality}, we infer the equality $\langle w_\ell,\beta_\tau\rangle = \langle w_\ell, \beta_{\underline\alpha}\rangle$, as required.

	This proves that $\nu_{x;\underline{w}}(f)=\nu_{\sigma,\underline \alpha}(f)$ for all $f\in \sO_{X,\eta_\sigma}$. Using the relation $\strop(f/g)=\strop(f)-\strop(g)$ for two elements $f,g \in \sO_{X,\eta_\sigma}$, we finally conclude that   $\nu_{x;\underline{w}}(f)=\nu_{\sigma,\underline \alpha}(f)$ for all $f\in K(X)^*$ and the theorem follows.
\end{proof}

\subsection{Flag valuations}
In this section, we discuss an alternative way for getting valuations of higher rank on $X$ based on flags of subvarieties, and explain the relation to our constructions above. More details  on valuations associated to flags of subvarieties can be found in~\cite{LM,KK12}, where they are used to define Newton-Okounkov bodies.

Consider a flag of subvarieties 
\[\mathcal{F}\colon \qquad  F_0\supsetneq F_1 \supsetneq \dots \supsetneq F_k\]
where $F_0=X$, and for each $1\leq \ell \leq k$, $F_i$ is a smooth  irreducible subvariety of $F_{\ell-1}$ with  $\mathrm{codim}_X(F_\ell)=\ell$.

Under these hypotheses, each $F_{\ell}$ defines a discrete valuation $\mathrm{ord}_{F_\ell}$ over the function field of $F_{\ell-1}$. We choose a uniformizer $t_\ell$ for $\mathrm{ord}_{F_\ell}$. Using these orders of vanishing, we can construct a higher rank valuation on $K(X)$ as follows.

\begin{prop}\label{eq:flg} Notations as above, consider the map
	\begin{align}
	\begin{split}\nu_{\mathcal{F}}\colon K(X)^*&\rightarrow \mathbb{R}^k\\
	f&\mapsto (\mathrm{ord}_{F_1}(f_1), \mathrm{ord}_{F_2}(f_2),\dots,\mathrm{ord}_{F_k}(f_k))\label{val-flag}
	\end{split}
	\end{align}
	where $f_1=f$ and $f_{\ell+1}$ is the restriction of $f_\ell\cdot t_{\ell}^{-\mathrm{ord}_{F_{\ell}}(f_\ell)}$ to $F_{\ell+1}$ viewed in the function field $K(F_{\ell+1})$. This is a rank $k$ valuation which is independent of the choice of uniformizers $t_\ell$.
\end{prop}
Given a nonempty SNC divisor $D=\sum_{i \in \mathcal I} D_i$, we can define a flag of subvarieties if we fix an ordered sequence of components $D_{i_1},\dots,D_{i_k}$ of $D$, for $i_1, \dots, i_k \in \mathcal I$ with non-empty intersection, and an irreducible component $S$ of the intersection $D_{i_1} \cap \dots \cap D_{i_k}$. In this case, we set $F_0 = X$ and for each $1\leq j \leq k$,  we define $F_j$ as the unique irreducible component of $D_{i_1}\cap \dots \cap D_{i_j}$ which contains $S$. Then we have automatically $F_j \subsetneq F_{j-1}$. 

Since $D$ is SNC, each $F_{j}$ is a smooth connected subvariety of codimension one inside $F_{j-1}$, and we get a flag of subvarieties
\begin{equation}
\mathcal{F}\colon \qquad X=F_0\supsetneq F_1\supsetneq\dots\supsetneq F_k=S \label{eq:F}
\end{equation}
which verifies the hypothesis of Proposition \ref{eq:flg}.

We now prove that $\nu_\mathcal{F}$ corresponds to a quasi-monomial valuation defined in terms of $D$. For this let $\sigma$ be the cone corresponding to the stratum $F_k$ of $D$. This cone has rays indexed by $I_\sigma=\{i_1, \dots, i_k\} \subseteq \mathcal I$. Consider the standard basis $e_{i_1}, \dots, e_{i_r}$ of $N_\sigma$ which is contained in $\sigma$ where $e_{i_j}$ is the primitive vector of the ray corresponding to $i_j$.

\begin{thm}\label{thm:flagvaluations} Let $\mathcal{F}$ the flag on $(\ref{eq:F})$ and $\nu_\mathcal{F}$ the valuation defined by Proposition $\ref{eq:flg}$. Then  $\nu_\mathcal{F}=\nu_{x,\underline{w}}$ for $(x;\underline w) = (e_{i_1};e_{i_2}, \dots, e_{i_k}) \in \TC^{k-1}\Sigma(D)$. 
\end{thm}

\begin{proof} Without loss of generality, we can assume that $D_{i_1}\cap\dots\cap D_{i_k}=\{p\}$ is a closed point of $X$. Indeed, if this is not the case, we can extend the flag in \ref{eq:F} to a complete flag, by adding, if needed, more components to the divisor $D$, and  then work with this complete flag. The result then follows by taking the projection to the first $k$ components of the valuation.

	Now take $z_1,\dots ,z_k$ equations for $D_{i_1},\dots ,D_{i_k}$ around $p$. Using these elements, for each $1\leq r\leq k$, we get a restriction map (called as well \emph{reduction map} in the literature)
	\begin{align*}
	\mathrm{res}_j\colon K(F_{j-1})&\longrightarrow K(F_j),\\
	f&\longmapsto f z_j^{-\mathrm{ord}_{F_j}(f)}{\big|}_{F_j}, 
	\end{align*}
	which satisfies $\mathrm{res}_j(\mathcal{O}_{F_{j-1},x})\subseteq \mathcal{O}_{F_j,x}$.

	The elements $z_{j},\dots,z_k$ give us a local system of parameters for the local ring $\mathcal{O}_{F_{j-1},x}$ and hence they induce an isomorphism $\widehat{\mathcal{O}}_{F_{j-1},x}\simeq k[[z_j,\dots,z_k]]$.
	
	In this way, we obtain  an extension to the power series ring for $\mathrm{res}_j$ as follows
	
	\[ \begin{tikzcd}
	\mathcal{O}_{F_{j-1},x} \arrow{r}{\mathrm{res}_{j}} \arrow[hookrightarrow]{d}{} & \mathcal{O}_{F_j,x} \arrow[hookrightarrow]{d}{} \\%
	\widehat{\mathcal{O}}_{F_{j-1},x}  \arrow["\text{\rotatebox[origin=c]{270}{$\simeq$}}",phantom]{d}&  \widehat{\mathcal{O}}_{F_{j},x}\arrow["\text{\rotatebox[origin=c]{270}{$\simeq$}}",phantom]{d}\\[-5mm] 
	k[[z_j,\dots,z_k]] \arrow{r}{\mathrm{res}_j}& k[[z_{j+1},\dots,z_k]]\\[-7mm]
	f=\sum_\beta a_\beta z^\beta \arrow[mapsto]{r}& \bigl(z_j^{-\mathrm{ord}_{z_j}(f)}\sum_\beta a_\beta z^\beta\bigr)\rest{z_j=0}
	\end{tikzcd}
	\]

	Now, given $f\in \mathcal{O}_{X,x}$, if we write $f=\sum_\beta a_\beta z^\beta \in \widehat{\mathcal{O}}_{X,x}$, then by Theorem \ref{thm:analytic}, we have
	\[\nu_{x,\underline{w}}(f)=\bigl(\langle e_{i_1},\beta_\alpha \rangle, \langle e_{i_2},\beta_\alpha \rangle,\dots,\langle e_{i_k},\beta_\alpha \rangle\bigr)\]
	where $\beta_{\underline \alpha} \in\minS{f}^\sigma$ is the exponent which minimizes the right hand side with respect to the lexicographic order in $\R^k$ given in the proof of that theorem. On the other hand, we have by definition
	\[\nu_\mathcal{F}(f)=\bigl(\mathrm{ord}_{z_1}(f_1),\mathrm{ord}_{z_2}(f_2),\dots,\mathrm{ord}_{z_k}(f_k)\bigr)\]
	where $f_1=f$ and $f_{r}=\mathrm{res}_r(f_{r-1})$.
	
	Now notice that
	\begin{align*}\mathrm{ord}_{z_1}(f)=&\min\bigl\{\langle e_{i_1},\beta\rangle \st \beta\in \mathrm{supp}(f)\bigr\}\\
	=&\min\bigl\{\langle e_{i_1},\beta\rangle \st \beta\in \minS{f}^\sigma\bigr\}\\
	=&\langle e_{i_1},\beta_{\underline \alpha} \rangle.
	\end{align*}
	This shows that the first coordinates of $\nu_{\mathcal F}(f)$ and $\nu_{x,\underline{w}}(f)$ are equal. Proceeding by induction, suppose the first $j$ coordinates of $\nu_\mathcal{F}(f)$ and $\nu_{x,\underline{w}}(f)$ are equal. We get
	\begin{align*}
	\mathrm{ord}_{z_{j+1}}(f_{j+1})=&\, \mathrm{ord}_{z_{j+1}}(\mathrm{res}_j(f_{j}))\\
	=&\, \mathrm{ord}_{z_{j+1}}\left.\left (z_j^{-\mathrm{ord}_{z_j}(f)}\sum_{\beta\in \mathrm{supp}(f_j)} a_\beta z^\beta \right)\right|_{z_j=0}\\
	=& \min\Bigl\{\langle e_{i_{j+1}},\beta\rangle \st \beta\in \mathrm{supp}(f_j) \text{ and } \langle e_{i_j}, \beta \rangle =\langle e_{i_j}, \beta_{\underline \alpha} \rangle \Bigr\}\\
	=& \langle e_{i_{j+1}},\beta_{\underline \alpha}\rangle \qquad \textrm{(by the definition of $\beta_{\underline \alpha}$)}.
	\end{align*}
	This shows the equality between the $j+1$-coordinates of $\nu_{\mathcal F}(f)$ and $\nu_{x,\underline{w}}(f)$. 
	The valuations are thus equal on $\mathcal{O}_{X,x}$, and so they coincide on $K(X)$, as required. 
\end{proof}

\section{Tropical weak approximation theorem}\label{sec:approximation}
The aim of this section is to prove the tropical weak approximation theorem.

\subsection{Statement of the theorem} Recall that a subset $A \subseteq \Z_+^I$ is called an \emph{antichain} for the partial order $\leq = \leqcw$ if any pair of distinct elements $\beta, \gamma \in A$ are incomparable, that is, $\beta \not\leq \gamma$ and $\gamma \not \leq \beta$. By the well-quasi-ordering property of $\leqcw$, this implies $A$ is finite. 

Let $\Sigma =\Sigma(X, D)$ be a dual cone complex. Using the notations previously introduced, for a pair of faces $\tau\prec \sigma$ of $\Sigma$, we denote by $\mathrm{pr}_{_{\hspace{-.05cm} \sigma \succ \tau}}$ the corresponding projection map $\R^{I_\sigma} \to \R^{I_\tau}$.

\begin{defn}[Coherent family of antichains associated to cones]  Suppose for any cone $\sigma$, we have an antichain $A^\sigma \subseteq \Z_+^{I_\sigma}$. We call the collection $\mathcal A = \{A^\sigma\,|\, \sigma \in \Sigma(X, D)\}$ \emph{coherent} if for any inclusion of faces $\tau \subseteq \sigma$, we have the relation 
	\[A^\tau=\min_{\leq}\, \mathrm{pr}_{_{\hspace{-.05cm} \sigma \succ \tau}}(A^\sigma).\]
\end{defn}

\begin{thm}[Tropical weak approximation theorem]\label{thm:approx}  Let $X$ be a smooth quasi-projective variety over a field $k$ and let $D$ be an SNC divisor on $X$. 
	Let $\mathcal{A}=\{A^\sigma\,|\, \sigma \in \Sigma(X, D)\}$ be a coherent family of antichains. 
	There exists then a rational function $f\in K(X)$ such that for each cone $\sigma$ of $\Sigma(X,D)$, we have $f\in \sO_{X,\eta_\sigma}$ and $A^\sigma=A^\sigma_f$.
\end{thm}

\begin{remark}  Stronger versions of this theorem might be true. For example, given admissible expansions $f_\sigma \in \widehat \sO_{X, \eta_\sigma}$ for each $\sigma \in \Sigma(X, D)$ such that each $f_\sigma$ has only finitely many non-zero terms, and such that for inclusion of faces $\tau \prec \sigma$, we have $\iota_{\sigma \succ \tau}(f_\sigma) = f_{\tau}$, one might expect the existence of a rational function $f\in K(X)$ such that $f - f_\sigma$ has an admissible expansion in $\widehat{\sO}_{X, \eta_\sigma}$ in which every monomial  is divisible by a monomial in $f_\sigma$. 
\end{remark}

A corollary of the theorem is the following.

\begin{cor}[Approximation theorem for tropical functions]\label{cor:approx} Let $X$ be a smooth quasi-projective variety over a field $k$ and let $D$ be an SNC divisor on $X$. For any tropical function $F\colon \Sigma(X,D)\rightarrow \R$, there is a rational function $f\in K(X)$ such that $\strop(f)=F$.
\end{cor}

The rest of this section is devoted to the proof of the above theorems. We first prove Theorem~\ref{thm:approx} and then later explain how to deduce the above corollary from this result.

\subsection{Proof of Theorem~\ref{thm:approx} in the toric case} It would be more instructive to first treat the case of a toric variety with the arrangement of the corresponding toric divisors. In this situation, we can drop the quasi-projectivity condition.

Let $\Sigma$ be a unimodular fan of dimension $d$ in the real vector space $N_\R$ of the same dimension, and let $\P_{\Sigma}$ be the corresponding toric variety. Each ray $\varrho$ in $\Sigma$ gives the corresponding divisor $D_\varrho$ in $\P_{\Sigma}$. By unimodularity assumption on $\Sigma$, the divisor $D = \cup_{\varrho \in \Sigma_1} D_\varrho$ is SNC.

Let $\sigma$ be a cone in $\Sigma$, and denote by $\varrho_1,\dots, \varrho_d$ the rays of $\sigma$. Denote the rays of the dual cone $\sigma^\vee$ by $\zeta_1, \dots, \zeta_d$. Let $n_1, \dots, n_d$ be the primitive vectors of $\varrho_1, \dots, \varrho_d$ and denote by $m_1,\dots,m_d$ the primitive vectors of the rays $\zeta_1,\dots,\zeta_d$, respectively. Note that $\langle m_j, n_i \rangle =\delta_{i,j}$, where $\langle .\,,. \rangle$ denotes the duality pairing between $N $ and $M$.

For each point $\underline a = (a_1, \dots, a_d) \in A^\sigma$, consider the rational function 
\[f_{\sigma, \underline a} : = \frac{(\chi^{m_1})^{a_1}\cdots(\chi^{m_d})^{a_d}}{(\chi^{m_1}+\dots+\chi^{m_d}+1)^\ell}\]
for a large enough integer $\ell$ (to be determined later). 

Let $\varrho$ be a ray of $\Sigma$ with primitive vector $n \in N$. The order of vanishing of $f_{\sigma, \underline a}$ along the component $D_{\varrho}$ of $D$ is given by
\[\ord_{D_\varrho}(f_{\sigma, \underline a}) = \langle a_1m_1+\dots+a_dm_d , n \rangle -\ell \cdot \min\{0,\langle m_1, n \rangle,\dots,\langle m_d, n \rangle\}.
\]

In particular, taking $\varrho = \varrho_j$, $j\in [d]$, we get $\ord_{D_{\varrho_j}}(f_{\sigma, \underline a}) = a_j$. Moreover, if $\varrho$ is different from $\varrho_1, \dots, \varrho_d$, then,  by duality, there exists $j \in [d]$ such that $\langle m_j, n\rangle <0$. Upon the choice of $\ell$, this imposes $f_{\sigma, \underline a}$ to have a large order of vanishing along $D_{\varrho}$.

Consider now the rational function  $f_\sigma$ in $K(\P_{\Sigma})$ defined as
\[f_\sigma \coloneqq \sum_{\underline a\in A^{\sigma}} f_{\sigma, \underline a} = \sum_{\underline a = (a_1,\dots,a_d)\in A^\sigma}\frac{(\chi^{m_1})^{a_1}\cdots(\chi^{m_d})^{a_d}}{(\chi^{m_1}+\dots+\chi^{m_d}+1)^\ell}.\]

In the completed local ring $\widehat{\sO}_{\P_\Sigma, x_\sigma}$ we have the equality
\[
\frac{(\chi^{m_1})^{a_1}\cdots(\chi^{m_d})^{a_d}}{(\chi^{m_1}+\dots+\chi^{m_d}+1)^\ell}
=
(\chi^{m_1})^{a_1}\cdots(\chi^{m_d})^{a_d}\cdot \left(1+\sum_{k\geq 1} (-1)^k(\chi^{m_1}+\dots+\chi^{m_d})^k\right)^\ell,
\]
which gives an admissible expansion of $f_\sigma$ with respect to the local parameters $\chi^{m_1}, \dots, \chi^{m_d}$ around $\eta_\sigma$, the point of intersection of $D_{\varrho_1}, \dots, D_{\varrho_d}$.

From this, we see that $A^{\sigma}_{f_{\sigma, \underline a}} = \{\underline a\}$, and since $A^{\sigma}$ is an antichain, it follows $A^\sigma_{f_\sigma}=A^\sigma$.

Now let $\tau$ be another facet and denote by $\{\rho_j\}_{j=1}^d$ its rays. They correspond to the components $D_{\rho_1}, \dots, D_{\rho_d}$ of $D$ with the torus-invariant point $
\eta_\tau$ as the point of intersection. From the preceding discussion, we infer that if $\rho_j$ is not a ray of $\sigma$, then, choosing $\ell$ large enough, we can ensure that $f_\sigma$ has a large order of vanishing along $D_{\rho_j}$. 

Since the order of vanishing of $f_\sigma$ along such a component $D_{\rho_j}$ is equal to the minimum $j$-th coordinate of any element  of $A^\tau_{f_\sigma}$, we see that all the elements of $A^\tau_{f_\sigma}$ have large $j$-th coordinates. On the other hand, on the intersection face $\delta = \tau \cap \sigma$, we have 
\[\mathrm{pr}_{\tau \succ \delta}(A^\tau_{f_\sigma}) = \mathrm{pr}_{\sigma \succ \delta}(A^\sigma_{f_\sigma}). \]

In particular, by the coherence of the collection $A^{\sigma}$, this shows that if $\ell$ is chosen to be large enough, then any element in $A^{\tau}_{f_{\sigma}}$ dominates an element of $A^{\tau}$, that is,
\[A^\tau = \min_{\leq}\Bigl(A^{\tau} \cup A^{\tau}_{f_\sigma}\Bigr).
\]

Now we choose $\ell$ large enough, and taking $f_\sigma$ as above for each facet $\sigma$ of $\Sigma$, for generic choices of $\lambda_\sigma$ in the base field, we set  
\begin{equation}\label{eq:function_approximation_toric} f\coloneqq\sum_{\sigma \in \Sigma_d}\lambda_\sigma f_\sigma.
\end{equation}

We observe that for any facet $\tau$, for any pair of rational functions $h, g$, and generic choice of scalars $\lambda, \mu$, we have  
\[A^\tau_{\lambda h+\mu g}=\min_{\leq} \bigl(A^\tau_h\cup A^\tau_g \bigr). \] 
We thus infer that for any facet $\tau$ of $\Sigma$ and for the function $f$ defined in \eqref{eq:function_approximation_toric}, we have
\[A^\tau_f=\min_{\leq}\left(\bigcup_{\sigma \in \Sigma_d} A^\tau_{f_\sigma}\right)= A^\tau_{f_\tau}=A^\tau.\] 
The result follows now by the coherence property which implies $A^\delta_f=A^\delta$ for any face $\delta$ of $\Sigma$.

\subsection{Proof of Theorem~\ref{thm:approx}}  We now treat the theorem in its full generality.   In the following, we will use the following terminology borrowed from lattice theory concerning the combinatorial structure of faces in a cone complex.

\begin{defn}[The (multivalued) meet and join operations $\wedge$ and $\vee$] \label{defn:sup_inf} \rm 	
	Given two faces $\tau$ and $\sigma$ in a cone complex $\Sigma$, we denote by $\tau\wedge \sigma$ the set of all maximal common faces between $\tau$ and $\sigma$. If $\tau$ and $\sigma$ are faces of a cone $\zeta$, we denote by $\tau\vee_\zeta \sigma$ the unique minimal face of $\zeta$ that contains both $\tau$ and $\sigma$. Notice that if $\Sigma$ does not have parallel faces, then $\tau\wedge \sigma$ is a single cone and $\tau\vee_\zeta \sigma$ is independent of $\zeta$. In this case, we denote this cone by $\tau\vee\sigma$.
\end{defn}

In the rest of this section we assume given an SNC divisor $D = \sum_{i\in \I} D_i$ in $X$, and we consider the dual cone complex $\Sigma(X, D)$.

\subsubsection{Adapted family of rational functions} Proceeding somehow similarly as in the proof of the toric case, we will prove the existence of a family of rational functions with nice properties depicted in the following theorem. 

\begin{thm} \label{thm:adapted} Let $\sigma$ be a face of $\Sigma(X, D)$. There exists a rational function $u_\sigma\in K(X)$ with the following properties:
	
	\begin{itemize}
		\item[(P1)] $u_\sigma$ belongs to all local rings $\sO_{x,\eta_\delta}$, for $\delta\in \Sigma(X, D)$, and is invertible in $\sO_{X,\eta_\sigma}$.

		\item[(P2)] It has a zero along the divisor $D_j$ for each $j\notin I_\sigma$.

		\item[(P3)] For any face $\tau$ of $\Sigma(X,D)$ the following holds. If $\zeta \in \tau \wedge \sigma$ is a maximal common face of $\tau$ and $\sigma$, then the restriction ${u_\sigma}_{|S_\zeta}$ of $u_\sigma$ on the stratum $S_\zeta$ has a zero along all the strata $S_{\zeta\vee_\tau \varrho_j} \subseteq S_\tau$ for any $j\in I_\tau \setminus I_\zeta$. 
	\end{itemize}
\end{thm}

\begin{defn}[Adapted family of rational functions] Given a dual cone complex $\Sigma(X, D)$, the collection of rational functions $u_\sigma$, $\sigma \in \Sigma(X, D)$, verifying the properties (P1), (P2), and (P3) in the above theorem is called an \emph{adapted family of rational functions} for the dual cone complex.
\end{defn}

In order to prepare for the proof of the above theorem, we start by stating two lemmas concerning the existence of rational functions with prescribed regularity on a given finite set of points. 

\begin{lem}\label{lem:1} Let $Y \subseteq X$ be a closed irreducible set and let $x$ be a (non-necessarily closed) point in $X\setminus Y$. Then there exists an irreducible divisor $E \subset X$ which contains $Y$ but not $x$.
\end{lem}
\begin{proof}
	Denote by $\eta_Y$ the generic point of $Y$. As $x\notin Y$, we have $\sO_{X,x}\setminus \sO_{X,\eta_Y}\neq \varnothing$. Take $f\in\sO_{X,x}\setminus \sO_{X,\eta_Y}$. Then, $\eta_Y$ is contained in the indeterminacy set of $f$. Since $X$ is smooth, the indeterminacy set is the support of the negative part of $\mathrm{div}(f)$, and so there is a component $E$ of this negative part which contains $Y$ but not $x$.
\end{proof}

\begin{lem}\label{lem:2} Suppose $X$ is quasi-projective and let $E \subset X$ be a reduced divisor.  Given points $x_1,\dots,x_n$ in $E$, and a point $x\notin E$, there is a rational function $u$ which vanishes on each component of $E$ with order of vanishing one, which belongs to each local ring $\mathcal{O}_{X,x_i}$, $i=1,\dots, n$, and which is invertible at $x$.
\end{lem}

\begin{proof} Taking a projective compactification, we can assume without loss of generality that $X$ is projective. Consider an ample divisor $H$ not containing any of the points $x,x_1,\dots,x_n$ and not sharing any component with $E$. Then, for some large integer number $n$, the divisor $nH-E$ is very ample, and so base point free. Therefore, there is a section $u$ of $\sO(nH-E)$ which does not vanish on $x$. Taking $u$ generic, the corresponding rational function satisfies all the required properties.
\end{proof}

We are now ready to prove the existence of adapted families of rational functions. 
\begin{proof}[Proof of Theorem~\ref{thm:adapted}]
	We first apply Lemma \ref{lem:1} to each stratum $S_\tau$ not contained in $S_\sigma$ to get an irreducible divisor $E_\tau \subset X$ which contains $S_\tau$ but not $S_\sigma$. Let 
	\[E \coloneqq \sum_{\tau: \,\,S_\tau \not\supset S_\sigma}E_\tau. \]
	We now apply Lemma \ref{lem:2} to $E$, the points $\eta_{\tau}$ for $\tau\in \Sigma(X, D)$ with $S_\tau \not\supset S_\sigma$, that is, $\tau \not\prec \sigma$, and the point $\eta_\sigma$, which clearly does not belong to $E$. We infer the existence of  a rational function $u_\sigma$ in $K(X)$ that vanishes on each component $E_\tau$ of $E$, which belongs to $\sO_{X,\eta_\delta}$ for any point $\delta \in \Sigma(X, \Sigma)$, and which is invertible in $\sO_{X,\eta_\sigma}$. We claim $u_\sigma$ verifies the claimed properties (P1), (P2), and (P3).

	The first claim (P1) is clearly satisfied by the construction of $u_\sigma$.

	Also, notice that if $j\notin I_\sigma$, then $S_{\varrho_j} = D_j$. Since both $D_j$ and $E_{\varrho_j}$ are irreducible, we must have $E_{\varrho_j}=D_j$. Since $j\notin I_\sigma$, we have $\eta_\sigma \notin D_j$, and so by the choice of $u_\sigma$, it should vanish on $E_{\varrho_j}$. This shows that $u_\sigma$ verifies Property (P2).

	Finally, let $\zeta$ be a maximal common face of $\tau$ and $\sigma$, and let $j \in I_\tau \setminus I_\zeta$. By the maximality of $\zeta$, the cone $\zeta \vee_\tau \varrho_j$ is not a face of $\sigma$. This means $\eta_\sigma$ does not belong to $S_{\zeta \vee_\tau \varrho_j }$, and so $u_\sigma$ vanishes on $E_{\zeta \vee_\tau \varrho_j}$.
	Notice as well that $ u_\sigma \in \sO_{X, \eta_\zeta}$ and $S_{\zeta \vee_\tau \varrho_j} \subseteq E_{\zeta \vee_\tau \varrho_j}\cap S_\zeta$. It follows that the restriction ${u_\sigma}_{|S_\zeta}$ of $u_\sigma$ to $S_\zeta$ vanishes on $S_{\zeta \vee_\tau \varrho_j}$, and so $u_\sigma$ verifies also (P3). 
\end{proof}

\subsubsection{Proof of the weak approximation theorem} Let $\sigma$ be a face of $\Sigma(X, D)$. Applying Lemma \ref{lem:2}, we find a local equation $z_i$ for $D_i$ around $\eta_\sigma$ for each $i\in I_\sigma$ with the additional property that $z_i\in \mathcal{O}_{X,\tau}$ for each cone $\tau$ that is not a face of $\sigma$. With this choice of local parameters, we  define
\begin{equation}\label{eq:f sigma} 
f_{\sigma}\coloneqq u_\sigma^\ell \sum_{\underline{a}\in A^\sigma}\prod_{i\in I_\sigma} z_i^{a_i},
\end{equation}
for a large enough number $\ell$ which will be precised  in a moment. Notice that $f_\sigma$ is defined in each local ring $\mathcal{O}_{X,\tau}$ for any cone $\tau\in \Sigma(X,D)$. We will prove the following.

\begin{prop} \label{prop:technical}
	Provided $\ell$ is large enough, $f_\sigma$ verifies the following two properties.
	\begin{enumerate}
		\item The set $A^\sigma_{f_\sigma}$ is equal to $A^\sigma$, and
		\item For each face $\tau$ of $\Sigma(X, D)$ different from $\sigma$, we have 
		\[\min_{\leq}(A^\tau_{f_{\sigma}}\cup A^\tau)=A^\tau.\]
	\end{enumerate}
\end{prop}
Using this, we can finish the proof of the approximation theorem.
\begin{proof}[Proof of Theorem~\ref{thm:approx}]
	Let \[f\coloneqq\sum_\sigma \lambda_\sigma f_\sigma\]
	where $\lambda_\sigma$ is a generic choice of coefficients for each face of $\Sigma(X,D)$. Applying the above proposition, we get for each face $\tau$ of $\Sigma(X, D)$, 
	\[A^\tau_f= \min_{\leq} \bigcup_\tau A^\tau_{f_\sigma}= A^\tau.\]
	In other words, $f$ is the rational function stated in the theorem.
\end{proof}

We are left with the proof of Proposition~\ref{prop:technical}.
\begin{proof}[Proof of Proposition~\ref{prop:technical}]
	We use the notations preceding the proposition. By invertibility of $u_\sigma$ in $\sO_{X, \eta_\sigma}$,  the expression \eqref{eq:f sigma} gives an admissible expansion of $f_\sigma$, and so we clearly have $A^\sigma_{f_\sigma} = A^\sigma$. This shows the assertion (1) in the proposition.

	Now, in order to prove Claim (2), let $\tau \not= \sigma$ be a face of $\Sigma(X, D)$, and take local parameters $w_i$ for each $D_i$ around $\eta_\tau$, for $i\in I_\tau$. The element $u_\sigma$ lives in $\sO_{X, \eta_\tau}$ and so in $\widehat{\sO}_{X,\eta_\tau}$. Consider an admissible expansion in $\widehat{\sO}_{X,\eta_\tau}$ for $u_\sigma$
	\begin{equation}\label{eq:adm exp unit}u_\sigma=\sum_{\beta} c_\beta w^\beta.
	\end{equation}
	Property (P2) above implies that for each $j\in I_\tau\setminus I_\sigma$ and for each $\alpha = (\alpha_i)_{i\in I_\tau}$ in the support of \eqref{eq:adm exp unit}, we should have 
	\[\alpha_j\geq \ord_{D_j}(u_\sigma)\geq 1.\]
	
	More generally, we claim the following.

	\begin{claim}\label{claim:bound}
		For each $\alpha$ in the support of the admissible expansion \eqref{eq:adm exp unit}, there is a maximal common face $\zeta$ of $\tau$ and $\sigma$ such that for each $j\in I_\tau \setminus  I_\zeta$, we have $\alpha_j\geq 1$.
	\end{claim}
	\begin{proof} Let $\alpha$ be an element in the support of the admissible expansion \eqref{eq:adm exp unit}, and 
		consider $J=\bigl\{i\in I_\tau\mid \alpha_i=0\bigr\}$. Denote by $\tau_J$  the face of $\tau$ corresponding to $J \subseteq I_\tau$. It will be enough to show that  $\tau_J \subseteq \sigma$. Indeed, in that case, $\tau_J$ will be a common face of $\tau$ and $\sigma$, and so there shall exist a face $\zeta \in \tau \wedge \sigma$ which contains $\tau_J$. For any $j \in I_{\tau} \setminus I_{\zeta}$, we will have $\alpha_j \geq 1$, and the claim will follow.
		
		For the sake of a contradiction, suppose  $\tau_J$ is not a face of $\sigma$, and let $\zeta$ be a maximal common face of $\sigma$ and $\tau_J$. In particular, $\zeta \subsetneq \tau_J$, which implies that $I_\zeta \subsetneq J$. 
		
		We have a projection 
		\begin{align*}\pi\colon \sO_{X,\eta_\tau}&\rightarrow \sO_{S_\zeta,\eta_\tau}\\
		h&\mapsto h_{|S_\zeta}
		\end{align*}
		that extends by continuity to a projection $\pi\colon \widehat{\sO}_{X,\eta_\tau}\rightarrow \widehat{\sO}_{S_\zeta,\eta_\tau}$. Using this, we obtain an admissible expansion for ${u_\sigma}_{|S_\zeta}$ in $\widehat{\sO}_{S_\zeta,\eta_\tau}$ in terms of local parameters ${w_i}_{|S_\zeta}$ for $S_{\zeta \vee_\tau \varrho_i}$, for $i\in I_\tau \setminus I_\zeta$. This is  obtained by applying the projection $\pi$ to both sides of \eqref{eq:adm exp unit}. Indeed, for each $\beta$, the restriction ${c_\beta}_{|S_\zeta}$ is still a unit in $\widehat{\sO}_{S_\zeta,\eta_\tau}$, and so we get 
		\[{u_\sigma}_{|S_\zeta}=\sum_{\beta} c_\beta w^\beta_{|S_\zeta},\]
		which is an admissible expansion in $\widehat{\sO}_{S_\zeta,\eta_\tau}$. In particular, since $I_\zeta \subseteq J$, and since $\alpha_i = 0$ for all $i\in J$, we get that $\pi_{I_\tau \setminus I_\zeta}(\delta)$ is in the support of the admissible expansion for ${u_\sigma}_{|S_\zeta}$.

		Take now $j\in J\setminus I_\zeta$. As ${u_\sigma}_{|S_\zeta}$ vanishes along the divisor $S_{\zeta\vee_\tau \varrho_j}$, and a local equation for this is given by ${w_j}_{|S_\zeta}$, we should have that ${w_j}_{|S_\zeta}$ divides ${u_\sigma}_{|S_\zeta}$ inside $\sO_{S_\zeta,\eta_\tau}$. This implies that for any $\beta$ in the admissible expansion ${u_\sigma}_{|S_\zeta}$, we must have $\beta_j >0$. In particular, this gives $\alpha_j>0$  which contradicts the definition of $J$, so the claim follows.
	\end{proof}
	
	Let $h_\sigma = \sum_{\underline a \in A^\sigma} \prod_{j\in I_\sigma} z_j^{a_j}$. We have $f_\sigma = u_\sigma^\ell h_\sigma$, from which we get the inclusion
	\[A^\tau_{f_\sigma} \subseteq A^{\tau}_{u_\sigma^\ell} + A^{\tau}_{h_\sigma}.\]

	Moreover, we have 
	\[A^\tau_{u_\sigma^\ell} \subseteq \underbrace{A^\tau_{u_\sigma} + \dots + A^{\tau}_{u_{\sigma}}}_{\ell \textrm{ times}}.\]

	Let now $\beta$ be an element of $ A^\tau_{f_\sigma}$. It follows that we can write $\beta$ as the sum of $\ell$ elements in $A^{\tau}_{u_\sigma}$ and an element $\gamma\in A^\tau_{h_\sigma}$.

	By what preceded, we have for each $j\in I_\tau\setminus I_\sigma$ and each element $\alpha$ in $A^{\tau}_{u_\sigma}$ that $\alpha_j \geq 1$. It follows that we have $\beta_j \geq \ell$ for all $j\in I_\tau \setminus I_\sigma$. For $\ell$ large enough, this is certainly larger than the $j$-coordinate of any element in $A^\tau$.

	We now show how to control the $j$-coordinates of $\beta$ for $j\in I_\tau\cap I_\sigma$. For this we write
	\[\beta = \alpha^1+ \dots+\alpha^\ell + \gamma\]
	for $\alpha^1, \dots, \alpha^\ell \in A^{\tau}_{u_\sigma}$ and $\gamma\in A^\tau_{h_\sigma}$. 
	
	Applying Claim~\ref{claim:bound}, for each $\alpha^i$, which is in the support of the admissible expansion  \eqref{eq:adm exp unit} of $u_\sigma$, we infer the existence of a maximal common face $\zeta_i$ of $\tau$ and $\sigma$ such that for each $j \in I_{\tau} \setminus I_{\zeta_i}$, we have $\alpha^i_{j} \geq 1$. Here $\alpha^i_{j}$ is the $j$-coordinate of $\alpha^i$.

	Let $r$ be the number of elements of $\tau \wedge \sigma$. By the pigeonhole principle, there is a maximal common face $\zeta$ of both $\sigma$ and $\tau$ such that we have $\zeta_i =\zeta$ for at least $\ell/r$ indices $i\in [\ell]$.  We thus get for each $j\in I_\tau\setminus I_\zeta$, the inequality
	\[\beta_j=\alpha^1_{j}+\dots + \alpha^\ell_{j} + \gamma_j\geq \ell/r+\gamma_j.\]
	We infer again that, $\ell$ being chosen large enough,  the $j$-coordinate of $\beta$ is larger than the $j$-coordinate of any element in $A^\tau$ provided that $j$ is in $I_\tau\setminus I_\zeta$.

	To finish the proof of Property (2), note that since $\zeta$ is a common face for $\tau$ and $\sigma$, we have by the coherence property that 
	\[\min_{\leq} \, \mathrm{pr}_{_{\hspace{-.05cm} \tau \succ \zeta}}(A^\tau_{f_\sigma})= \min_{\leq}\, \mathrm{pr}_{_{\hspace{-.05cm} \sigma \succ \zeta}}(A^\sigma_{f_\sigma}) \qquad \textrm{and} \qquad \min_{\leq} \mathrm{pr}_{_{\hspace{-.05cm} \tau \succ \zeta}}(A^\tau)=\min_{\leq}\mathrm{pr}_{_{\hspace{-.05cm} \sigma\succ \zeta}}(A^\sigma).\]
	Since $A^\sigma_{f_\sigma} = A^\sigma$, this shows that 
	\[\min_{\leq} \mathrm{pr}_{_{\hspace{-.05cm} \tau \succ \zeta}}(A^\tau_{f_\sigma}) = \min_{\leq} \mathrm{pr}_{_{\hspace{-.05cm} \tau \succ \zeta}}(A^\tau).\]
	
	Hence, given that $\beta \in A^{\tau}_{f_\sigma}$, there is an element $\beta' \in A^\tau$ such that $\mathrm{pr}_{_{\hspace{-.05cm} \tau \succ \zeta}}(\beta)\geq \mathrm{pr}_{_{\hspace{-.05cm} \tau \succ \zeta}}(\beta')$. Moreover, by what we discussed above, choosing $\ell$ large enough, all the $j$-coordinates of $\beta$ for $j$ outside $I_\zeta$ will be also larger than the corresponding $j$-coordinates of $\beta'$.  
	This shows that we actually have $\beta \geq \beta'$ and the claim in (2) follows, that is, 
	\[\min_{\leq}(A^\tau_{f_{\sigma}}\cup A^\tau)=A^\tau.\]
\end{proof}

\subsection{Proof of Corollary~\ref{cor:approx}} \label{sec:cor_approx_proof}
The proof is based on the following proposition 

\begin{prop}  \label{prop:trop_difference} Any tropical function $F$ on the cone complex $\Sigma(X, D)$ can be written as the difference of two tropical functions $F_1$ and $F_2$ such that the restriction of $F_i$ to each cone $\sigma\in \Sigma(X, D)$ is convex and non-negative. 
\end{prop}
\begin{proof}
	Let $\Sigma$ be a subdivision of $\Sigma(X, D)$ such that $F$ is conewise linear on $\Sigma$.	In order to prove the existence of $F_1, F_2$, it will be enough to note that
	\begin{itemize}
		\item[$(i)$] there exists a subdivision $\Sigma'$ of $\Sigma$ and a non-negative tropical function $G$  that is strictly convex on the subdivision given by $\Sigma'$ of each cone of the original dual complex $\Sigma(X, D)$, and 
		\item[$(ii)$] if $G$ is such a function, then for any large enough integer $\ell$, the function $F+\ell G$ is a non-negative convex tropical function. 
	\end{itemize}
	
	In fact, given this, we can write $F = (F+\ell G) - \ell G$  and take $F_1=F+\ell G$ and $F_2=\ell G$ which verify the convexity condition. 
	
	We omit the proof of $(i)$ and only prove $(ii)$.  Let $G$ be a tropical function on $\Sigma(X, D)$ which is conewise linear on $\Sigma'$ and strictly convex on the subdivision of each cone $\sigma$ of $\Sigma(X, D)$ provided by $\Sigma'$. Since $F$ is conewise linear on $\Sigma'$, it follows that for large enough integer $\ell$, $F_1=F+\ell G$ remains strictly convex on each subdivided cone $\sigma$. Since $\ell G$ is obviously convex on the subdivided cones of $\Sigma(X, D)$, the claim follows. 
\end{proof}
\begin{proof}[Proof of Corollary~\ref{cor:approx}] By the proposition above,  there are tropical functions $F_1,F_2$ such that each $F_i$ is non-negative and convex on each facet of $\Sigma(X,D)$ and $F_1-F_2=F$. Each $F_i$, for $i=1,2$, is given by a coherent family of antichains $A_i^\sigma$ in the sense that for each cone $\sigma$ we have
	\[F_i(x)=\min \{\langle x,\beta \rangle \mid \beta \in A^\sigma_i\}.\]
	By approximation theorem, there are rational functions $f_1,f_2$ such that $A_{f_i}^\sigma=A^\sigma_i$ for $i=1,2$ and fo each cone $\sigma$ on $\Sigma(X,D)$. We therefore get $F = \strop(f_1/f_2)$ and the theorem follows.
\end{proof}

\section{Tropical topology on tangent cone bundles}\label{sec:tropical-topology}
In this section we study the \emph{tropical topology} on the space of quasi-monomial valuations. By the duality and approximation theorems proved in the previous sections, this coincides with the coarsest topology on the tangent cone which makes the directional derivatives of tropical functions all continuous. Therefore, we define the tropical topology in the general framework of cone complexes and their tangent cones.

\subsection{Definition of the topology} 
In order to motivate what follows, we first observe that tangent cone bundles inherit a natural Euclidean topology defined as follows.

\begin{defn}[Euclidean Topology] Let $\Sigma$ be a cone complex and $k$ a non-negative integer number. The Euclidean topology on the tangent cone $\TC^k\Sigma$ is the topology defined by the inclusion  \[\TC^k\Sigma\hookrightarrow \bigcup_{\sigma \in \Sigma} N_{\sigma,\mathbb{R}}^{k}.\]
Here, the space on the right hand side is obtained by gluing the vector spaces $N_{\sigma,\mathbb{R}}^{k}$. It is endowed with the quotient topology induced by the Euclidean topology on $N^k_{\sigma,\mathbb{R}}$, $\sigma \in \Sigma$.
\end{defn}

This topology turns out however to be not properly adapted to the higher rank context. This is suggested by  the following example of  a tropical function whose derivative is not continuous with respect to the Euclidean topology. 

\begin{example}\rm Let $\sigma=\R_+\times \R_+$ and consider the tropical function 
	\begin{align*} F\colon \sigma&\longrightarrow \R \\
	(x_1,x_2)&\longmapsto \min\{x_1,x_2\}.
	\end{align*}
	For the first directional derivative of $F$, we have
	\begin{align*} 
	\Der F\colon \TC\sigma&\longrightarrow \R^2 \\
	((x_1,x_2);(y_1,y_2))&\longmapsto 
	\begin{cases} 
	(x_1,y_1) & \text{ if either, } x_1<x_2, \text{ or, } x_1=x_2 \text{ and } y_1<y_2\\
	(x_2,y_2) & \text{ if either, } x_1>x_2 \text{ or, } x_1=x_2 \text{ and } y_2<y_1.
	\end{cases} 
	\end{align*}
	This map is not continuous with respect to the Euclidean topology. To see this, consider the map  \[t\mapsto \Der F((t,1-t);(y_1,y_2))\] for $y_1\neq y_2$. This function has a discontinuity at $t=\frac{1}{2}$.
\end{example}

The analytic description of higher rank quasi-monomial valuations leads the following natural topology.

\begin{defn}[Tropical topology] Let $\Sigma$ be a cone complex and $\TC^k\Sigma$ its tangent cone. We consider $\R^{k+1}$ with its Euclidean topology and define the \emph{tropical topology} on $\TC^k\Sigma$ as the coarsest topology which makes all the maps
	\[\Der^kF\colon \TC^k\Sigma\longrightarrow \R^{k+1}\]
	continuous for any tropical function $F\colon \supp{\Sigma}\rightarrow \R$.
\end{defn}

\begin{remark} In the case $\Sigma=\Sigma(X, D)$ for a smooth quasi-projective variety $X$ and an SNC divisor $D$ on $X$, the tropical topology on $\TC^{k}\Sigma$ is the coarsest one making the directional derivative $\Der^{k} \strop(f)$ a continuous function from $\TC^k\Sigma \to \R^{k+1}$, for any $f \in K(X)^\times$.
\end{remark}

\subsection{Description of the topology}

The aim of this section is to give a description of this topology by introducing a basis of open sets. This will be based on the following definition. 

\begin{defn}[$\widetilde\Sigma$-open sets] Let $\Sigma$ be a cone complex and $\widetilde{\Sigma}$ a rational subdivision of it. A set $U\subset \TC^k\Sigma$ is called a \emph{$\widetilde \Sigma$-open set} if for any cone $\sigma \in \widetilde\Sigma$, the intersection $U\cap \TC^k\sigma$ is open in $\TC^k\sigma$ with respect to the Euclidean topology.
\end{defn}

\begin{figure}[h]
		\centering
		\includegraphics[width=1\linewidth]{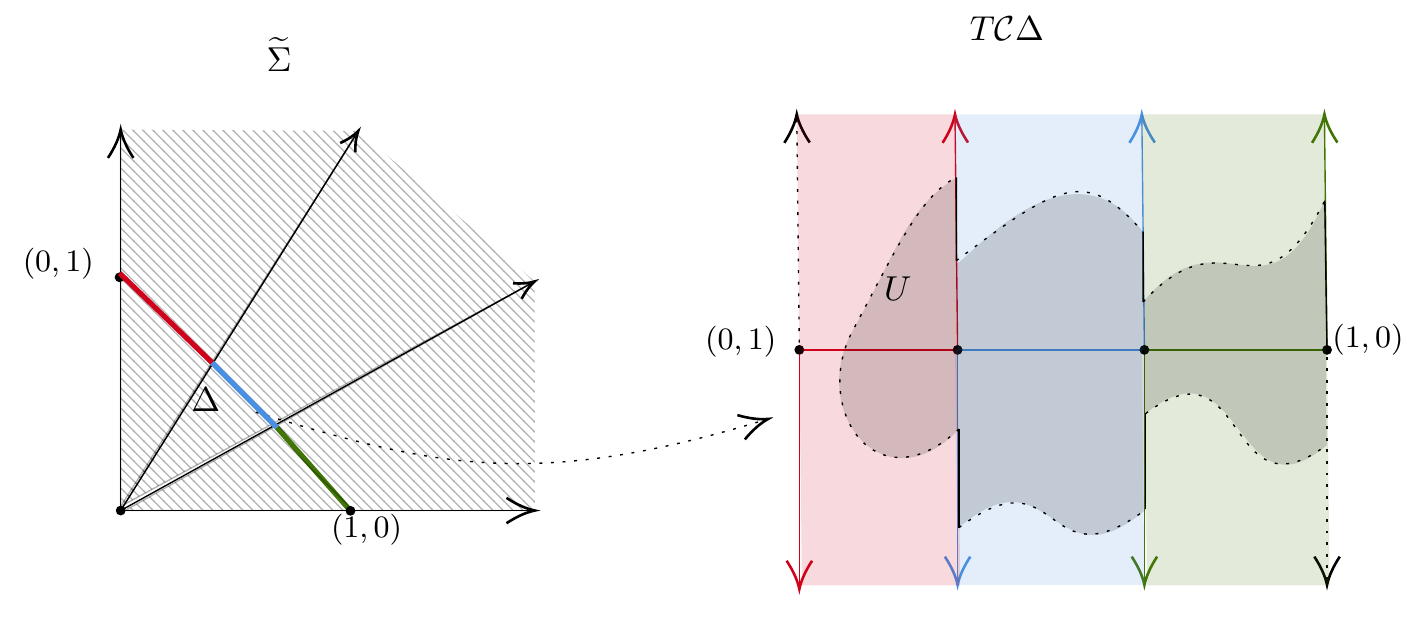}
		\caption{The restriction of a $\widetilde{\Sigma}$-open set $U$ to the set $\TC \Delta$.}
		\label{fig:sigma-open}
	\end{figure}

\begin{example} In order to give an idea of the notion of  $\widetilde{\Sigma}$-open set, we consider the example depicted in Figure~\ref{fig:sigma-open}. We have the positive orthant $\sigma=\bigl\{(x,y)\in \R^2\,\bigl| \, x\geq 0, y\geq 0\bigr\}$ together with a subdivision $\widetilde{\Sigma}$ of it. The section $\Delta=\operatorname{conv}((0,1),(1,0))$ gives a segment subdivided in three smaller segments colored in red, blue and green. The tangent cone $\TC \Delta$ can be written as the union of the tangent cones of each of these smaller segments. Drawing the fibers vertically, we get the picture on the right hand side of Figure \ref{fig:sigma-open}. A $\widetilde{\Sigma}$-open set restricted to $\TC\Delta$ looks like the union of open sets one in the induced topology of each of the red, blue and green parts. Note in particular that it is not necessarily an open set in the Euclidean topology of $\TC\Delta$.

\end{example}

We have the following characterization of the tropical topology. 
\begin{thm} \label{thm:characterization-tropical-topology} Let $\Sigma$ be a cone complex and consider its tangent cone $\TC^k\Sigma$ for $k$ a non-negative integer number. Then,
	\begin{enumerate}
		\item For each subdivision $\widetilde\Sigma$ of $\Sigma$, the $\widetilde\Sigma$-open sets of $\TC^k\Sigma$ are open with respect to the tropical topology. 
		\item The union of all $\widetilde\Sigma$-open sets for $\widetilde \Sigma$ a rational subdivision of $\Sigma$ forms a basis of open sets for the tropical topology.
	\end{enumerate}
\end{thm}

The proof of this theorem is given in the next section. We state the following corollary.

	\begin{cor} Let $\Sigma$ be a cone complex and $\TC^k\Sigma$ its tangent cone endowed with the tropical topology. Then,
		\begin{enumerate}
			\item The tropical topology is second countable and is finer than the Euclidean topology. In particular, it is both Hausdorff and normal.

			\item A set is dense in $\TC^k\Sigma$ with respect to the tropical topology if and only if it is dense with respect to the Euclidean topology.

			\item $\TC^k\Sigma$  is not locally compact in general \emph{(}in fact, as soon as $k>0$ and the dimension of $\Sigma$ is at least two\emph{)}.
		\end{enumerate}
	\end{cor}

	\begin{proof} Let us start by (1). Given a countable basis $\mathcal{U}$ of $\R^{k+1}$, we can construct a countable subbasis of $\TC^{k}\Sigma$ by considering all the open sets of the form $\Der F^{-1}(U)$ were $F$ is a tropical function and $U\in \mathcal{U}$. This shows that the tropical topology is second countable. Moreover, as any Euclidean open set is $\widetilde\Sigma$-open for any subdivision $\widetilde\Sigma$ of $\Sigma$, we get that the tropical topology is finer than the Euclidean topology. Hence, it is Hausdorff and normal.

		For point (2), it is enough to notice that for any subdivision $\widetilde\Sigma$ of $\Sigma$ and any $\widetilde\Sigma$-open set $U$, there is an Euclidean open set $V$ contained in $U$.

		For point (3), arguing by contradiction, suppose $U$ is open and the closure $\overline{U}$ in the tropical topology is compact. We then find an open set $V\subseteq U$ such that $V$ is open in the Euclidean topology, as in (2). The closure $\overline{V}$ with respect to the tropical topology is then compact. Denoting by $\overline{V}_{\mathrm{Trop}}$ and $\overline{V}_{\mathrm{Euc}}$ the set $\overline{V}$ endowed with its tropical and Euclidean topologies, respectively, we see that the identity map
		\[\mathrm{id}\colon \overline{V}_{\mathrm{Trop}}\rightarrow \overline{V}_{\mathrm{Euc}}\]
		is continuous. Moreover, as $\overline{V}_{\mathrm{Trop}}$ is compact and $\overline{V}_{\mathrm{Euc}}$ is Hausdorff, the identity becomes then a homeomorphism. This implies that the tropical and the Euclidean topologies should agree on $V$. However, this is not possible if we have both $k>0$ and $\dim \Sigma>1$ (we can choose a subdivision $\widetilde\Sigma$ of $\Sigma$ subdividing $V$, and a $\widetilde \Sigma$-open set $W\subseteq V$ such that $W\cap \TC\sigma$ is $\widetilde\Sigma$-open but not Euclidean open).
	\end{proof}

\subsection{Proof of Theorem~\ref{thm:characterization-tropical-topology}}

We adapt the following terminology in the sequel. For a cone $\sigma$ of a cone complex $\Sigma$, and a point $(x;\underline w) \in \TC^k \Sigma$, by saying \emph{$\sigma$ supports the point $(x; \underline w)$} we mean {the point $(x;\underline w)$ belongs to $\TC^k\sigma$.

	 We first prove the second point assuming the first.

	\begin{proof}[Proof of (2)]	Let $F\colon \supp{\Sigma} \rightarrow \R$ be a tropical function and consider a rational subdivision $\widetilde \Sigma$  of $\Sigma$ such that $F$ is linear on each cone of $\widetilde \Sigma$. Let $V \subset \R^{k+1}$ be an open set for the Euclidean topology. We show that $\bigl(\Der^k F\bigr)^{-1}(V)$ is $\widetilde \Sigma$-open. This proves the result.

		Let $\delta$ be a cone of $\widetilde \Sigma$. By the choice of $\widetilde \Sigma$, there is a linear function $F_\delta$ on $N_{\delta, \R}$ such that for any point $(x;\underline w)\in \TC^k\delta$, with $\underline w = (w_1, \dots, w_k)$, we have
		\[\Der^k F(x; \underline{w})=(F_\delta(x),F_\delta(w_1), \dots,F_\delta(w_k)).\]
		The intersection 
		\[\bigl(\Der^k F\bigr)^{-1}(V)\cap \TC^k \delta= \bigl(\underbrace{F_\delta \times \dots \times F_{\delta}}_{(k+1) \textrm{ times}}\bigr)^{-1}(V)\cap \TC^k\delta\] 
		is clearly an open set in  $\TC^k \delta$ for the Euclidean topology, and the claim follows.
		
			\end{proof}

		\begin{proof}[Proof of (1)] Let $\widetilde \Sigma$ be a rational subdivision of $\Sigma$ and let $U$ be a $\widetilde{\Sigma}$-open set. We have to show that $U$ is open for the tropical topology of $\TC^k\Sigma$.

		We first observe that if $\widetilde \Sigma'$ is a subdivision of $\widetilde \Sigma$, any $\widetilde \Sigma$-open set is also $\widetilde \Sigma'$-open. Therefore, in order to prove the above claim, we can assume that $\widetilde \Sigma$ is simplicial.

		Take $(x;\underline{w})\in U$. We will prove that $(x;\underline{w})$ is an interior point of $U$ for the tropical topology by explicitly constructing a neighborhood of $(x;\underline w)$ for the tropical topology included in $U$.
		
		Let $\zeta$ be the minimal face of $\widetilde\Sigma$ which supports $(x;\underline w)$. For each facet $\delta$ of $\widetilde \Sigma$ we find a rational subdivision $\widetilde \Sigma_\delta$ of $\widetilde \Sigma$ with the following properties.

		\begin{enumerate}[label=(\alph*)]
			\item $\widetilde \Sigma_\delta$ is simplicial.
			 
			\item There is a unique facet in $\widetilde \Sigma$ denoted by $\gamma=\gamma_{\zeta, \delta}$ which contains $\zeta$ and which is contained in $\delta$.
			 
			\item  For each pair $(\delta, \varrho)$ consisting of $\delta$ and a ray $\varrho$ of $\delta$, there is a tropical function $F^{\delta,\varrho}$ on $\supp{\Sigma}$ such that...			
			
			\begin{itemize}
				\item[(1)] $F^{(\delta,\varrho)}$ is linear on each cone of $\widetilde \Sigma_\delta$.
				
				\item[(2)] Over the facet $\gamma$ of $\widetilde \Sigma_\delta$, we have $F^{\delta,\varrho}\rest\gamma=\chi^m\rest\gamma$ where $m$ is the primitive element in the ray dual to $\varrho$ in $\delta^\vee$ and $\chi^m$ is the linear function induced by this vector. That is, $F^{\delta,\varrho}$ takes value one on the primitive vector of $\varrho$ and value zero on all the rays of $\delta$.
			\end{itemize}
		\end{enumerate}
		
		For such a facet $\delta$ of $\widetilde \Sigma$ which contains $\zeta$,  consider the function $\Phi_\delta$  defined by the collection of functions $F^{\delta, \varrho}$, $\varrho$ a ray of $\delta$,
		\[\Phi_\delta \coloneqq (F^{\delta,\varrho})_{\varrho \text{ ray of }\delta}\colon\, \Sigma \rightarrow \R^{\dim(\delta)}. \]
		Let $\gamma = \gamma_{\zeta, \delta}$ be the facet of $\widetilde \Sigma(\zeta)$ which contains $\zeta$ and which is contained in $\delta$. 
		By  (c), the linear functions $F^{\delta, \varrho}$ for $\varrho$ a ray of $\delta$ are linearly independent on $\gamma$, and so $\Phi_\delta$ restricted to $\gamma$ is a homeomorphism to its image in $\R^{\dim(\delta)}$.  
		
		The directional derivative map
		\begin{equation}\label{map DF} \Psi_\delta=(\Der^k F^{\delta,\varrho})_{\varrho \text{ ray of }\delta}\colon  \TC^k\Sigma\rightarrow (\R^{\dim(\delta)})^{k+1}
		\end{equation}
		is a homeomorphism onto its image when restricted to $\TC^k \gamma$, when we put on $\TC^k\gamma$ its Euclidean topology. Indeed, restricted to $\TC^k \gamma$, $ \Psi_\delta$ can be identified with the restriction to $\TC^k\gamma$ of the invertible linear map 
		\[ (\underbrace{\Phi_\delta \times \Phi_\delta \times \dots \times \Phi_\delta}_{(k+1) \textrm{ times }}) \colon  \R^{\dim(\delta) \times (k+1)} \to \R^{\dim(\delta) \times (k+1)}.\]
		
		Hence, there is an open set $U_\delta \subseteq (\R^{k+1})^{\dim \delta}$ such that its preimage under the map \eqref{map DF} satisfies
		\[\Psi_\delta^{-1}(U_\delta)\cap \TC^k \gamma = U\cap \TC^k\gamma.\]
		Note that 
		\[\Psi_\delta^{-1}(U_\delta)  = \bigcap_{\delta, \varrho \textrm{ ray in } \delta} \bigl(\Der^kF^{(\delta, \varrho)}\bigr)^{-1}(U_\delta),\]
		and so $\Psi_\delta^{-1}(U_\delta)$ is an open set in the tropical topology and therefore a neighborhood of $(x;\underline{w})\in U$. This proves that $(x;\underline{w})$ is an interior point of the intersection $\bigcap_{\delta} \Psi_\delta^{-1}(U_\delta)$ for $\delta$ running over all facets which contain $\zeta$, which is an open set for the tropical topology. Denote by $W$ this intersection. Note that we have $W \cap \TC^k\gamma_{\zeta, \delta} \subset U \cap \TC^k\gamma_{\zeta, \delta}$ for each facet $\delta$ which contains $\delta$.

		Let now $\widetilde \Sigma'$ be a rational subdivision of $\Sigma$ with the following properties:
		\begin{itemize}
			\item $\widetilde \Sigma'$ is finer than $\widetilde \Sigma_\delta$ for all facets $\delta$ which contain $\zeta$.
			\item there exists a tropical function $G$ which is linear on each cone of $\widetilde \Sigma'$, and 
			which is strictly positive on the relative interior of each cone $\tau$ of $\widetilde \Sigma'$ which supports $(x; \underline w)$ and which is non-positive everywhere else.  
		\end{itemize}
		
		Then, $(\Der^kG)^{-1}(\R_{>0}\times \R^k)$ is an open set in the tropical  topology, and moreover, it is contained in the union $\bigcup_\tau \TC^k \tau$ where the union goes over all cones $\tau$ of $\widetilde \Sigma'$ which support $(x; \underline{w})$. It follows from these constructions that 
		$$(x; \underline{w})\in (\Der^kG)^{-1}(\R_{>0}\times \R^k)\cap \bigcap_\delta \Psi_\delta^{-1}(U_\delta) \subseteq \bigcup_\delta U \cap \TC^k\delta$$
		where $\delta$ runs over facets of $\widetilde \Sigma$ which contain $(x; \underline w)$. We infer that $U$ is a neighborhood of $(x; \underline{w})$ in the tropical topology, and the theorem follows.
	\end{proof}

	\section{Spaces of valuations and the retraction map} \label{sec:spaces-valuations}
	For a given variety $X$, we introduce several spaces of valuations and show how tangent cones endowed with their tropical topology naturally fit inside them. 
	
	\subsection{Higher rank analytification and its centroidal filtration}
	\begin{defn} Given a variety $X$, we define the \emph{birational analytification of $X$ of rank bounded by $k$} as the set
		\[X^{\bir,k}\coloneqq\bigl\{\nu\colon K(X)^*\rightarrow \R^k\st  \nu \text{ is a valuation}\bigr \}\]
	 endowed with the coarsest topology which makes continuous all the evaluation maps
		\begin{align*}
		\mathrm{ev}_f\colon X^{\bir,k}&\longrightarrow \R^k \\
		\nu &\longmapsto \nu(f),
		\end{align*}
		for any $f\in K(X)^*$, $\R^k$ given its Euclidean topology. We define moreover the subspaces 
		\begin{align*}X^{\bethd,k}\coloneqq\bigl\{\nu\in X^{\bir,k}\st \nu \text{ has center in } X\bigr \} \\ 
		X^{\dalethd,k} \coloneqq \bigl\{\nu \in X^{\bir,k}\st \nu \text{ does  not have center in } X\bigr\}
		\end{align*}
		of $X^{\bir, k}$, and endow them with the topology induced by that of $X^{\bir, k}$.	\end{defn}
	
	\begin{remark}  Notice that $X^{\bir,k}=X^{\bethd,k}\sqcup X^{\dalethd,k}$ and $X^{\bir,k}=X^{\bethd,k}$ if $X$ is proper. In the terminology of~\cite{FR16}, the space $X^{\bir,k}$ coincides with the subspace of all valuations defined over the generic point in the Hahn analytification of $X$ endowed with the extended Euclidean topology. Moreover, the notation $X^{\bethd,k}$ is used in analogy with the analytic space $X^\bethsm$ of Berkovich~\cite{Ber-formal} and Thuillier~\cite{Th}, where we have used a dot as a reminder that we are considering only the birational parts.
	\end{remark}
	
	We now introduce a flag of subspaces on $X^{\bir,k}$ which interpolate between $X^{\bethd,k}$ and $X^{\bir,k}$. Most of the constructions we will do in the following will be compatible or can be extended to this centroidal filtration.

	\begin{defn}[The centroidal filtration] We define $\mathscr{F}^0X^{\bir,k}:= X^{\bir,k}$, and for $1\leq r\leq k$, we define
		\[\mathscr{F}^rX^{\bir,k}:=\bigl\{\nu\in X^{\bir,k}\st \mathrm{proj}_r(\nu) \text{ has center in }X\bigr\}\]
		where $\mathrm{proj}_r(v)$ is the composition of $v$ with the projection $\R^k\rightarrow \R^r$ to the first $r$ coordinates. In other words, 
		\[\mathscr{F}^rX^{\bir,k}=\mathrm{proj}_r^{-1}X^{\bethd,r}.\]
		This gives a decreasing filtration
		\[X^{\bir,k}=\mathscr{F}^0X^{\bir,k}\supseteq \mathscr{F}^1X^{\bir,k} \supseteq \dots \supseteq \mathscr{F}^kX^{\bir,k}=X^{\bethd,k}.\]
	\end{defn}

	\subsection{Inclusion of tangent cones in the analytification}	
		Given an SNC divisor $D$ on the variety $X$, by Theorem~\ref{thm:analytic} we get an inclusion map
		\begin{align*} \iota\colon \TC^{k-1}\Sigma(D)&\lhook\joinrel\longrightarrow X^{\bethd,k}\\
		(x;\underline{w})&\longmapsto \nu_{x;\underline{w}}.
		\end{align*}
		
	\begin{prop}[The inclusion map]\label{prop:inclusion} The above map $\iota$ induces a homeomorphism between $\TC^{k-1}\Sigma(D)$ endowed with the tropical topology and its image with the  topology induced by $X^{\bethd,k}$. In the case $X$ is proper, we can restrict the codomain to an inclusion
		\begin{align*} \iota\colon \TC^{k-1}\Sigma(D)&\lhook\joinrel\longrightarrow (X\setminus D)^{\dalethd,k}.
		\end{align*}
	\end{prop}
	\begin{proof} By Proposition \ref{thm:monomial} the center of a quasi-monomial valuation defined by $D$ is always on $D$, so we can restrict to the codomain above in each case. Moreover, the fact that the map induces a homeomorphism with its image is a direct consequence of the approximation theorem. 
	\end{proof}
	
	Regarding the center, the following result will be useful later. For a valuation $\nu$ in $X^{\bethd,k}$, we denote by $\centre{\nu}$ the center of $\nu$ in $X$.

	\begin{prop}\label{prop:anticontinuous} Let $X$ be a variety, then the map
		\begin{align*} c_X\colon  X^{\bethd,k}&\longrightarrow X\\
		\nu&\longmapsto \centre{\nu}
		\end{align*}
		that assigns to each valuation its center in $X$ is anticontinuous.
	\end{prop}
	\begin{proof} The proof will use the Topology-Mixing Lemma~\ref{lem:mixing-lemma} proved below. Let $U=\Spec(A)\subseteq X$ be an affine open set. Then for a valuation $\nu\in X^{\bethd,k}$, we have
		\begin{align*}
		\centre{\nu}\in U \iff \nu_{\rest{A}}\geq 0 \iff \nu\in \bigcap_{\substack{f\in A\\f\neq 0}} \ev{f}^{-1}[0,\infty).
		\end{align*}
		
		By Lemma~\ref{lem:mixing-lemma} (applied to $1/f$), the set $\ev{f}^{-1}[0,\infty)$ is closed. We infer that $c_X^{-1}(U)=\bigcap_{f\in A} \mathrm{ev}_f^{-1}[0,\infty)$ is closed. If $V$ is an arbitrary open set, then, since $X$ is Noetherian, we have a finite cover $V=\bigcup_i U_i$ by open affine subsets and so $c_X^{-1}(V)=\bigcup_i c_X^{-1}(U_i)$ is closed.
	\end{proof}

	\subsection{The retraction map}\label{sec:retraction}
	Let $X$ be a smooth variety and $D$ an SNC divisor on $X$. Endowing the tangent cone with the tropical topology, Proposition~\ref{prop:inclusion} gives an inclusion of $\TC^{k-1}\Sigma(D)$ as a topological subspace of $X^{\bethd,k}$. In this section we construct a retraction of $X^{\bethd,k}$ onto $\TC^{k-1}\Sigma(D)$ for this inclusion, and study its basic properties. This generalizes the picture from rank one to higher rank. 
	
	\subsubsection{Definition of the retraction map}
	
	We start by recalling how to apply a valuation to a divisor when the valuation has a center in a variety $X$.
	
	\begin{defn} Let $E$ be a Cartier divisor in  $X$. Given a valuation $\nu\in X^{\bethd,k}$ with center $\centre{\nu}$ in $ X$, we define $\nu(E)\coloneqq \nu(z)$ where $z\in \mathcal{O}_{X,\centre{\nu}}$ is a local equation for $E$ around the point $\centre{\nu}$. 
	\end{defn}
	
	As two local equations differ by a unit, this is well-defined.
	We identify $\M^k(D)$ with $\TC^{k-1}\Sigma(D)$ using the duality Theorem.
	
	\begin{defn}[Retraction]\label{defn:retraction}  Let $D$ be an SNC divisor on a variety $X$. \emph{The retraction to $\TC^{k-1}\Sigma(D)$} is the map
		\begin{equation}\label{eq:retraction} \ret{D}\colon X^{\bethd,k}\rightarrow \TC^{k-1}\Sigma(D)
		\end{equation}
		given  by sending any valuation $\nu\in X^{\bethd,k}$ to the unique pair $(x;\underline w)\in \TC^{k-1}\Sigma(D)$, with corresponding quasi-monomial valuation $\nu_{x, \underline w}$, which verifies for any component $D_i$ of $D$, the equality 
		\begin{equation}\label{retraction defn} \nu_{x,\underline w}(D_i) =\nu(D_i).
		\end{equation}
	\end{defn} 
	\begin{prop}\label{prop:retraction} The map $\ret{D}$ verifies the following properties: 
		\begin{enumerate}
			\item It is well-defined.
			\item It is continuous.
			\item It is a retraction for the inclusion $\iota$ from Proposition \ref{prop:inclusion}, that is, $\ret{D}\circ \iota=\id$. 
		\end{enumerate}
	\end{prop}	
	\begin{proof} 	(1) and (3) are clear. The proof of (2) will be based on the Topology-Mixing Lemma~\ref{lem:mixing-lemma} below, and will be given in Section~\ref{sec:cont-retraction}.
	\end{proof}

	\subsection{Topology-Mixing Lemma}  We prove the following lemma.

	\begin{lem}[Topology-Mixing Lemma]\label{lem:mixing-lemma} Let $X$ be an algebraic variety and fix an element $f\in K(X)^*$. The set 
		\[\ev{f}^{-1}\bigl((-\infty,0]\bigr)=\bigl\{\nu \in X^{\bethd,k}\mid \nu(f)\leqlex 0\bigr\}\]
		is a closed set inside $X^{\bethd,k}$.
	\end{lem}
	
	\begin{remark} Notice that the interval $(-\infty,0]$ in $\R^k$ constructed with the lexicographic order is not closed with respect to the Euclidean topology. Therefore, the lemma does not follow directly from the definition of the tropical topology on $X^{\bethd, k}$. It might appear to be somehow unexpected in this regard, as it happens to mix the Euclidean and ordered topologies (where the name given to the result). The statement might not be true when the interval $(-\infty,0]$ is replaced by other half intervals. 
	\end{remark}
	
	\begin{proof}	This is equivalent to showing that \[\ev{f}^{-1}\bigl((0,\infty)\bigr)=\bigl\{\,\nu \in X^{\bethd,k}\, \bigl|\, \nu(f)\succ_{\lex}0\, \bigr\}\]
		is an open set inside $X^{\bethd,k}$. For this, we will show that any element $\nu\in \ev{f}^{-1}(0,\infty)$ is an interior point of $\ev{f}^{-1}(0,\infty)$. 
		
		So let $\nu$ be such an element. First, notice that since $\nu(f)\succ 0$, we must have $\nu(f+1)=0$. This shows $\nu(f) \neq \nu(f+1)$. 
		Take now two disjoint open neighborhoods $U$ and $V$ of $\nu(f)$ and $\nu(f+1) =0$ in $\R^k$,  respectively, so that we have $\nu \in \ev{f}^{-1}(U)\cap\ev{f+1}^{-1}(V)$. For any valuation $\tilde{\nu}\in \ev{f}^{-1}(U)\cap\ev{f+1}^{-1}(V)$, we have $\tilde{\nu}(f)\neq \tilde{\nu}(f+1)$. This implies that
		\[0=\tilde{\nu}(1)=\min_{\leqlex}\bigl\{\tilde{\nu}(f+1),\tilde{\nu}(f)\bigr\}.\]
		Since $\tilde{\nu}(f)\in U$, and $0 =\nu(f+1)\notin U$ (by the choice of $U$ and $V$), we get $\tilde{\nu}(f)\neq 0$. This implies that $\tilde{\nu}(f)\succ 0$ and so $\tilde \nu \in \ev{f}^{-1}\bigl((0,\infty)\bigr)$. We infer that the open neighborhood $\ev{f}^{-1}(U)\,\cap\,\ev{f+1}^{-1}(V)$ of $\nu$ in $X^{\bethd,k}$ is contained in $\ev{f}^{-1}\bigl((0,\infty)\bigr)$, from which the result follows.
	\end{proof}
	
\subsection{Continuity of the retraction map}\label{sec:cont-retraction} In this section, we prove 	
	part (2) of Proposition~\ref{prop:retraction} by using the Topology-Mixing Lemma.

	By the definition of the tropical topology, in order to prove the continuity of $\ret{D}$, it will be enough to show that for each tropical function $F\colon \Sigma(D)\rightarrow \R$, the composition \[\mathcal{F}\coloneqq\Der^{k-1}F \circ \ret{{D}}\colon X^{\bethd,k}\rightarrow \R^k\] is continuous. Using the approximation theorem, we can find a rational function $f$ such that $F=\strop(f)\colon \Sigma(D)\rightarrow \R$ and then $\mathcal{F}=\mathscr{D}^{k-1}\strop(f)\circ \ret{{D}}$. We will fix such a function. 
	
We will construct a sequence of covers of the form $X^{\bethd,k}=\bigcup_i G_i$ by finitely many closed subsets, and reduce to showing the continuity of the restriction $\mathcal F \rest{G_i}$ to each $G_i$.

	We start by taking a finite affine open cover $X=\bigcup_j U_j$ with the property that each component $D_i$ of $D$ is a principal divisor over each $U_j$. We then have $X^{\bethd,k}=\bigcup_i U_i^{\bethd,k}$ where \[U_i^{\bethd,k}:=\bigl\{\, \nu\in X^{\bethd,k}\, \bigl|\, c_\nu\in U_i\, \bigr\}\] is closed by Proposition~\ref{prop:anticontinuous}. We thus get a finite closed cover. By the observation above, it will be enough to prove that $\mathcal{F}$ is continuous restricted to each $U^{\bethd,k}_i$. This means, we can assume that $X$ is affine and each divisor $D_i$ is principal and so of the form $\div(z_i)=D_i$ for some regular function $z_i$. Moreover, we can assume as well that the rational function $f$ defining $\mathcal F$ is a regular function on $X$.

	Now, for each cone $\sigma\in \Sigma(D)$, consider the set
	\[G_\sigma\coloneqq\bigl\{\,\nu\in X^{\bethd,k}\, \bigl|\, \centre{\nu}\in X\setminus \bigcup_{j\notin I_\sigma} D_j\bigr\}.\] 
	That is, $G_\sigma$ is the set of all valuations $\nu \in X^{\bethd,k}$ whose center $\centre{\nu}$ does not belong to any component $D_i$ with $i\notin I_\sigma$. As we have $G_\sigma = c_X^{-1}\bigl(\,X\setminus \bigcup_{i\notin I_\sigma} D_i\,\bigr)$, by Proposition~\ref{prop:anticontinuous} this set is a closed subset of $X^{\bir,k}$. 
	
	Given $\nu\in X^{\bethd,k}$, either, $\centre{\nu}\notin D$, in which case we get $\nu\in G_\sigma$ for all facets $\sigma$ of $\Sigma(D)$. Or, $\centre{\nu}\in D$, in which case, there exists a component $D_i$ such that  $\centre{\nu}\in D_i$. Taking $\sigma$ the face of $\Sigma(D)$ whose associated stratum $S_\sigma$ contains $\centre{\nu}$, we get that $\nu\in G_\sigma$. This shows that 
	
\begin{prop} The family $G_\sigma$, $\sigma$ a face of $\Sigma(D)$, forms a closed cover of $X^{\bethd,k}$. 
	\end{prop}
	
	Hence, it is enough to prove the continuity over each $G_\sigma$. Without loss of generality, assume that $D_1, \dots, D_r$ are all the components of $D$ which contain $\eta_\sigma$.
	
	As $f$ is regular, we have $f\in \sO_{\eta_\sigma}$. Consider now an admissible expansion 
	\[f=\sum_\beta c_\beta z^\beta\in \mathcal{O}_{X,\eta_\sigma}\]
	of $f$ around the point $\eta_\sigma$ where $\div(z_i)=D_i$.  
	Then, $\mathcal{F}(\nu)=\min\bigl\{\, \nu(z^\beta)\,\bigl|\, \beta\in \minS{f}^\sigma\,\bigl\}$ for each $\nu\in G_\sigma$. In order to prove the continuity of $\mathcal{F}$ on $G_\sigma$, we further decompose $G_\sigma$ as a finite union of closed sets as follows. For each $\beta\in \minS{f}^\sigma$, consider the set 
	\[G_{\sigma, \beta}\coloneqq\bigl\{\, \nu\in G_\sigma\, \bigl|\, \mathcal{F}(\nu)=\nu(z^\beta)\,\bigr\}.\]
	 We have 
		\begin{align*}
		G_{\sigma, \beta} = \bigl\{\,\nu\in G_\sigma\, \bigl|\, \nu(z^\beta)\leqlex \nu(z^{\beta'}) \; \forall \beta'\in \minS{f}^\sigma\,\bigr\}
		= G_\sigma \cap \bigcap_{\beta'\in \minS{f}^\sigma} \bigl\{\, \nu\in G_\sigma\, \bigl|\, \nu(z^{\beta-\beta'})\leqlex 0\, \bigr\}
		\end{align*}
		which is closed by the Topology-Mixing Lemma applied to rational functions $z^{\beta -\beta'}$. We get

	\begin{prop} The family $G_{\sigma, \beta}$, $\beta \in \minS{f}^\sigma$, is a closed cover of $G_\sigma$.
	\end{prop}

	We can now finish the proof of the continuity of the retraction map.
	
	\begin{proof}[Proof of part (2) of Proposition \ref{prop:retraction}] By the above discussion, we are reduced to show that $\mathcal{F}$ is continuous over each $G_{\sigma, \beta}$. But this is clear since the restriction of $\mathcal{F}$ to $G_{\sigma, \beta}$ equals $\ev{z^\beta}$, which is continuous by the definition.
	\end{proof}

	\section{Log-smooth pairs}\label{sec:log-smooth}
	In the previous section, given a variety $X$ we constructed a retraction from $X^{\bethd,k}$ to the tangent cone $\TC^{k-1}\Sigma(D)$ associated to an SNC divisor $D$ on $X$. In this section,  we will introduce other instances   for which we can construct dual complexes, tangent cones, and corresponding retractions. The results will be of use in the subsequent section in order to prove the limit formulae.  We start by the following definition. 
	
	\begin{defn} Let $X$ be a smooth variety
		\begin{enumerate}
			\item A \emph{log-smooth pair over $X$} is the data of a pair $\mathbf{Y}=(Y,D)$ consisting of a smooth variety $Y$ and an SNC divisor $D$ on $Y$ together with a proper morphism $\varphi\colon Y\rightarrow X$ such that the restriction \[\varphi\rest{Y\setminus D}\colon Y\setminus D\longrightarrow X \setminus \varphi(D)\]
			is an isomorphism.	The morphism $\varphi$ is called the \emph{structure morphism of the log-smooth pair.}
			
			Given log-smooth pairs $\mathbf{Y}'=(Y',D')$ and $\mathbf{Y}=(Y,D)$, a morphism $\mathbf{Y}'\rightarrow \mathbf{Y}$ between them is a proper morphism 
			\[f\colon Y'\longrightarrow Y\]
			that commutes with the structure map of $Y$ and $Y'$ and such that $\Supp\bigl(f^*(D)\bigr)\subseteq \Supp \bigl(D'\bigr)$. 
			
			We denote by $\LSP_X$ the category of log-smooth pairs over $X$.

			\item A \emph{log-smooth compactification of $X$} is a proper variety $Y$ containing $X$ as an open subvariety such that $Y\setminus X$ is an SNC divisor on $Y$. 
			
			A morphism between log-smooth compactifications $Y'$ and $Y$ is a morphism $f\colon Y'\rightarrow Y$ between the underlying varieties such that $f^{-1}(X)=X$ and $f\rest{X}$ is an isomorphism. The category of log-smooth compactifications of $X$ will be denoted by $\LSC_X$. 
			
			(Notice that for a morphism as above we have $f^{-1}(Y\setminus X)=Y'\setminus X$.)

			\item A \emph{compactified log-smooth pair} is the data of a pair $\overline{\mathbf{Y}}=(Y,D)$ consisting of a proper variety $Y$ and an SNC divisor $D\subset Y$ together with a birational map $\varphi\colon Y\dashrightarrow X$  such that
			the divisor $D$ can be decomposed as $D = D^\circ + D^\infty$ where $D^\circ$ and $D^\infty$ do not have any component in common, and such that
			\begin{itemize}\item[$(i)$]  the domain of definition of $\varphi$ is $Y\setminus D_\infty$, that is, 
				\[\varphi\colon  Y\setminus D^\infty\longrightarrow X\]
				is well-defined and $Y\setminus D^\infty$ is the maximum open set with this property.
				\item[$(ii)$] the pair $(Y\setminus D^\infty, D^\circ\rest{Y\setminus D^\infty})$ is a log-smooth pair for $X$, i.e., $\varphi\rest{Y\setminus D^\infty}$ is a proper morphism from $Y\setminus D^\infty$ to $X$ and the restriction
				\[Y\setminus (D^\circ\cup D^\infty)\longrightarrow X \setminus \varphi(D^\circ)\]
				is an isomorphism.
			\end{itemize}

			A morphism $\overline{\mathbf{Y}}'\rightarrow \overline{\mathbf{Y}}$ between compactified log-smooth pairs $\overline{\mathbf{Y}}'=(Y',D')$ and $\overline{\mathbf{Y}}=(Y,D)$ is a proper morphism $f\colon Y'\rightarrow Y$ which commutes with the structure map $\varphi$ and such that $f^*(D)\subseteq D'$. Notice that in this case we have $f^*(Y\setminus X)=Y'\setminus X$. The category of compactified log smooth pairs will be denoted by $\CLSP_X$.

			\item Given compactified log-smooth pairs $\overline{\mathbf{Y}}'$ and $\overline{\mathbf{Y}}$, we say that $\overline{\mathbf{Y}}'$ dominates $\overline{\mathbf{Y}}$ if there is a morphism $\overline{\mathbf{Y}}'\rightarrow \overline{\mathbf{Y}}$, and we will denote this by $\overline{\mathbf{Y}}'\geq \overline{\mathbf{Y}}$. Similar notations are given for log-smooth pairs and log-smooth compactifications.
		\end{enumerate}
		
	\end{defn}
	
	\begin{prop}\label{prop:filtered} The categories $\CLSP_X$, $\LSP_X$ and $\LSC_X$ are filtered. That is, for any pair of objects $Y_1,Y_2$, there is a third object $Y_3$ together with morphisms $Y_3\rightarrow Y_1$ and $Y_3\rightarrow Y_2$.
	\end{prop}
	\begin{proof} This is obtained by standard arguments using resolution of singularities. 
	\end{proof}
	
	\begin{defn} Given a compactified log-smooth pair $\overline{\mathbf{Y}}=(Y,D)$, we denote by $\M^k({\bf Y}) = \M^k(Y,D)$ the set consisting of all quasi-monomial valuations of rank bounded by $k$ on $Y$ relative to the divisor $D$, denote by $\Sigma(\overline{\mathbf{Y}})=\Sigma(Y,D)$ the dual cone complex to the divisor $D$ on $Y$ and by $\TC^{k-1}\Sigma(\overline{\mathbf{Y}})$ its tangent cone. Similar notations will be used for log-smooth pairs and log-smooth compactifications.
	\end{defn}
	
	\subsection{The retraction map revisited} 
	\begin{prop}\label{prop:retractions} Let $X$ be a smooth variety.
		\begin{enumerate}
			\item For each log-smooth pair $\mathbf{Y}=(Y,D)$ over $X$, there is a continuous retraction
			\[\ret{{\mathbf{Y}}}\colon X^{\bethd,k}\longrightarrow \TC^{k-1}\Sigma(\mathbf{Y}).\]
			\item For each log-smooth compactification $Y$ of $X$, there is a continuous retraction 
			\[\ret{{Y}}\colon X^{\dalethd,k}\longrightarrow \TC^{k-1}\Sigma(Y)\setminus \{0\}.\]
			\item For each compactified log-smooth pair $\overline{\mathbf{Y}}=(Y,D)$ over $X$  there is a continuous retraction
			\[\ret{{\overline{\mathbf{Y}}}}\colon X^{\bir,k}\longrightarrow \TC^{k-1}\Sigma(\overline{\mathbf{Y}}).\]
		\end{enumerate}
	\end{prop}
	\begin{proof} (1) The structure map $\varphi\colon Y\rightarrow X$ is birational, so $K(Y)\cong K(X)$ and we get $Y^{\bir,k}\cong X^{\bir,k}$.  Moreover,  $\varphi\colon Y\rightarrow X$ is a proper map, and so by the valuative criterion of properness, a valuation $\nu$ has a center in $Y$ if and only if it has a center in $X$. Therefore, $Y^{\bethd,k}\cong X^{\bethd,k}$.
		
		The retraction $\ret{{\bf Y}}$ is given by the composition
		\[X^{\bethd,k}\overset{\varphi^*}{\longrightarrow} Y^{\bethd,k}\overset{\ret{D}}{\longrightarrow} \TC^{k-1}\Sigma(\bf Y),\]
		where $\ret{D}$ is given by Definition~\ref{defn:retraction} applied to the pair $(Y,D)$. This proves the claim.

		\noindent
		(2) If $Y$ is a log smooth compactification of $X$, then $D=Y\setminus X$ is an SNC divisor on $Y$. Applying  Definition~\ref{defn:retraction} to this SNC divisor and using that $X^{\bir,k}=Y^{\bir,k}$, we get a map
		\[X^{\bir,k}\longrightarrow \TC^{k-1}\Sigma(Y).\]
		Moreover, a valuation $\nu$ goes to $0\in \TC^{k-1}\Sigma(Y)$ iff $\nu$ is centered outside $D$, hence by restriction we get a map
		\[X^{\dalethd,k}\longrightarrow \TC^{k-1}\Sigma(Y)\setminus\{0\}.\]

		\noindent
		(3) As in (1), the retraction is given by the composition 
		\[X^{\bethd,k}\overset{\varphi^*}{\longrightarrow} Y^{\bethd,k}\overset{\ret{D}}{\longrightarrow} \TC^{k-1}\Sigma(\overline{\bf Y})\]
		where $\ret{D}$ is Definition~\ref{defn:retraction} applied to the divisor $D$ inside $Y$ and $\varphi^*$ is the pullback along the rational map $\varphi\colon Y\rightarrow X$. The claim follows.
	\end{proof}
	
	\begin{prop}\label{prop:compatibility}The retractions from Proposition \ref{prop:retractions} are compatible in the sense that if we have a morphism ${\overline{\bf Y'}}\rightarrow {\overline{\bf Y}}$ of compactified log-smooth pairs, then we have
		\[\ret{{\overline{\bf Y}}}\circ \ret{{\overline{\bf Y}'}}=\ret{{\overline{\bf Y}}}.\]
		Similar statements hold for log-smooth pairs and log-smooth compactifications.
	\end{prop}
	
	\begin{proof} Let $\nu \in Y^{\bir, k}$ be a valuation. Consider a compactified log-smooth pair ${\overline{\bf Y}'} \geq {\overline{\bf Y}}$ above $\overline{\bf Y}$ and denote by $D'_1, \dots, D'_l$ all the components of $D'$ in $Y'$.   Let $D_i$ be a component of $D$ in $Y$. There exists a subset $J_i \subseteq [l]$ such that
		$\pi^*(D_i)=\sum_{j\in J_i} n_jD'_j$.
		
		Let $h_j$ be local parameters for $D'_j$ around the center $\centrebis{\nu}$ of $\nu$ in $Y'$. The product $\prod_{j\in J_i} h_j$ is a local equation for $D_i$ around the center $\centre{\nu}$ of $\nu$ in $Y$. We have 
		\[\nu(D_i)=\nu(\prod_{j\in J_i} h^{n_j}_j)=\sum_{j\in J_i} n_j\nu(h_j)=\sum_{j\in J_i} n_j\ret{\overline{\bf Y}'}(\nu)(h_j)=\ret{\overline{\bf Y}'}(\nu)(\prod h^{n_j}_j)=\ret{\overline{\bf Y}'}(\nu)(D_i).\]
		This implies that the two valuations $\nu$ and $\ret{\overline{\bf Y}'}(\nu)$ are mapped to the same point by the retraction map $\ret{\overline{\bf Y}}$, that is,  $\ret{\overline{\bf Y}}(\nu)=\ret{\overline{\bf Y}}(\ret{\overline{\bf Y}'}(\nu))$. Since this holds for all valuations $\nu$, the compatibility of the retraction maps $\ret{\overline{\bf Y}}$ and $\ret{\overline{\bf Y}'}$ follows. 
	\end{proof}
	
	\subsection{The retraction inequality}
	
	We finish this section by recalling the following useful statement from~\cite{JM12} called \emph{the retraction inequality}. This will be used in the next section. 
	\begin{prop}[Retraction inequality]\label{prop:val-inequality} Let ${\bf Y}=(Y,D)\in \CLSP_X$ be a compactified log-smooth pair and let $\nu\in X^{\bir,k}$ be  a valuation with center a point $x$ of $Y$. Then, for each $f\in \mathcal{O}_{Y,x}$ we have the inequality 
		\begin{equation}\nu(f)\succeq(\ret{{\bf Y}}(\nu))(f)
		\end{equation}
		with equality when the zero set $V(f)$ of $f$ is included in $D$ locally around $x$.	
	\end{prop}
	\begin{proof} Denote by $D_1, \dots, D_m$ the components of $D$ which pass through $x$ and take local equations $z_i$ for each component $D_i$ around $x$. The family $\{z_i\}_{i=1}^m$ can be extended to a set of local parameters $\{z_i\}_{i=1}^r$ for $Y$ at $x$. By Corollary \ref{cor:finite} there is a finite admissible expansion 		\[f=\sum_{\beta \in \minS{f}} a_\beta u_\beta z^\beta.\]
	We get 
		\begin{equation*}
		\nu(f)
		\geqlex\; \min_{\beta \in \minS{f}}\left\{\nu(a_\beta u_\beta z^\beta)\right\}
		\geqlex\; \min_{\beta \in \minS{f}} \{\nu(z^\beta)\}
		=\; (\ret{{\bf Y}}(\nu))(f),
		\end{equation*}
		which is the stated inequality.  Suppose now that $V(f)\subseteq D$ locally around $x$. We can write $f$ in $\mathcal{O}_{Y,x}$ as $f=u\prod_{i=1} ^m z_i^{n_i}$ for a unit $u\in \mathcal O_{Y, x}$ and non-negative integers $n_i$.
		We conclude
		\[\nu(f)=\sum_{i=1}^m n_i \nu(z_i)=\sum_i n_i \ret{{\bf Y}}(\nu)(z_i)=(\ret{{\bf Y}}(\nu))(f).\]
	\end{proof}
	
	\section{Limit formulae}\label{sec:limit}
	Let $X$ be a smooth variety over an algebraically closed field $k$. In this section we will see how it is possible to reconstruct the space $X^{\bethd,k}$ of valuations with center inside $X$ in terms of the spaces $\TC^{k-1} \Sigma(\mathbf{Y})$ for log-smooth pairs $\mathbf{Y}$ studied in the previous section (Theorem \ref{thm:bir limit} below). We will give a similar result for the set $X^{\dalethd,k}$ of valuations whose center is outside $X$ in terms of a limit $\TC^{k-1} \Sigma(Y)$ over log-smooth compactifications of $X$ (Theorem \ref{thm:center limit} below). Similarly, we show how to reconstruct the centroidal filtration $\mathscr{F}^r X^{\bir,k}$ in terms of a centroidal filtration for $\TC^{k-1}\Sigma(\overline{\bf{Y}})$ over compactified log-smooth pairs $\overline{\bf Y}$.
	
	\subsection{Limit formula for $X^{\bethd,k}$}
	
	The compatibility shown in \ref{prop:compatibility} for the retraction maps presented in Proposition \ref{prop:retractions} implies that there exist natural continuous maps 
	\begin{align}
	 r\colon X^{\bethd,k}& \longrightarrow \lim_{\substack{\longleftarrow\\ {\bf Y}\in \LSP_X}} \TC^{k-1}\Sigma({\bf Y})
	\label{map:limit1} \\
	r\colon X^{\dalethd,k}&\longrightarrow \lim_{\substack{\longleftarrow\\Y\in \LSC_X}} \TC^{k-1}\Sigma(Y)\setminus\{0\}
	\label{map:limit2}\\ 
	r\colon X^{\mathrm{bir},k}&\longrightarrow \lim_{\substack{\longleftarrow\\ \overline{\bf Y}\in \CLSP_X}} \TC^{k-1}\Sigma(\overline{{\bf Y}}).
	\label{map:limit3}
	\end{align}
	
	The objective of this section is to prove the following theorem.
	
	\begin{thm}[Limit formula]\label{thm:bir limit} The maps \refeq{map:limit1}, \refeq{map:limit2} and \refeq{map:limit3} above are all homeomorphisms. 
	\end{thm}

	It will be enough to construct an inverse for each of these maps. We treat the case of \refeq{map:limit3}, the proofs in the other cases will be essentially the same. The inverse for this map is 
	
	\begin{align}\label{map:inverse} 
	\begin{split} q\colon \lim_{\substack{\longleftarrow\\ \clsp{Y}\in \CLSP_X}} \TC^{k-1}\Sigma(\clsp{Y}) &\longrightarrow X^{\mathrm{bir},k} \\ s=[(x;\underline{w})]_{\clsp{Y}}&\longmapsto \nu_s
	\end{split}
	\end{align}
	where $\nu_s$ is the valuation defined by
	\[\nu_s(f)=\sup_{\clsp{Y}}\nu_{s, _{\clsp{Y}}}(f), \qquad  \text{ for every $f\in \bigcap_{\clsp{Y}} \sO_{X,\eta_{s, \clsp{Y}}}$}.\]	
	Here,  $\nu_{s, _{\clsp{Y}}}\coloneqq \nu_{x,\underline{w}}$ for the element $(x;\underline{w})$ in the position indexed by $\clsp{Y}$ in the sequence $s$, and $\eta_{s, _{\clsp{Y}}}$ is the center of $\nu_{s, _{\clsp{Y}}}$.
	
	\begin{prop} The map \refeq{map:inverse} is well defined and is the inverse for the map \refeq{map:limit3}.
	\end{prop}
	
	\begin{proof} We first note that if $\clsp{Y}'\geq \clsp{Y}$, then $r_{D}(\nu_{s,\clsp{Y}'})=\nu_{s,\clsp{Y}}$ and therefore $\eta_{s,\clsp{Y}'}\in\overline{\{\eta_{s,\clsp{Y}'}\}}$. Hence, the sequence of points $\eta_{s,\clsp{Y}}$ is decreasing for the order given by specialization. Therefore, it eventually becomes constant equal to some $\eta$. We get the equality $\bigcap_{\overline{\bf Y}} \sO_{X,\eta_{s, \bf Y}}=\sO_{X,\eta}$.

		By Proposition \ref{prop:val-inequality} the sequence 
		\[[\nu_{s,\clsp{Y}}(f)]_{\clsp{Y}}\]
		is increasing. Moreover, if we fix $\clsp{Y}=(Y,D)$, a compactified log smooth pair, and we take ${\clsp{Y}'}=(Y',D'_\red)$ where $(Y',D')$ is an embedded resolution of singularities for $V(f)\cup D^\infty\subseteq Y$, we get that
		\[V(f)\cup\Supp(D^\infty)\subseteq \Supp(D_{\red}).\]
		
	For any $\clsp{Y}''\geq \clsp{Y}'$, we get $D''\supseteq V(f)$. By the equality part of Proposition \ref{prop:val-inequality}, we infer  \[\nu_{s,\clsp{Y}''}(f)=\nu_{s,\clsp{Y}'}(f).\]
	This means the sequence eventually becomes constant, and so, for $f\in \sO_{X,\eta}$, the supremum is attained. Moreover, Proposition \ref{prop:filtered} shows that given $f,g\in \sO_{X,\eta}$, there is a compactified log-smooth pair ${\clsp{Y}}$ in which the sequences $[\nu_{s,{\clsp{Y}}}(f)]_{\clsp{Y}}$, $[\nu_{s,\clsp{Y}}(g)]_{\clsp{Y}}$, $[\nu_{s,\clsp{Y}}(f+g)]_{\clsp{Y}}$, $[\nu_{s, \clsp{Y}}(fg)]_{\clsp{Y}}$ are all constant at the same time from ${\clsp{Y}}$ onwards. Hence, for this $\clsp{Y}$, we have
		\begin{align*} \nu_s(fg)=\nu_{s,\clsp{Y}}(fg)=\nu_{s,\clsp{Y}}(f)+\nu_{s, \clsp{Y}}(g)=\nu_s(f)+\nu_s(g) .
		\end{align*}
Similarly, $\nu_s(f+g)\geq \min\{\nu_s(f),\nu_s(g)\}$, so $\nu_s$ is a valuation. This shows that $q$ is well defined.

		We now show that $q$ is the inverse of $r$. For this, we need to check that the composition over each side gives the identity. This translates into the equalities
		\begin{enumerate}
			\item $[r_{\clsp{Y}}(\nu_s)]_{\clsp{Y}}=s$ for any  $\displaystyle s=[(x;w)]_{\clsp{Y}}\in \lim_{\substack{\longleftarrow\\ {\clsp{Y}}\in \overline{\CLSP}_X}} \TC^{k-1}\Sigma(\clsp{Y})$, and
			\item $\nu_s=\nu$ for any $\nu\in X^{\bir,k}$ with $s=[r_{\clsp{Y}}(\nu)]_{\clsp{Y}}.$
		\end{enumerate}
		
		For the first equality, we have to prove that for any ${\clsp{Y}}\in \CLSP_X$ we have $r_{\clsp{Y}}(\nu_s)=(x;w)$ for $(x;w)$ in the $\clsp{Y}$ instance of the sequence $s$. In order to do this consider $\{z_i\}_i$ local equations for the components $D_i$ of the divisor $D$ defining ${\clsp{Y}}$ around the center of $\nu_s$ in $Y$. There is a compactified log-smooth pair ${\clsp{Y}}'=(Y',D')$ in which we simultaneously have the equalities $\nu_s(z_i)=\nu_{s,\clsp{Y}'}(z_i)$ for each $i$. In this case, we get $r_{\clsp{Y}}(\nu_s)=r_{\clsp{Y}}(\nu_{s, \clsp{Y}'})$ and by the compatibility of the retraction maps, we infer that $r_{\clsp{Y}}(\nu_{s, \clsp{Y}'})=(x;w)$.

		For the second equality, it is enough to prove that for each $f\in \sO_{X,\eta}$, we have $\nu_s(f)=\nu(f)$. This follows directly by the equality part in Proposition \ref{prop:val-inequality}.
	\end{proof}
	
	\begin{prop} The map  \refeq{map:inverse} is continuous.
	\end{prop}
	\begin{proof} The topology of $X^{\bir,k}$ is generated by open sets of the form
		\[U=\{\nu \in X^{\bir,k}\mid \nu(f)\in A\}\]
		for $A\subseteq \mathbb{R}^k$ an Euclidean open set and $f\in K(X)^*$ a rational function. Hence, given a fixed sequence  $s=[(x;w)]_{\clsp{Y}}$ such that $q(s)\in U$, it is enough to find a neighborhood $V$ of $s$ such that $q(V)\subseteq U$. For this, take a compactified log-smooth pair ${\clsp{Y}}$ and consider $f=\frac{g}{h}$ where $g,h\in \mathcal{O}_{Y,c_Y(q(s))}$. Take a compactified log-smooth pair $\clsp{Y}'$ by choosing an embedded resolution of $V_Y(g)\cup V_Y(h)\cup D$. Consider then 
		\[V_g=\biggl\{t\in \lim_{\substack{\longleftarrow\\ {\clsp{Y}}\in \CLSP_X}} \TC^{k-1}\Sigma({\clsp{Y}})\,\, \Bigl|\,\, \nu_{t,{\clsp{Y}}} \text{ has center inside } V_{Y'}(g)\biggr\}, \textrm{ and }\]
		\[V_h=\biggl\{t\in \lim_{\substack{\longleftarrow\\ {\clsp{Y}}\in \CLSP_X}} \TC^{k-1}\Sigma({\clsp{Y}})\,\, \Bigl|\,\, \nu_{t,{\clsp{Y}}} \text{ has center inside } V_{Y'}(h)\biggr\}.\]
		
		By Proposition \ref{prop:anticontinuous}, for the center map $\TC^{k-1}\Sigma(\clsp{Y})\rightarrow Y'$, we see that both $V_g$ and $V_h$ are open neighborhoods of $s$ in the direct limit. By Proposition \ref{prop:val-inequality}, for each $t\in V_g\cap V_h$ we have $\nu_t(g)=\nu_{t,\clsp{Y}}(g)$ and $\nu_t(h)=\nu_{t,\clsp{Y}}(h)$. We thus get $\nu_t(f)=\nu_{t,\clsp{Y}}(f)$. Consider \[V'=\biggl\{t\in \lim_{\substack{\longleftarrow\\ {\clsp{Y}}\in \CLSP_X}}\TC^{k-1}\Sigma({\clsp{Y}})\mid \nu_{t,\clsp{Y}}(f)\in A\biggr\}.\]
		This is another open neighborhood of $s$ in the direct limit. To conclude, note that for $V=V_g\cap V_h \cap V'$ and $t\in V$, we have
		\[\nu_{t,\bf Y}(f)=\nu_{t}(f)\in A\]
		which shows that $q(V)\subseteq U$. This proves the continuity.
	\end{proof}
As a consequence of the limit formulae and weak factorization theorem, we get the following.
\begin{cor}[Density of flag valuations] For $1\leq k\leq \dim X$, the valuations which are equivalent to a flag valuation of rank $k$ in some birational model of $X$ are dense in $X^{\mathrm{bir},k}$.
\end{cor}
\begin{proof} The isomorphism \eqref{map:limit3} shows that the family of open sets of the form $r_{\overline{\bf Y}}^{-1}(U)$ where  $\overline{\bf Y}$ moves through $\CLSP$, and $U\subseteq \TC^{k-1}\Sigma({\bf \overline{Y}})$ is an open set in the tropical topology, form a basis of $X^{\bir,k}$. Hence, as $r(U)\subseteq r_{\overline{\bf Y}}^{-1}(U)$, it is enough to show that $r(U)$ contains a valuation which is a flag valuation in some birational model of $X$. For this, take $(x;w_1,\dots,w_k)\in U$ with $x,w_1,\dots,w_k$ rational and linearly independent. We will show that $\nu_{x,\underline{w}}\in r(U)$ is equivalent to a flag valuation. For this, take a rational subdivision $\widetilde{\Sigma}$ of $\Sigma$ containing $x,w_1,\dots,w_k$ in different rays of a single cone $\sigma$. By weak factorization theorem \cite{Wlo97}, there is a way to obtain $\widetilde{\Sigma}$ as a sequence of blow-ups and blow-downs on the fans $\Sigma$. The same sequence of blow-ups and blow-down gives rise to a birational model $X'$ and a SNC divisor $D'\subseteq X'$ in which $\nu_{x;\underline{w}}$ is a diagonal scalar multiple of a flag valuation by \ref{thm:flagvaluations}.
\end{proof}
There are similar corollaries for $X^{\bethd,k}$ and $X^{\dalethd,k}$.

	\subsection{Refined limit formula}
	\label{subsection:refined}
	We define the centroidal filtration on tangent cones. 
	\begin{defn}[Centroidal filtration] Given a compactified log-smooth pair $\clsp{Y}=(Y,D)$ over $X$, we have a decomposition of $D$ as $D=D^\circ\cup D^\infty$. This gives the subcomplex $\Sigma(D^\circ)$ inside $\Sigma(\clsp{Y})$ which we denote by $\Sigma(\clsp{Y}^\circ)$.

		We define the \emph{centroidal filtration} of $\TC^{k-1}{\clsp{Y}}$ to be the filtration
		\[\mathscr{F}^0\TC^{k-1}\Sigma(\clsp{Y})\supseteq \mathscr{F}^0\TC^{k-1}\Sigma(\clsp{Y})\supseteq \dots \supseteq \mathscr{F}^k\TC^{k-1}\Sigma(\clsp{Y})\]
		given for $0\leq r \leq k$ by
		\begin{equation*} \mathscr{F}^r\TC^{k-1}\Sigma(\clsp{Y})\coloneqq\Bigl\{\,(x;\underline{w}_{k-1})\in \TC^{k-1}\Sigma(\clsp{Y})\,\,\Bigl|\,\, (x;\underline{w}_{r-1})\in \TC^{r-1}\Sigma(\clsp{Y}^\circ)\, \Bigr\}.
		\end{equation*}
	\end{defn}
	
	\begin{remark} For $i<j$, property $(x;\underline{w}_{j})\in \TC^{j}\Sigma(\clsp{Y}^\circ)$ implies $(x;\underline{w}_{i})\in \TC^{i}\Sigma(\clsp{Y}^\circ)$. Therefore, the sequence $\bigl(\mathscr{F}^r\TC^{k-1}\Sigma(\clsp{Y})\bigr)_r$ is indeed decreasing. Moreover, we have \[\mathscr{F}^0\TC^{k-1}\Sigma(\clsp{Y})=\TC^{k-1}\Sigma(\clsp{Y}) \text{ and } \mathscr{F}^k\TC^{k-1}\Sigma(\clsp{Y})=\TC^{k-1}\Sigma(\clsp{Y}^\circ).\]
	This is similar to the centroidal filtration on $X^{\bir,k}$.
	\end{remark}
	
	The limit formula for compactified log-smooth pairs in the previous subsection preserves the centroidal filtrations, and we obtain a limit description of each term of the filtration.
	
	\begin{thm}\label{thm:center limit} For each $0\leq r \leq k$, the isomorphism of Theorem \ref{thm:bir limit} restricts to a homeomorphism
		\begin{equation*} 
		\mathscr{F}^rX^{\mathrm{bir},k}\longrightarrow \lim_{\substack{\longleftarrow\\ \clsp{Y}\in \CLSP_X}} \mathscr{F}^r\TC^{k-1}\Sigma(\clsp{Y}).
		\end{equation*}
	\end{thm}
	
	\begin{proof}
		We first note that given $(x;\underline{w})\in \mathscr{F}^r\TC^{k-1}\Sigma(\clsp{Y})$, the center of the valuation $\mathrm{proj}_r(\nu_{x,\underline{w}})$ is the point $\eta_\sigma$ where $\sigma$ is the smallest face such that $(x,\underline{w}_r)\in\TC^{r}\sigma$. Hence, if $(x;\underline{w})\in \mathscr{F}^r\TC^{k-1}\Sigma(\clsp{Y}^\circ)$, then the center of $\nu_{x,\underline{w}}$ is inside $D^\circ\subseteq X$. This proves that the inclusion $\TC^{k-1}\Sigma(\clsp{Y})\hookrightarrow X^{\bir,k}$ restricts to an inclusion $\mathscr{F}^r\TC^{k-1}\Sigma(\clsp{Y})\hookrightarrow \mathscr{F}^rX^{\bir,k}$
		
		Moreover, since for any compactified log-smooth pair, the center of $\nu$ is a specialization of the center of $r_{\clsp{Y}}(\nu)$, we see that for each pair $\clsp{Y}'\geq \clsp{Y}$, the retraction map in \ref{prop:retractions} induces a map
		\[\mathscr{F}^r\TC^{k-1}\Sigma(\clsp{Y}')\longrightarrow \mathscr{F}^r\TC^{k-1}\Sigma(\clsp{Y})\]
		and these maps are still compatible. 
		This implies that once we take the inverse limit, we obtain a natural map
		\[r\colon \mathscr{F}^rX^{\bir,k}\longrightarrow \lim_{\substack{\longleftarrow\\ \mathbf{Y}\in \CLSP_X}} \mathscr{F}^r\TC^{k-1}\Sigma(\clsp{Y}).\]
		
		On the other hand, the inverse map $q$ defined in \ref{map:inverse} also restricts to a map
		\[q\colon \mathscr{F}^r\TC^{k-1}\Sigma(\clsp{Y})\longrightarrow \mathscr{F}^rX^{\bir,k}.\]
		Indeed, if $s=[(x;w)]_{\bY}$ is a sequence of elements in $\mathscr{F}^r\TC^{k-1}\Sigma(\clsp{Y})$, the center of the elements $\mathrm{proj}_r(\nu_{s,\clsp{Y}})$ in $X$ are points $\eta_{s,\clsp{Y}}$ with the property that $\eta_{s,\clsp{Y}'}$ specializes $\eta_{\clsp{Y}}$ if $\clsp{Y}'$ dominates $\clsp{Y}$. $X$ being Noetherian, the sequence $[(\eta_{s,\clsp{Y}})]_{\clsp{Y}}$ is eventually constant, and hence $\eta_{s,\clsp{Y}}$ is in the projection of a stratum of $D^\circ$ onto $X$. The center of $\mathrm{proj}_r(q(s))$ is  thus on $X$, and so $q(s)\in \mathscr{F}^rX^{\bir,k}$. The maps $r$ and $q$ are still inverse to each other. This finishes the proof.
	\end{proof}
	
	\begin{cor} The limit above can be restricted to each stratum in the centroidal filtration, that is, for each $r$, we have a homeomorphism
		\begin{equation*} 
		\mathscr{F}^rX^{\mathrm{bir},k}\setminus  \mathscr{F}^{r+1}X^{\mathrm{bir},k} \;\longrightarrow \!\lim_{\substack{\longleftarrow\\ \clsp{Y}\in \CLSP_X}} \! \mathscr{F}^r\TC^{k-1}\Sigma(\clsp{Y}) \setminus \mathscr{F}^{r+1}\TC^{k-1}\Sigma(\clsp{Y}).
		\end{equation*}
	\end{cor}

	\section{Variations of Newton-Okounkov bodies} \label{sec:heuristic}

Let $X$ be a smooth projective variety of dimension $d$ over an algebraically closed field $\kappa$, with function field $K(X)$. 
Consider a big line bundle $L = \mathcal O(E)$ on $X$ with the corresponding graded algebra
\[H_\bullet=\bigoplus_{n\geq 0} H_n\]
where each $H_n = H^0(X, \mathcal O(nE))$ is a $\kappa$-vector subspace of $K(X)$ of finite dimension.

Consider the space $\mathrm{BC}(\mathbb{R}^d)$ of compact convex subsets of $\mathbb{R}^d$ endowed with the Hausdorff distance. We consider the map
\begin{align}\begin{split}\label{map:okounkov}
\Delta\colon X^{\bir,d}&\longrightarrow \mathrm{BC}(\mathbb{R}^d) \\
\nu &\longmapsto \Delta_\nu=\overline{\bigcup_{n\geq 0} \left \{\frac{\nu(f)}{n}\mid f\in H_n \right \}}
\end{split}\end{align}
which attaches to each valuation the corresponding Newton-Okounkov body.

\begin{conj} The restriction of the map $\Delta$ to each higher rank skeleton is continuous.
\end{conj}

We provide a heuristic argument for the validity of this conjecture. First, we show 

\begin{thm}\label{thm:weak convergence} Let $(C_n)_n\in \mathrm{BC}(\mathbb{R}^d)$ be a sequence of full dimensional compact convex subsets of $\mathbb{R}^d$. Then we have that \[C_n\overset{n\rightarrow \infty}{\longrightarrow} C\in \mathrm{BC}(\mathbb{R}^d)\] in the Hausdorff distance if and only if for each continuous function $f\colon \mathbb{R}^d\longrightarrow \R$ with compact support, we have
	\[\int_{C_n}f(x)dx\overset{n\rightarrow \infty}{\longrightarrow} \int_C f(x)dx.\]
\end{thm}
\begin{proof}
	
	\noindent
	{$(\Rightarrow)$} If we denote by $B(C;\varepsilon)$ the $\varepsilon$-neighbourhood of $C$, we have the inequality
	\[\mathrm{Vol}(C_n\Delta C)\leq \mathrm{Vol}(B(C_n;d_H(C_n,C))\setminus C_n)+\mathrm{Vol}(B(C;d_H(C_n,C))\setminus C).\]
	Here $C_n\Delta C$ is the symmetric difference of $C_n$ and $C$. The left hand side goes to $0$ as $n$ goes to infinity. We thus get $\mathrm{Vol}(C_n\Delta C)\longrightarrow 0$. Hence we have the almost everywhere convergence $f\cdot\mathbbm{1}_{C_n}\rightarrow f\cdot\mathbbm{1}_{C}$, and so by the dominated convergence theorem the integrals converge.
	
	{$(\Leftarrow)$} 
	If $C_n$ does not converge to $C$, passing to a subsequence if necessary, we get the existence of $\varepsilon>0$ such that the Hausdorff distance $d_H(C_n,C)\geq\varepsilon$ for all $n$. We thus have either \[\sup_{x\in C_n}\{d(x,C)\}\geq\varepsilon \qquad \text{ or }\qquad \sup_{y\in C}\{d(C_n,y)\}\geq\varepsilon\] happen infinitely many times.

	In the first case, since $C_n$ is compact for each $n$, there is an $x_n\in C_n$ for which the supremum is attained. Let $x'$ be an accumulation point of $(x_n)$.  Consider an open ball of the form $B(y;\varepsilon)$ contained in $C$. If for each continuous function $f$ the integrals above converge, by taking $f$ as a bump function supported exactly on $B(y;\varepsilon)$, we get that $\mathrm{Vol}(B(y;\varepsilon)\setminus C_n)\rightarrow 0$. Therefore, since $C_n$ and $B(y, \varepsilon)$ are convex, for $n$ big enough, we would get that $C_n$ contains the ball $B(y;\varepsilon/2)$. This implies that for each such $n$, we have
	\[C_n\setminus C\supseteq \operatorname{conv}(B(y;\varepsilon/2)\cup\{x\})\setminus C.\]
	The set appearing in the right hand side is independent of $n$, and has nonempty interior. Taking again  a bump function supported in this set, we see that the integrals could not converge which would give a contradiction.

	In the second case, for each $n$, there is a $y_n\in C$ for which the supremum is attained. By compactness, and passing to a subsequence if necessary, we can assume that $y_n$ converge to a point in $C$ that we denote by $y'$. For $n$ big enough, we have $d(C_n,y')\geq \varepsilon/2$, as $d(C_n,y')\geq |d(C_n,y_n)-d(y_n,y')|$. It follows that $B(y';\varepsilon/2)\cap C$ is disjoint from $C_n$ for infinitely many $n$. By taking $f$ as a bump function supported on $B(y';\varepsilon/2)\cap C$ we see that $C_n$ would not converge to $C$ in a weak sense, which would be a contradiction. 
\end{proof}

\subsubsection*{Heuristic argument for the validity of the conjecture} We will use the well-known fact that the family of sets  \[\Delta^n_\nu=\left\{\frac{\nu(f)}{n}\mid f\in H_n\right\}\]
is equidistributed in the set $\Delta_\nu$ endowed with Lebesgue measure~\cite{Bouck12}. That is, for each continuous function $h$ of compact support, we have 

\[\lim_{n\rightarrow \infty}\frac 1{N_n}\sum_{x\in \Delta^n_\nu}h(x)=\int_{\Delta_\nu}h(x) dx\]
with $N_n$ the dimension of $H_n$ over $\kappa$ (equal to the size of $\Delta_\nu^n$).

If $\nu_j$ is a sequence of valuations converging to $\nu$, we would use Theorem \ref{thm:weak convergence} to get $\Delta_{\nu_k}\rightarrow \Delta_\nu$ in the Hausdorff distance if we could interchange the limits in the following 
\begin{align*}
\lim_{j\rightarrow \infty}\int_{\Delta_{\nu_k}} h(x)dx&=\lim_{j\rightarrow \infty}\lim_{n\rightarrow \infty}\frac{1}{N_n}\sum_{x\in \Delta^n_{\nu_j}}h(x)\overset{?}{=}\lim_{n\rightarrow \infty}\lim_{j\rightarrow \infty}\frac{1}{N_n}\sum_{x\in \Delta^n_{\nu_j}}h(x)\\
&=\lim_{n\rightarrow \infty}\frac{1}{N_n}\sum_{x\in \Delta^n_{\nu}}h(x)=\int_{\Delta_\nu}h(x) dx.
\end{align*}

		\bibliographystyle{alpha}
	\bibliography{bibliographie}
	
	\end{document}

%% file: dualconecomplex.pdf_t
\begin{picture}(0,0)%
\includegraphics{dualconecomplex.pdf}%
\end{picture}%
\setlength{\unitlength}{3947sp}%
\begingroup\makeatletter\ifx\SetFigFont\undefined%
\gdef\SetFigFont#1#2#3#4#5{%
  \reset@font\fontsize{#1}{#2pt}%
  \fontfamily{#3}\fontseries{#4}\fontshape{#5}%
  \selectfont}%
\fi\endgroup%
\begin{picture}(8324,7431)(2911,-7470)
\put(10501,-4861){\makebox(0,0)[lb]{\smash{{\SetFigFont{25}{30.0}{\familydefault}{\mddefault}{\updefault}{\color[rgb]{0,0,0}$\zeta$}%
}}}}
\put(3226,-5161){\makebox(0,0)[lb]{\smash{{\SetFigFont{25}{30.0}{\familydefault}{\mddefault}{\updefault}{\color[rgb]{0,0,0}$\sigma$}%
}}}}
\put(7276,-5311){\makebox(0,0)[lb]{\smash{{\SetFigFont{25}{30.0}{\familydefault}{\mddefault}{\updefault}{\color[rgb]{0,0,0}$\eta$}%
}}}}
\put(2926,-1111){\makebox(0,0)[lb]{\smash{{\SetFigFont{41}{49.2}{\familydefault}{\mddefault}{\updefault}{\color[rgb]{0,0,0}$\Sigma$}%
}}}}
\put(5626,-5686){\makebox(0,0)[lb]{\smash{{\SetFigFont{25}{30.0}{\familydefault}{\mddefault}{\updefault}{\color[rgb]{0,0,0}$\tau$}%
}}}}
\end{picture}%